\documentclass[11pt]{article}
\usepackage{amssymb,amsmath,amsthm,hyperref}
\numberwithin{equation}{section} \allowdisplaybreaks
\newtheorem{thm}{Theorem}[section]
\newtheorem{lem}{Lemma}[section]

\theoremstyle{definition}

\theoremstyle{remark}
\newtheorem{rem}{Remark}[section]

\begin{document}

\title{Existence and uniqueness results for viscous, heat-conducting  3-D
fluid with vacuum\thanks{This work is supported by NSFC 10571158,
and Zhejiang Association for international exchange of personal.}}

\author{Ting Zhang\thanks{Email: TZ: zhangting79@hotmail.com; DF:
dyf@zju.edu.cn}  and  Daoyuan Fang\\
\textit{\small Department of Mathematics,  Zhejiang University,
Hangzhou 310027, PR China}}

\date{}

\maketitle
\begin{abstract}
We consider the 3-D full Navier-Stokes equations whose the
viscosity coefficients and the thermal conductivity coefficient
depend on the density and the temperature. We prove the local
existence and uniqueness of the strong solution in a domain
$\Omega\subset\mathbb{R}^3$. The initial density may vanish in an
open set and $\Omega$ could be  a bounded or unbounded domain. We
also prove a blow-up criterion for the solution. Finally, we show
the blow-up of the smooth solution to the compressible
Navier-Stokes equations in $\mathbb{R}^n$ ($n\geq1$) when the
initial density has compactly support and the initial
total momentum is nonzero. \\
\textbf{Keywords:} {Compressible Navier-Stokes equations;
Existence; Uniqueness; Blow-up.} \\ \textbf{AMS Subject
Classification:} 35Q30; 76N10.
\end{abstract}

%\tableofcontents

\section{Introduction}\label{2-b-Sec1}
The motion of a viscous, heat-conducting fluid in a domain
$\Omega\subset \mathbb{R}^3$ can be described by the system of
equations, known as the Navier-Stokes equations:
    \begin{equation}
      \rho_t+\mathrm{div } (\rho u)=0,
      \label{2-b-E1.1}
    \end{equation}
             \begin{equation}
              (\rho u)_t+\mathrm{div} (\rho u\otimes u)
              -\mathrm{div}(\mu(\nabla u+\nabla u^\top))-\nabla(\lambda\mathrm{div} u)
              +\nabla p=\rho f
              \label{2-b-E1.3}
            \end{equation}
    \begin{equation}
          (\rho e)_t+\mathrm{div }(\rho e u)- \mathrm{div}(\kappa\nabla \theta)
          +p\mathrm{div }u=\frac{\mu}{2}|\nabla u+\nabla
          u^\top|^2+\lambda(\mathrm{div } u)^2+\rho h,
          \label{2-b-E1.2}
    \end{equation}
in $(0,+\infty)\times \Omega\subset
(0,+\infty)\times\mathbb{R}^3$, and the initial and boundary
conditions:
    \begin{equation}
          (\rho,\theta,u)|_{t=0}=(\rho_0,\theta_0,u_0)
          \textrm{ in }\Omega,
          \label{2-b-E1.4-1}
        \end{equation}
            \begin{equation}
              (\theta,u)=(\theta_b,0)
              \textrm{ on }(0,\infty)\times\partial\Omega,
            \end{equation}
%when $\Omega$ is   a bounded domain in $\mathbb{R}^3$ with smooth
%boundary, the half space $\mathbb{R}^3_+$ or an exterior domain in
%$\mathbb{R}^3$ with smooth boundary.
                \begin{equation}
                  (\rho, \theta, u)(t,x)
                  \rightarrow(\rho^\infty,\theta^\infty,0)\textrm{ as }
                  |x|\rightarrow\infty.
                  \label{2-b-E1.4-3}
                \end{equation}
%when $\Omega$ is  the whole space $\mathbb{R}^3$, the half space
%$\mathbb{R}^3_+$ (need $\theta_b=\theta^\infty$) or an exterior
%domain in $\mathbb{R}^3$ with smooth boundary.
In this paper, $\Omega$ is either a bounded domain in
$\mathbb{R}^3$ with smooth boundary or an usual unbounded domain
such as the whole space $\mathbb{R}^3$, the half space
$\mathbb{R}^3_+$ ($\theta_b=\theta^\infty$) and an exterior domain
in $\mathbb{R}^3$ with smooth boundary.  The known fields $h$ and
$f$ denote a heat source and an external force per unit mass.
$\rho$, $\theta$, $e$, $p$ and $u$ denote the unknown density,
temperature, internal energy, pressure and velocity fields of the
fluid, respectively. The pressure $p$ and internal energy $e$ are
related to the density and temperature via the equations of state:
        $$
        p=p(\rho,\theta),\ e=e(\rho,\theta).
        $$
In accordance with the basic principles of classical
thermodynamics,  $p$ and $e$ are interrelated through the
following relationship:
    $$
    \partial_\rho e=\frac{1}{\rho^2}(p-\theta \partial_\theta p).
    $$
In view of classical thermodynamics (see the book \cite{Lions96}),
the viscosity coefficient $\mu=\mu(\rho,\theta)$, the second
viscosity coefficient   $\lambda=\lambda(\rho,\theta)$ and the
thermal conductivity coefficient $\kappa=\kappa(\rho,\theta)$ are
required to satisfy the natural restriction
    $$
    \mu>0,\ \ 3\lambda+2\mu\geq  0
    \ \ \textrm{ and }\ \kappa\geq0.
    $$
In this paper, we assume that $(e,p,\kappa,\mu,\lambda)$ satisfy
    \begin{equation}
   \partial_\theta e(\rho,\theta)\geq \underline{c}>0,
\ 0\leq p(\rho,\theta)\leq\rho^\frac{1}{2} p_0(\rho)(1+\theta^2),
    \label{2-b-E1.7}
    \end{equation}
    \begin{equation}
  |\partial_\theta p(\rho,\theta)|\leq \rho^\frac{1}{2} p_1(\rho)(1+\theta),
    \ |\partial_\rho p(\rho,\theta)|\leq  p_2(\rho)(1+\theta^2),
    \label{2-b-E1.8}
\end{equation}
    \begin{equation}
    \mu(\rho,\theta)\geq\mu^0>0,\ \ 3\lambda(\rho,\theta)+2\mu(\rho,\theta)\geq  0,
    \ \ \kappa(\rho,\theta)\geq \kappa^0>0,
    \label{2-b-E1.4}
    \end{equation}
        \begin{equation}
          \rho^{-\frac{1}{4}}\partial_{\theta}(\partial_\theta e,\partial_\rho e,\mu,\lambda,\kappa)(\rho,\theta)
        \in L^\infty_{loc}([0,\infty)
        \times[0,\infty)),\label{2-b-E1.4.1}
        \end{equation}
    \begin{equation}
      \partial_\theta e,\partial_\rho e,p,\mu,\lambda,\kappa\in
      C^1([0,\infty)\times[0,\infty)),
      \ p_0,p_1,p_2\in C([0,\infty)),
      \label{2-b-E1.5}
    \end{equation}
and
    \begin{equation}
      \partial_\theta p,\partial_\rho p\in
      C^\alpha([0,\infty),L^\infty_{loc}([0,\infty))),
    \end{equation}
for all $\rho,\theta\geq0$ and a constant
$\alpha\in(\frac{1}{2},1)$. It is easy to see that the ideal flow
   ( $
    p=R\rho\theta,\ e=C_v\theta,
    \ \textrm{ and }\ C_v,\kappa,\mu,\lambda=\textrm{constants}
    $) satisfies the above assumptions with $p_0=p_1=\rho^\frac{1}{2}$ and $p_2=1$. The above
assumptions are motivated by the facts in \cite{Becker,Zeldovich}
where it is pointed out that  $\mu$,  $\lambda$ and $\kappa$ of a
real gas are vary with the temperature and   density, $e$ grows as
$\theta^{1+\gamma}$ with $\gamma\approx0.5$ and $\kappa$ increases
like $\theta^q$ with $q\in[4.5,5.5]$.

In the case that the data $(\rho_0,\theta_0,u_0,f,h)$ are
sufficiently regular and the initial density $\rho_0$ is bounded
away from zero, there exits a unique local strong solution to the
problem (\ref{2-b-E1.1})-(\ref{2-b-E1.4-3}), and the solution
exists globally in time provided that the initial data are small
in some sense. For details, we refer the readers to   papers
\cite{Danchin00,Kawashita02,Marsumura83,Solonnikov80} and the
references therein. When $\mu,\lambda,\kappa$ are three constants,
there have been some existence results on the strong solutions for
the general case of nonnegative initial density. R. Salvi and I.
Stra\u{s}kraba showed in \cite{Salvi93} that if $\Omega$ is a
bounded domain, $p(\cdot)\in C^2[0,\infty)$, $\rho_0\in H^2$,
$u_0\in H^1_0\cap H^2$ and the compatibility condition:
    \begin{equation}
    -\mu\mathrm{div}(\nabla u_0+\nabla u_0^\top)-
    \lambda\nabla\mathrm{div}u_0+\nabla p(\rho_0)
    =\rho_0^\frac{1}{2}g,
    \textrm{ for  } g\in L^2,
    \label{2-b-E1.3-1}
    \end{equation}
is satisfied, then there exists a unique local strong solution
$(\rho,u)$ to the initial boundary value problem for the
isentropic fluids. Independently of their work, H. J. Choe and H.
Kim \cite{Choe03} proved a similar existence result when $\Omega$
is either a bounded domain or the whole space,
$p=a\rho^\gamma(a>0,\gamma>1)$, $\rho\in L^1\cap H^1\cap W^{1,6}$,
$u_0\in D^1_0\cap D^2$ and the condition (\ref{2-b-E1.3-1}) is
satisfied. Later, there are some results in \cite{Cho04,Cho05-1}.

Considering the physically important case that   coefficients
$\mu$, $\lambda$ and $\kappa$ are not constants, some results had
obtained in one-dimension, see
\cite{Jiang98,Straskraba03,yang3,zhang2006-2} etc.;
Feireisl\cite{Feireisl04} proved the existence of globally defined
variational solutions to the compressible Navier-Stokes equations
in $\mathbb{R}^N$ with temperature-dependent viscosity; and some
results had obtained for the incompressible Navier-Stokes
equations, see \cite{Cho05-2,Choe03-2,Lions96} etc.. In
\cite{Cho05-2}, Cho-Kim study the case that
$(\mu,\lambda,\kappa)=(\mu,\lambda,\kappa)(\rho,\rho\theta)$, and
this method is fail to our case that
$(\mu,\lambda,\kappa)=(\mu,\lambda,\kappa)(\rho,\theta)$, since
the density will be degenerate. In this paper, we shall combine
the ideas in \cite{Cho05-1,Cho05-2,Choe03-2,Feireisl04,Lions96} to
study strong solutions to the initial boundary value problem
(\ref{2-b-E1.1})-(\ref{2-b-E1.4-3}) in a domain
$\Omega\subset\mathbb{R}^3$, with nonnegative initial density and
$(e,p,\kappa,\mu,\lambda)=(e,p,\kappa,\mu,\lambda)(\rho,\theta)$.
The lack of simplicity of real flow
$(e,p,\kappa,\mu,\lambda)=(e,p,\kappa,\mu,\lambda)(\rho,\theta)$
makes the analysis significantly different from the ideal flow (
    $
    p=R\rho\theta,\ e=C_v\theta,
    \ \textrm{ and }\ C_v,\kappa,\mu,\lambda=\textrm{constants}
    $).

 Throughout this paper, we use the following
simplified notations for the standard homogeneous and
inhomogeneous Sobolev spaces:
    $$
    L^r=L^r(\Omega),
    \ W^{k,r}=W^{k,r}(\Omega),
    \ H^k=W^{k,2},
    $$
        $$
        D^{k,r}=\{v\in L^1_{loc}(\Omega):
        \|\nabla^kv\|_{L^r}<\infty\},
        \ D^k=D^{k,2},
        $$
    $$
    D^1_0=\{v\in L^6:\|\nabla v\|_{L^2}<\infty\textrm{ and }
    v=0\textrm{ on }\partial \Omega\},
    \ H^1_0=D^1_0\cap L^2.
    $$
It follows from the classical Sobolev embedding results that
    $$
    \|w\|_{L^6}\leq C\|w\|_{D^1_0},
    \ \|w\|_{L^\infty}\leq C\|w\|_{W^{1,q}},
    \ \|w\|_{L^\infty}\leq C\|w\|_{D^{1}_0\cap D^{1,q}}
    $$
and
    $$
  \|w\|_{L^\infty}\leq C\|w\|_{D^1_0\cap D^2},
  \ \|w\|_{L^\infty}\leq C\|w\|_{L^6\cap D^{1,q}},
    $$
provided that $q>3$. Hereafter, we denote by $C$ a generic
positive constant depending only on the fixed constant  $q$ and
the norms of $(h,f)$. When $\Omega$ is   a bounded domain in
$\mathbb{R}^3$ with smooth boundary, we choose $\theta_c\in
C^\infty(\mathbb{R}^3)$ satisfying $\theta_c\equiv\theta_b$ near
$\partial\Omega$. When $\Omega$ is the half space $\mathbb{R}^3_+$
or an exterior domain in $\mathbb{R}^3$ with smooth boundary, we
choose $\theta_c\in C^\infty(\mathbb{R}^3)$ satisfying
$\theta_c\equiv\theta_b$ near $\partial\Omega$ and
$\theta_c|_{B_R^c}\equiv\theta^\infty$ where $B_R^c=\{|x|\geq R\}$
and $R$ is large enough. When $\Omega=\mathbb{R}^3$, we choose
$\theta_c\in C^\infty(\mathbb{R}^3)$ satisfying
$\theta_c|_{B_R^c}\equiv\theta^\infty$.

Our existence and uniqueness result can be stated as follows.

\begin{thm}\label{2-b-T1}
  Let $\rho^\infty,\theta_b,\theta^\infty\geq0$ and $3<q< 6$ be fixed constants.
  Assume that $(\rho_0,\theta_0,u_0,h,f)$ satisfy the
  regularity conditions
    \begin{equation}
    \rho_0\geq0,\ \rho_0-\rho^\infty\in L^6\cap D^1\cap D^{1,q},
    \ (\theta_0-\theta_c,u_0)\in D^1_0\cap D^2,\label{2-b-E3.7}
    \end{equation}
            \begin{equation}
            (h,f)\in C([0,\infty);L^2)\cap L^2(0,\infty;L^q)
            \textrm{ and }(h_t,f_t)\in L^2(0,\infty;H^{-1})
            \label{2-b-E3.8}
            \end{equation}
  and the compatibility condition
    \begin{equation}
    \left\{\begin{array}{l}
      -\mathrm{div}(\kappa_0\nabla
      \theta_0)-Q(\rho_0,\theta_0,u_0)=\rho_0^{\frac{1}{2}}g_1,\\
      L(\rho_0,\theta_0,u_0)+\nabla p(\rho_0,\theta_0)=\rho^{\frac{1}{2}}_0g_2,
      \end{array}   \right.   \textrm{ in }\Omega
      \label{2-b-E3.9-1}
    \end{equation}
  for some $(g_1,g_2)\in L^2$,where $\kappa_0=\kappa(\rho_0,\theta_0)$, $Q(\rho,\theta, u)=\frac{\mu(\rho,\theta)}{2}|\nabla
u+\nabla u^\top|^2+\lambda(\rho,\theta)(\mathrm{div}u)^2$ and
$L(\rho,\theta,u)=-\mathrm{div}(\mu(\rho,\theta)(\nabla u+\nabla
u^\top))-\nabla(\lambda(\rho,\theta)\mathrm{div}u)$. Then
  there exist a small time $T_*>0$ and a unique strong solution
  $(\rho,\theta,u)$ to the initial boundary value problem
  (\ref{2-b-E1.1})-(\ref{2-b-E1.4-3}) such that
        \begin{equation}
          \rho-\rho^\infty\in C([0,T_*];L^6\cap D^1\cap D^{1,q}),
          \ \rho_t\in C([0,T_*];L^2\cap L^q),
          \label{2-b-E3.8-1}
        \end{equation}
            \begin{equation}
              (\theta-\theta_c,u)\in C([0,T_*];D^1_0\cap D^2)\cap
              L^2(0,T_*;D^{2,q}),
              \label{2-b-E3.8-2}
            \end{equation}
        \begin{equation}
          (\theta_t,u_t)\in L^2(0,T_*;D^1_0)
          \textrm{ and
          }(\rho^{\frac{1}{2}}\theta_t,\rho^\frac{1}{2}u_t)\in
          L^\infty([0,T_*];L^2),
          \label{2-b-E3.8-3}
        \end{equation}
            \begin{equation}
               \|\rho(t,\cdot)-\rho_0(\cdot)\|_{L^6\cap D^1\cap D^{1,q}}
    +\|(\theta(t,\cdot)-\theta_0(\cdot),u(t,\cdot)-u_0(\cdot))\|_{D^1_0\cap D^{2}}
    \rightarrow0,\ \textrm{ as }\ t\rightarrow0.
    \label{2-b-E3.8-4}
            \end{equation}
Furthermore, we have the following blow-up criterion: If $T^*$ is
the maximal existence time of the strong solution
$(\rho,\theta,u)$ and $T^*<\infty$, then
    \begin{eqnarray}
      &&\limsup_{t\nearrow T^*}(\|\rho(t,\cdot)-\rho^\infty\|_{L^6\cap D^1\cap D^{1,q}}
      +\|\sqrt{\rho}\theta_t(t,\cdot)\|_{L^2}\nonumber\\
      &&+\int^t_0\|\nabla \theta_t(s,\cdot)\|_{L^2}^2ds+\|\theta(t,\cdot)-\theta_c\|_{ D^{1}_0\cap D^{1,6}}
      +\|u(t,\cdot)\|_{D^1_0})=\infty.
      \label{2-b-E3.13-1}
    \end{eqnarray}
\end{thm}
\begin{rem} The $L^6\cap D^1\cap D^{1,q}$-regularity of $\rho_0-\rho^\infty$
seems inevitable to prove the local well-posedness in the
framework of Sobolev spaces for any compressible fluid model in
three dimensions, because the Sobolev embedding $L^6\cap
D^{1,q}\hookrightarrow L^\infty$ holds only for $q>3$  in
$\mathbb{R}^3$.  Moreover, we allow $\rho_0$ and $\rho^\infty$ to
vanish and so we may consider both interior vacuum and vacuum at
infinity. Here, we only assume that $\rho_0-\rho^\infty\in L^6\cap
D^1\cap D^{1,q}$, which is weaker than the assumption  in
\cite{Cho05-2}.
\end{rem}
\begin{rem}
  To simplify the presentation, we only give the proof of the case
  that $\theta_b=\theta^\infty=\theta_c=0$.
\end{rem}
\begin{rem}
  As has been observed in
\cite{Cho04,Cho05-1,Cho05-2,Choe03,Choe03-2,Salvi93}, the lack of
a positive lower bound of $\rho_0$ should be compensated with the
compatibility condition (\ref{2-b-E3.9-1}). If $(\rho,\theta,u)$
is a sufficiently smooth solution to
(\ref{2-b-E1.1})-(\ref{2-b-E1.2}), then letting $t\rightarrow0$ in
the equations (\ref{2-b-E1.3})-(\ref{2-b-E1.2}), we easily derive
a natural condition: there exists a pair $(g_1,g_2)$ of scalar and
vector fields such that
        \begin{equation}
          -\mathrm{div}(\kappa_0\nabla  \theta_0)
          -Q(\rho_0,\theta_0, u_0)=\rho_0g_1
          \textrm{ and }
          L(\rho_0,\theta_0,u_0)+\nabla p(\rho_0,\theta_0)=\rho_0g_2
          \textrm{ in }\Omega.
          \label{2-b-E1.4-2}
        \end{equation}
But it turns out that the weaker condition (\ref{2-b-E3.9-1}) is
sufficient to prove the local existence and uniqueness of strong
solution. Roughly speaking, the compatibility condition
(\ref{2-b-E3.9-1}) is equivalent to the $L^2$-integrability of
$\sqrt{\rho}\theta_t$ and $\sqrt{\rho}u_t$ at $t=0$. Consider the
elliptic system with variable coefficients, we obtain some
regularity results in Section \ref{2-b-Sec4}. Under these elliptic
regularity results and the condition (\ref{2-b-E3.9-1}), we could
deduce that $(\theta_t,u_t)\in L^2(0,T;D^1_0)$ and
$(\sqrt{\rho}\theta_t,\sqrt{\rho}u_t)\in L^\infty(0,T;L^2)$ for a
small time $T>0$. Here, the compatibility condition
(\ref{2-b-E3.9-1}) is satisfied automatically for all initial data
$(\rho_0,\theta_0,u_0)$ with the regularity (\ref{2-b-E3.7})
whenever $\rho_0$ is bounded away from zero.
\end{rem}
\begin{rem}
Consider the isentropic case and $\mu,\lambda$ being constants, B.
Desjardins \cite{Desjardins97} prove the local existence of a weak
solution $(\rho,u)$ with a bounded nonnegative density to the
periodic problem as long as
    $$ \sup_{0\leq t\leq
    T^*}(\|\rho(t)\|_{L^\infty(\mathbb{T}^3)}+ \|\nabla
    u(t)\|_{L^2(\mathbb{T}^3)})<\infty;$$
 Y. Cho, H. J. Choe and H. Kim
\cite{Cho04} proved the local existence of a weak solution
$(\rho,u)$ as long as
        $$
        \limsup_{t\nearrow T^*}(\|\rho(t)\|_{H^{1}\cap W^{1,q}}
      +\|u(t)\|_{D^1_0})<\infty.
      $$
Using   similar arguments in our proof, we can see that in the
special case that the function $e_\theta=e_\theta(\rho)$ is
independent of the function $\theta$, (\ref{2-b-E3.13-1}) may be
replaced by
        \begin{equation}
      \limsup_{t\nearrow T^*}(\|\rho(t,\cdot)-\rho^\infty\|_{L^6\cap D^1\cap D^{1,q}}
      +\|(\theta-\theta_c,u)(t,\cdot)\|_{ D^{1}_0\cap D^{1,6}})=\infty.
        \end{equation}
\end{rem}

Furthermore, since the local existence time $T_*$ and all
regularity   estimates  of  the strong solution  in Theorem
\ref{2-b-T1} depend only on  $\|(g_1,g_2)\|_{L^2}$ and the norms
of the initial data $(\rho_0,\theta_0,u_0,f,h)$. The proof of
Theorem \ref{2-b-T1} also shows that the continuous dependence of
the solution on the initial data holds at least for a small time
interval. We may state the following result without a proof.

\begin{thm}\label{2-b-T2}
  For each $i=1,2$, let $(\rho_i,\theta_i,u_i)$ be the local strong
  solution to the initial boundary value problem
   (\ref{2-b-E1.1})-(\ref{2-b-E1.4-3}) with the initial data
  $(\rho_{0i},\theta_{0i},u_{0i})$, which satisfy regularity
  conditions (\ref{2-b-E3.7})-(\ref{2-b-E3.8}) and the compatibility
  condition (\ref{2-b-E3.9-1}) with $(g_1,g_2)=(g_{1i},g_{2i})$.
  Assume also that for each $i=1,2$,
        \begin{eqnarray*}
    &&\|\rho_{0i}-\rho^\infty\|_{L^6\cap D^1\cap D^{1,q}
    }+\|(\theta_{0i}-\theta_c,u_{0i})\|_{D^1_0\cap D^2}
    +\|(g_{1i},g_{2i})\|_{L^2}\\
            &&+\|(h,f)\|_{L^\infty_tL^2_x\cap L^2_tL^q_x}
            +\|(h_t,f_t)\|_{L^2_tH^{-1}_x}\leq K
        \end{eqnarray*}
  and
$\overline{\rho}_0=\rho_{01}-\rho_{02}\in X$,
$(\sqrt{\rho_{01}}(\theta_{01}-\theta_{02}),\sqrt{\rho_{01}}(u_{01}-u_{02}))\in
L^2$, where  $X=L^2\cap L^\frac{3}{2}$, if $\rho^\infty=0$, $X=
L^2$, if $\rho^\infty\neq0$.
  Then there exist a small time $T_{**}$ and a positive constant
  $C_K$, such that
        \begin{eqnarray*}
      &&\sup_{t\in[0,T_{**}]}\left(\|\overline{\rho}\|^2_{X}
      +\|\sqrt{\rho_1}\overline{\theta}\|_{L^2}^2
      +\|\sqrt{\rho_1}\overline{u}\|_{L^2}^2
      \right)+\int^{T_{**}}_0\left(\|\overline{\theta}\|_{D^1_0}
      +\|\overline{u}\|_{D^1_0}\right)dt\\
            &\leq&C_K\left(\|\overline{\rho}_0\|^2_{X}
      +\|\sqrt{\rho_{01}}\overline{\theta}_0\|_{L^2}^2
      +\|\sqrt{\rho_{01}}\overline{u}_0\|_{L^2}^2
      \right),
    \end{eqnarray*} where $\overline{\rho}=\rho_1-\rho_2$,
$\overline{\theta}=\theta_1-\theta_2$ and $\overline{u}=u_1-u_2$.
\end{thm}

The global result of the initial boundary value problem
   (\ref{2-b-E1.1})-(\ref{2-b-E1.4-3}) with $\rho_0\geq0$ is very difficult to
obtained. For a special case that the temperature vanish in the
vacuum state, we can give a blow-up result in  Section
\ref{2-c-s1}.

In Section \ref{2-c-s1}, we consider the following equations for a
compressible fluid in $\mathbb{R}^n\times\mathbb{R}_+(n\geq1)$:
\begin{equation}
      \rho_t+\mathrm{div } (\rho u)=0,
      \label{2-c-E1.1}
    \end{equation}
            \begin{equation}
              (\rho u)_t+\mathrm{div} (\rho u\otimes u)
              +\nabla p=\mathrm{div}(\mu(\nabla u+\nabla u^\top))+
              \nabla(\lambda\mathrm{div} u),
              \label{2-c-E1.3}
            \end{equation}
            \begin{equation}
          (\rho e)_t+\mathrm{div }(\rho e u)
          +p\mathrm{div }u=\frac{\mu}{2}|\nabla u+\nabla
          u^\top|^2+\lambda(\mathrm{div } u)^2+\mathrm{div}(\kappa\nabla e),
          \label{2-c-E1.2}
        \end{equation}
where $(\kappa,\mu,\lambda)=(\kappa,\mu,\lambda)(\rho,\theta)$. In
Section \ref{2-c-s1}, we will consider only polytropic gas, so
that the equations of the state of the gas is given by
    \begin{equation}
      p=R\rho \theta,
      \ e=c_v\theta,
      \ \textrm{ and }
      \ p=A\exp(S/c_v)\rho^\gamma,
      \label{2-c-E1.4-1}
    \end{equation}
where $R>0$ is the gas constant, $A$ is a positive constant of
absolute value, $\gamma>1$ is the ratio of specific heat,
$c_v=\frac{R}{\gamma-1}$ is  the specific heat at constant volume
and $S$ is the entropy. In Section \ref{2-c-s1}, we only assume
that
    \begin{equation}
    \mu(0,0)>0,\ \ n\lambda(0,0)+2\mu(0,0)= \lambda_0> 0
    \ \textrm{ and }\ \kappa(s_1,s_2)\geq 0,
    \label{2-c-E1.4}
    \end{equation}
 for all $s_1,s_2\geq0$. The Navier-Stokes system is supplemented
 with the initial data
    $$
    (\rho,S,u)|_{t=0}=(\rho_0,S_0,u_0)(x)\in H^k(\mathbb{R}^n),
    \ k>[\frac{n}{2}]+2.
    $$

The blow-up of smooth solutions of compressible Euler equations
$(\mu=\lambda=\kappa=0)$ has been studied by several
mathematicians. T.C. Sideris \cite{Sideris85} showed that the life
span $T$ of the $C^1$ solution to the compressible Euler equations
is finite when the initial data have compact support and the
initial velocity is sufficiently large (super-sonic) in some
region. When the initial density and velocity have compact
supports, T. Makino, S. Ukai and S. Kawashima \cite{Makino86}
shown the blow-up of smooth solutions to the compressible Euler
equations in  $\mathbb{R}^3$ without external force and heat
source. In a different way from \cite{Makino86} and
\cite{Sideris85}, Z. Xin \cite{xin} showed   the same blow-up
result for the compressible Euler equations. Considered the case
that $\mu,\lambda$
 are constants and $\kappa= 0$, Xin\cite{xin}
 also showed similar results for
the compressible Navier-Stokes equations for polytropic gas, when
the initial density has compact support.  In \cite{Cho2004},
authors extended Z. Xin's result\cite{xin} to the case of
$\kappa>0$ and $\mu$, $\lambda$ being constants. But their point
of view in \cite{Cho2004,xin} cannot be applied for the problem in
Section \ref{2-c-s1}, since the fact that
    $$
    \int_{\mathbb{R}^n}    [\mathrm{div}(\mu(\nabla u+\nabla u^\top))+
              \nabla(\lambda\mathrm{div} u)]\cdot xdx=0
    $$
is strongly necessary in their argument, which seems hard to be
obtained here.  Here, we extend the Xin's blow-up result to the
case that $\mu,\lambda,\kappa$ depend on the density and the
temperature.

Before stating our blow-up result, we introduce some notations. We
denote by $B_r=B_r(0)$ the ball in $\mathbb{R}^n$ of the radius
$r$ centered at the origin. We will use several physical
quantities:
    $$
    m_0=\int_{\mathbb{R}^n}\rho_0dx
    \textrm{ (initial total mass)},
    $$
        $$
        m_1=\int_{\mathbb{R}^n}\rho_0(x)xdx,
        $$
    $$
    m_2=\int_{\mathbb{R}^n}\rho_0(x)u_0(x)dx
    \textrm{ (initial total momentum)}.
    $$
We always assume that $m_0>0$.

For the proof of the blow-up, we have only to prove the following
theorem.
\begin{thm}\label{2-c-T1.1}
We assume that $(\mu,\lambda,\kappa)$ satisfies (\ref{2-c-E1.4})
and $T>0$. Suppose that $(\rho,S,u)\in
C^1([0,T];H^k(\mathbb{R}^n))$, $k>[\frac{n}{2}]+2$, is a solution
to the Cauchy problem (\ref{2-c-E1.1})-(\ref{2-c-E1.4-1}) with the
initial data $(\rho_0,S_0,u_0)$, where the initial density
$\rho_0$ is compactly supported in a ball $B_{r_0}$. Then we have
    $$
    m_0r_0\geq m_{1i}+m_{2i}T\geq -m_0r_0,
    \ i=1,\ldots,n.
    $$
\end{thm}
\begin{rem}
If $T^*$ is the maximal existence time of the solution
$(\rho,S,u)$. Then
  since $m_0$ is strictly positive, this theorem implies
  that $T^*$ should be finite when $m_2\neq0$. This theorem also
  shows the relationship between the size of support and the life
  span. For example, the range of life span can be extended as the
 support of the initial density become larger. Especially, we can
expect the global existence of the smooth solution to the
compressible Navier-Stokes equations in the case that the initial
density is positive but has decay at infinity.
\end{rem}
\begin{rem}
As in \cite{Cho2004,xin}, our method need that the temperature
vanish in the vacuum state, so we consider the entropy $S$ and
equation (\ref{2-c-E1.4-1}) which guarantee the temperature vanish
in the vacuum state.
\end{rem}

The rest of this paper is organized as follows. First, based on
some regularity results on elliptic system in Section
\ref{2-b-Sec4}, we obtain some \textit{a priori} estimates for the
linearized problem in Section \ref{2-b-Sec2}. In Section
\ref{2-b-Sec2-1}, we will prove the existence results for the
linearized problem. Then, using the usual iteration argument, we
will prove
 Theorem \ref{2-b-T1} in Section \ref{2-b-Sec3}. In Section
\ref{2-c-s1}, we will prove a  blow-up  result.
 \section{\textit{A priori} estimates for the linearized problem}\label{2-b-Sec2}\setcounter{equation}{0}
In this section, we consider the following linearized problem
    \begin{equation}
      \rho_t+\mathrm{div } (\rho v)=0,\label{2-b-E2.1}
    \end{equation}
             \begin{equation}
              \rho u_t+\rho v\cdot\nabla u+L(\sigma,\eta,u)
              +\nabla p(\rho,\theta)=\rho f,
             \ \textrm{in }(0,T)\times\Omega,
              \label{2-b-E2.3}
            \end{equation}
        \begin{eqnarray}
          &&e_\theta(\sigma,\eta)\left(\rho \theta_t+\rho v\cdot\nabla\theta
    \right)- \mathrm{div}(\kappa(\sigma,\eta)
          \nabla \theta)
          +p(\rho,\theta)\mathrm{div }v\nonumber\\
    &=&Q(\sigma,\eta,v)+\rho h+e_\rho(\sigma,\eta)\rho^2\mathrm{div}v,\label{2-b-E2.2}
        \end{eqnarray}
        \begin{equation}
          (\rho,\theta,u)|_{t=0}=(\rho_0,\theta_0,u_0)
          \textrm{ in }\Omega,\label{2-b-E2.4}
        \end{equation}
            \begin{equation}
              (\theta,u)=(0,0)
              \textrm{ on }(0,T)\times\partial\Omega,
            \end{equation}
                \begin{equation}
                  (\rho, e, u)(t,x)
                  \rightarrow(\rho^\infty,0,0)\textrm{ as }
                  |x|\rightarrow\infty,
                  \label{2-b-E2.6}
                \end{equation}
where $\sigma,\eta,v$ are three known functions on
$(0,T)\times\Omega$. Recall that $Q(\sigma,\eta,
v)=\frac{\mu(\sigma,\eta)}{2}|\nabla v+\nabla
v^\top|^2+\lambda(\sigma,\eta)(\mathrm{div }v)^2$ and
$L(\sigma,\eta,u)=-\mathrm{div}(\mu(\sigma,\eta)(\nabla u+\nabla
u^\top))-\nabla(\lambda(\sigma,\eta)\mathrm{div }u)$.

Throughout this section, we assume that the known data satisfy the
following regularities
    \begin{equation}
      \rho_0\geq0,\ \rho_0-\rho^\infty\in L^6\cap D^1\cap D^{1,q},
      \ (\theta_0,u_0)\in D^1_0\cap D^2,\label{2-b-E2.21-1}
    \end{equation}
        \begin{equation}
          (h,f)\in C([0,\infty);L^2)\cap L^2(0,\infty;L^q),
          \ (h_t,f_t)\in L^2(0,\infty;H^{-1}),
          \label{2-b-E2.21-2}
        \end{equation}
    \begin{equation}
      (\eta, v)\in C([0,T]; D^1_0\cap D^2)\cap L^2(0,T;D^{2,q}),
      \ (\eta_t,v_t)\in L^2(0,T;D^1_0),
      \label{2-b-E2.21-3}
    \end{equation}
        \begin{equation}
          \sigma-\sigma^\infty\in C([0,T];L^6\cap D^1\cap D^{1,q}),
    \ \sigma_t\in C([0,T];L^2\cap L^q),\\
    \sigma^\frac{1}{2}\eta_t
    \in L^\infty([0,T];L^2),\label{2-b-E2.21-5}
        \end{equation}
for three constants $\rho^\infty,\sigma^\infty\geq0$ and
$q\in(3,6)$. Then the existence of a unique strong solution
$(\rho,\theta,u)$ to the linearized problem
(\ref{2-b-E2.1})-(\ref{2-b-E2.6}) can be proved by standard
methods at least for the case that $\rho_0$ is bounded away from
zero, i.e., $\rho_0\geq\delta$. Using  similar arguments in
\cite{Cho04,Cho05-2,Choe03-2}, applying the classical method of
characteristics, a semi-discrete Galerkin method and the usual
domain expansion technique, we can obtain the following two
lemmas.

\begin{lem}\label{2-b-L2.1}
  Assume that $\rho_0$ and $v$satisfy the regularities
    $$
    \rho_0\geq0,\ \rho_0-\rho^\infty\in C_0
    \textrm{ and } v\in L^\infty([0,T];D^1_0\cap D^2)\cap L^2(0,T;D^{2,q})
    $$
for two constants $\rho^\infty\geq0$ and $q\in(3,6)$, where the
space $C_0$ consists of continuous functions on $\Omega$ vanishing
at infinity. Then there exists a unique weak solution $\rho$ in
$\rho^\infty+C([0,T];C_0)$ to the linear hyperbolic problem
(\ref{2-b-E2.1}), (\ref{2-b-E2.4}) and (\ref{2-b-E2.6}). Moreover
the solution $\rho$ can be represented by
    \begin{equation}
      \rho(t,x)=\rho_0(U(0,t,x))\exp \left\{-\int^t_0
      \mathrm{div }v(s,U(s,t,x))ds\right\},
      \label{2-b-E2.7}
    \end{equation}
where $U\in C([0,T]\times[0,T]\times\overline{\Omega})$ is the
solution to the initial value problem
    \begin{equation}
      \left\{\begin{array}{rlll}
      \frac{\partial}{\partial
      t}U(t,s,x)&=&v(t,U(t,s,x)),&t\in[0,T]\\
            U(s,s,x)&=&x,&s\in[0,T],x\in\overline{\Omega}.
      \end{array}
      \right.\label{2-b-E2.8}
    \end{equation}
Assume in addition that $\rho_0-\rho^\infty\in L^6\cap D^1\cap
D^{1,q}$. Then, we also have
    \begin{equation}
      \rho-\rho^\infty\in C([0,T];L^6\cap D^1\cap D^{1,q})
      \textrm{ and } \rho_t\in L^\infty([0,T];L^2\cap L^q),
    \end{equation}
and
    \begin{eqnarray}
      \|\rho(t,\cdot)-\rho^\infty\|_{L^6\cap D^1\cap D^{1,q}}
      &\leq& \left(\|\rho_0-\rho^\infty\|_{L^6\cap D^1\cap D^{1,q}}
      +C\rho^\infty\int^t_0\| v(s,\cdot)\|_{D^{2}\cap D^{2,q}}ds
      \right)\nonumber\\
            &&\times\exp\left(C\int^t_0\|\nabla
            v(s,\cdot)\|_{W^{1,q}\cap D^1}ds
            \right)\label{2-b-E2.20}
    \end{eqnarray}
\end{lem}

\begin{lem}\label{2-b-L2.2}
  Assume in addition to
(\ref{2-b-E2.21-1})-(\ref{2-b-E2.21-5})
  that $\rho_0\geq \delta$ in $\Omega$ for a constant
  $\delta>0$. Then there exists a unique strong solution
  $(\rho,\theta,u)$ to the initial boundary value problem
  (\ref{2-b-E2.1})-(\ref{2-b-E2.6}) in $[0,T]$ such that
        \begin{equation}
          \rho-\rho^\infty\in C([0,T];L^6\cap D^1\cap D^{1,q}),
    \ \rho_t\in C([0,T];L^2\cap L^q),
          \label{2-b-E2.22-1}
        \end{equation}
            \begin{equation}
              (\theta,u)\in C([0,T];D_0^1\cap D^2)\cap L^2(0,T;D^{2,q}),
              \label{2-b-E2.22-2}
            \end{equation}
        \begin{equation}
          (\theta_t,u_t)\in L^2(0,T;D^1_0)\cap  C([0,T];L^2),
          \label{2-b-E2.22-3}
        \end{equation}
                \begin{equation}
                  \rho\geq \underline{\delta}
                  \textrm{ on }[0,T]\times\overline{\Omega}
                  \textrm{ for a constnat
                  }\underline{\delta}>0,
                \end{equation}
        \begin{equation}
          \|\rho(t,\cdot)-\rho_0(\cdot)\|_{L^6\cap D^1\cap D^{1,q}}
          +\|(\theta(t,\cdot)-\theta_0(\cdot),u(t,\cdot)-u_0(\cdot))\|_{D^1_0\cap
          D^2}\rightarrow 0,
          \ \textrm{ as }t\rightarrow0.
        \end{equation}
\end{lem}

Assume that the initial data $(\rho_0,\theta_0,u_0)$  satisfy the
hypotheses of Lemma \ref{2-b-L2.2} and the compatibility condition
    \begin{equation}
         \left\{\begin{array}{l}
          -\mathrm{div}(\kappa(\sigma(0),\eta(0))\nabla  \theta_0)
          -Q(\sigma(0),\eta(0), v(0))=\rho_0^\frac{1}{2}g_1,\\
          L(\sigma(0),\eta(0),u_0)+\nabla p(\rho_0,\theta_0)=\rho_0^\frac{1}{2}g_2
          \textrm{ in }\Omega.
          \end{array}\right.\label{2-b-E2.20-1}
        \end{equation}
 Let us choose a fixed $c_0$ so that
    $$
    c_0\geq 1+\rho^\infty+\|\rho_0-\rho^\infty\|_{L^6\cap D^1\cap D^{1,q}}
    +\|(\theta_0,u_0)\|_{D^1_0\cap D^2}+\|(g_1,g_2)\|^2_{L^2},
    $$
 and
    $$
    c_0\geq\sigma^\infty+ \|\sigma(0,\cdot)-\sigma^\infty\|_{L^6\cap D^1\cap D^{1,q}}
    +\|(\eta(0,\cdot),v(0,\cdot))\|_{D^1_0\cap D^2},
    $$
where $3<q< 6$ and $\sigma^\infty, \rho^\infty\geq0$ are three
constants. Moreover, assume that $\sigma,\eta,v$ satisfy the
regularity stated in Lemma \ref{2-b-L2.2} and
    \begin{equation}
      \sup_{t\in[0,T]}\|v(t,\cdot)\|_{D^1_0}
      +\int^{T}_0\left(\|v_t(t,\cdot)\|_{D^1_0}^2+\|v(t,\cdot)\|^2_{D^{2,q}}
      \right)dt\leq c_1,
      \label{2-b-E2.24}
    \end{equation}
        \begin{equation}
          \sup_{t\in[0,T]}\|v(t,\cdot)\|_{D^2}\leq c_2,
        \end{equation}
        \begin{equation}
          \sup_{t\in[0,T]}\left(\|(\sqrt{\sigma}\eta_t)(t,\cdot)\|_{L^2}+\|\eta(t,\cdot)\|_{D^2}
         \right)+\int^{T}_0\|\eta_t(t,\cdot)\|_{D^1_0}^2dt\leq
          c_3,\label{2-b-E2.24.2}
        \end{equation}
        \begin{equation}
    \sup_{t\in[0,T]}\left(\|\sigma_t(t,\cdot)\|_{L^2\cap L^{q}}
         +\|\sigma(t,\cdot)-\sigma^\infty\|_{L^6\cap D^1\cap D^{1,q}}^2
         \right)\leq c_4,\label{2-b-E2.24.4}
        \end{equation}
for some fixed constants $c_1$, $c_2$, $c_3$ and $c_4$ such that
    $$
    1<c_0<c_1<c_2,\ 1<c_0<c_3,c_4.
    $$
Then we could derive some \textit{a priori} estimates for the
solution $(\rho,\theta,u)$ which are independent of $\delta$.

\textbf{I. Estimates for the density and the coefficients
($\mu,\lambda,\kappa,e_\theta,e_\rho$).}

To estimate the density $\rho$, we first recall from
(\ref{2-b-E2.20}) that
    $$
    \|\rho(t,\cdot)-\rho^\infty\|_{L^6\cap D^1\cap D^{1,q}}\leq Cc_0\exp\left(C\int^t_0\|\nabla
            v(s,\cdot)\|_{H^1\cap W^{1,q}}ds
            \right)
    $$
for $t\in[0,T]$. Then, using H\"{o}lder's inequality, we obtain
    $$
    \int^t_0\|\nabla  v\|_{H^1\cap W^{1,q}}ds\leq
    t^{\frac{1}{2}}\left[\int^t_0\|\nabla  v\|^2_{H^1\cap W^{1,q}}ds
    \right]^{\frac{1}{2}}\leq C(c_2t+(c_2t)^{\frac{1}{2}}),
    $$
 and conclude that
    \begin{equation}
      \|\rho(t,\cdot)-\rho^\infty\|_{L^6\cap D^1\cap D^{1,q}}
      \leq Cc_0,
      \label{2-b-E2.25}
    \end{equation}
for $0\leq t\leq \min(T,T_0)$, where $T_0=c_2^{-1}<1$. Moreover,
it follows from (\ref{2-b-E2.7})  and (\ref{2-b-E2.25}) that
    \begin{equation}
      C^{-1}\delta \leq \rho(t,x)\leq Cc_0
      \textrm{ for }0\leq t\leq \min(T,T_0),
      x\in\overline{\Omega}.
    \label{2-b-E2.26}
    \end{equation}

From (\ref{2-b-E2.24})-(\ref{2-b-E2.24.2}), using  the
Gagliardo-Nirenberg inequality, we have
    \begin{equation}
      \|\eta(t)-\eta(0)\|_{D^1_0}\leq
      C\int^t_0\|\eta_t\|_{D^1_0}ds\leq
      C\left(\int^t_0\|\eta_t\|^2_{D^1_0}ds\right)^{\frac{1}{2}}t^{\frac{1}{2}}
      \leq Cc_3^{\frac{1}{2}}t^{\frac{1}{2}},
      \label{2-b-E2.23}
    \end{equation}
      \begin{equation}
      \|\eta(t)-\eta(0)\|_{D^{1,q}}\leq C
      \|\eta(t)-\eta(0)\|_{D^{1}_0}^{\frac{6-q}{2q}}
      \|\eta(t)-\eta(0)\|_{D^{2}}^{\frac{3q-6}{2q}}
      \leq Cc_3^{\frac{5q-6}{4q}}t^{\frac{6-q}{4q}},
    \end{equation}
    \begin{equation}
      \|v(t)-v(0)\|_{D^1_0}\leq
      C\int^t_0\|v_t\|_{D^1_0}ds
      \leq Cc_1^{\frac{1}{2}}t^{\frac{1}{2}},
       \label{2-b-E2.23.2}
    \end{equation}
      \begin{equation}
      \|v(t)-v(0)\|_{D^{1,q'}}\leq C
      \|v(t)-v(0)\|_{D^{1}_0}^{\frac{6-q'}{2q'}}
      \|v(t)-v(0)\|_{D^{2}}^{\frac{3q'-6}{2q'}}
      \leq Cc_1^{\frac{6-q'}{4q'}}c_2^{\frac{3q'-6}{2q'}}t^{\frac{6-q'}{4q'}},
    \end{equation}
for $q'=\max(q,4)$. Let
$T_1=\min(T_0,c_0^2c_3^{-1},c_0^{\frac{4q}{6-q}}c_3^{-\frac{5q-6}{6-q}}
,c_0^2c_1^{-1},c_0^{\frac{4q'}{6-q'}}c_1^{-1}c_2^{-\frac{6q'-12}{6-q'}})$,
we have
    \begin{equation}
    \|\eta(t,\cdot)\|_{D^1_0\cap D^{1,q}}+    \|v(t,\cdot)\|_{D^1_0\cap D^{1,q'}}\leq Cc_0,
    \ \forall 0\leq t\leq \min(T,T_1).\label{2-b-E2.24.3}
    \end{equation}
Using Sobolev's embedding theorem, we have
    $$
    \|(\eta,v)(t,\cdot)\|_{L^\infty}\leq Cc_0,
    \ \forall\ 0\leq t\leq \min(T,T_1).
    $$
From  (\ref{2-b-E2.1}), (\ref{2-b-E2.25}) and (\ref{2-b-E2.24.3}),
we obtain
    \begin{equation}
    \|\rho(t,\cdot)-\rho^\infty\|_{L^6\cap D^1\cap D^{1,q}}^2+
    \|\rho_t\|_{L^2\cap L^q}\leq A_0c_0^2,\label{2-b-E2.25.1}
    \end{equation}
where $0\leq t\leq \min(T,T_1)$. Here, the constant $A_0>1$ is
independent of any one of $c_i$, $i=0,\ldots,4$. Using the
interpolation inequality, we have
    $$
    \|\rho(t,\cdot)-\rho_0(\cdot)\|_{L^6}\leq C\|\rho(t,\cdot)-\rho_0(\cdot)\|_{L^2}^{\frac{1}{3}}
    \|\rho(t,\cdot)-\rho_0(\cdot)\|_{L^\infty}^{\frac{2}{3}}\leq
    Cc_0^{2}t^{\frac{1}{3}},
    $$
where $0\leq t\leq \min(T,T_1)$.

 Choosing $c_4=A_0c_0^2$,
from (\ref{2-b-E1.4.1})- (\ref{2-b-E1.5}),
(\ref{2-b-E2.24.2})-(\ref{2-b-E2.24.4}) and (\ref{2-b-E2.24.3}),
we have, for any  $0\leq t\leq \min(T,T_1)$,
    \begin{eqnarray*}
    &&\|\nabla\mu(\sigma(t,\cdot),\eta(t,\cdot))\|_{L^2\cap L^q}\\
    &\leq& \|\mu_{\rho}(\sigma(t,\cdot),\eta(t,\cdot))\nabla\sigma(t,\cdot)\|_{L^2\cap L^q}
    +\|\mu_{\theta}(\sigma(t,\cdot),\eta(t,\cdot))\nabla\eta(t,\cdot)\|_{L^2\cap L^q}\\
    &\leq& C(c_0)\left(\|\nabla\sigma(t,\cdot)\|_{L^2\cap L^q}
    +\|\nabla\eta(t,\cdot)\|_{L^2\cap L^q}\right)\leq C(c_0),
    \end{eqnarray*}
and
    \begin{eqnarray*}
    &&\|\mu_t(\sigma(t,\cdot),\eta(t,\cdot))\|_{L^3}\\
    &\leq& \|\mu_{\rho}(\sigma(t,\cdot),\eta(t,\cdot))\sigma_t(t,\cdot)\|_{L^3}
    +\|\mu_{\theta}(\sigma(t,\cdot),\eta(t,\cdot))\eta_t(t,\cdot)\|_{L^3}\\
    &\leq& C(c_0)\left(\|\sigma_t(t,\cdot)\|_{L^3}
    +\|\sigma^{\frac{1}{4}}\eta_t(t,\cdot)\|_{L^3}\right)\\
                &\leq& C(c_0)\|\sigma_t(t,\cdot)\|_{L^3}
    +C(c_0)\|\sigma^{\frac{1}{2}}\eta_t(t,\cdot)\|_{L^2}^{\frac{1}{2}}
        \|\eta_t(t,\cdot)\|_{L^6}^{\frac{1}{2}}\\
                 &\leq&
                C(c_0)(c_3+c_3^\frac{1}{2}\|\eta_t(t,\cdot)\|_{D^1_0}^{\frac{1}{2}}).
        \end{eqnarray*}
Similarly, we obtain
    \begin{equation}
    \sup_{0\leq t\leq \min(T,T_1)}\left(
    \|(\mu,\lambda,\kappa,e_\theta,e_\rho)(\sigma(t,\cdot),\eta(t,\cdot))\|_{C^\beta}+
    \|(p_0,p_1,p_2)(\rho(t,\cdot))\|_{L^\infty}
    \right)
    \leq c_5,\label{2-b-E2.26-0}
    \end{equation}
    \begin{equation}
    \sup_{0\leq t\leq \min(T,T_1)}
    \|\nabla(\mu,\lambda,\kappa,e_\theta,e_\rho)(\sigma(t,\cdot),\eta(t,\cdot))\|_{L^2\cap L^q}
    \leq c_5,
        \label{2-b-E2.26-1}
    \end{equation}
            \begin{equation}
    \sup_{0\leq t\leq \min(T,T_1)}
    \|\partial_t(\mu,\lambda,\kappa,e_\theta,e_\rho)(\sigma(t,\cdot),\eta(t,\cdot))\|_{L^3}
    \leq c_5\left(c_3+c_3^\frac{1}{2}\|\eta_t(t,\cdot)\|_{D^1_0}^{\frac{1}{2}}\right).
    \label{2-b-E2.26-2}
    \end{equation}
for a constant $c_5>1$ which depend only on $c_0$, and
$\beta=1-\frac{3}{q}\in(0,1)$.

\textbf{II. Estimates for the internal energy and the pressure.}

Taking the operator $\partial_t$ to the equation (\ref{2-b-E2.2}),
we obtain
    \begin{eqnarray*}
    &&e_\theta\rho \theta_{tt}-\mathrm{div}(\kappa\nabla \theta_t)
    -\mathrm{div}(\kappa_t\nabla \theta)+(p\mathrm{div}v)_t\\
        &=&Q(\sigma,\eta,v)_t+(\rho h)_t-(e_\theta\rho v\cdot \nabla \theta)_t
        -(e_\theta\rho)_t\theta_t
        +(e_\rho \rho^2\mathrm{div}v)_t.
    \end{eqnarray*}
Then multiplying this by $\theta_t$, integrating over $\Omega$ and
using (\ref{2-b-E2.1}), we have
    \begin{eqnarray*}
      &&\frac{1}{2}\frac{d}{dt}\int e_\theta\rho \theta^2_tdx+\int\kappa|\nabla
      \theta_t|^2dx+\int(p\mathrm{div}v)_t\theta_tdx
      -\int \theta_t\mathrm{div}(\kappa_t\nabla \theta)dx\nonumber\\
            &=&\int\left(-\frac{1}{2}(e_\theta\rho)_t\theta_t^2
            +\theta_t\left(-(e_\theta \rho v\cdot\nabla
            \theta)_t+
            Q(\sigma,\eta,v)_t+(\rho h)_t+(e_\rho \rho^2\mathrm{div}v)_t\right)
            \right)dx
    \end{eqnarray*}
and
    \begin{eqnarray}
      &&\frac{1}{2}\frac{d}{dt}\int e_\theta\rho \theta^2_tdx+\kappa^0\int|\nabla
      \theta_t|^2dx\nonumber\\
            &\leq&C\int (|p_t||\nabla v||\theta_t|+p|\nabla v_t||\theta_t|
            +|\kappa_t||\nabla \theta|
                    |\nabla \theta_t|
            +\rho|\partial_t e_\theta||\theta_t|^2\nonumber\\
                &&+|\rho v \cdot\nabla (e_\theta\theta_t^2)|
                +|(e_\theta\rho v\cdot\nabla \theta)_t||\theta_t|+(\mu+|\lambda|)|\nabla v||\nabla v_t||\theta_t|
                +(|\mu_t|+|\lambda_t|)|\nabla v|^2|\theta_t|\nonumber\\
                    &&+|\rho_t||h||\theta_t|+|h_t|\rho|\theta_t|
                    +|\theta_t||(e_\rho\rho^2\mathrm{div}v )_t|)dx:=\sum^{11}_{j=1}I_j.
                    \label{2-b-E2.28}
    \end{eqnarray}
Using (\ref{2-b-E2.24}), (\ref{2-b-E2.25})-(\ref{2-b-E2.26}) and
(\ref{2-b-E2.25.1})-(\ref{2-b-E2.26-2}), applying the
Cauchy-Schwartz, Gagliardo-Nirenberg, H\"{o}lder and Sobolev's
inequalities, we can estimate each term $I_j(t)$  for $0\leq t\leq
\min(T,T_1)$ as follows:
        \begin{eqnarray*}
      I_1&\leq&C\int\left(|p_\theta \theta_t|+|p_\rho \rho_t|\right)
      |\nabla v||\theta_t|dx\\
      &\leq& C\int p_1(1+\theta)\rho^\frac{1}{2}|\theta_t|^2|\nabla v|dx
      +C\int p_2(1+\theta^2)|\rho_t||\nabla v||\theta_t|dx\\
            &\leq&Cc_5\left(\|\theta\|_{L^6}\|\nabla
            v\|_{L^6}+\|\nabla v\|_{L^3}
            \right)\|\sqrt{\rho}\theta_t\|_{L^2}\|\theta_t\|_{D^1_0}\\
            &&+Cc_5\left(\|\theta\|_{L^6}^2\|\nabla
            v\|_{L^6}+\|\nabla v\|_{L^2}
            \right)\|\rho_t\|_{L^3}\|\theta_t\|_{D^1_0}\\
      &\leq& Cc_2^6c_5^2\left(1+\|\theta\|_{D^1_0}^4
      \right)\left(1+\|\sqrt{\rho}\theta_t\|_{L^2}^2
      \right)+\frac{\kappa^0}{20}\|\nabla \theta_t\|_{L^2}^2,
    \end{eqnarray*}
        \begin{eqnarray*}
      I_2&\leq&C \int p_0(1+\theta^2)\rho^{\frac{1}{2}}|\nabla
      v_t||\theta_t|dx\\
            &\leq& Cc_5\|\sqrt{\rho}\theta_t\|_{L^2}\|\nabla
            v_t\|_{L^2}(1+\|\theta\|_{L^\infty}^2)\\
                &\leq& \|\nabla v_t\|_{L^2}^2
                +Cc_5^2\|\sqrt{\rho}\theta_t\|_{L^2}^2(1+\|\theta\|_{D^1_0}^4
                +\|\theta\|_{D^2}^4),
    \end{eqnarray*}
    \begin{eqnarray*}
    I_{3}&\leq& C\|\kappa_t\|_{L^3}\|\nabla \theta\|_{L^6}\|\nabla
    \theta_t\|_{L^2}\leq Cc_5^2(c_3^2+c_3\|\eta_t\|_{D^1_0})\|\theta\|^2_{D^2}+\frac{\kappa^0}{20}\|\nabla
    \theta_t\|^2_{L^2}\\
    &\leq&
    Cc_5^2c_3^2\|\theta\|^2_{D^2}+c_3^{-1}\|\eta_t\|^2_{D^1_0}+
    Cc_5^4c_3^3\|\theta\|^4_{D^2}+\frac{\kappa^0}{20}\|\nabla
    \theta_t\|^2_{L^2},
    \end{eqnarray*}
           \begin{eqnarray*}
    I_4&\leq&C \|\partial_t
    e_\theta\|_{L^3}\|\sqrt{\rho}\theta_t\|_{L^2}\|\rho\|_{L^\infty}^\frac{1}{2}
    \|\theta_t\|_{D^1_0}\\
            &\leq&
            Cc_0c_5^2(c_3^2+c_3\|\eta_t\|_{D^1_0})\|\sqrt{\rho}\theta_t\|_{L^2}^2
            +\frac{\kappa^0}{40}\|\nabla \theta_t\|_{L^2}^2,
            \end{eqnarray*}
    \begin{eqnarray*}
      I_5&\leq& C\int \rho |v||\nabla
      e_\theta||\theta_t|^2+\rho|v||e_\theta||\theta_t||\nabla\theta_t|dx\\
            &\leq&C\left(
            \|\rho\|_{L^\infty}^\frac{1}{2}\|v\|_{L^\infty}\|\nabla
            e_\theta\|_{L^3}+\|\rho\|_{L^\infty}^\frac{1}{2}\|v\|_{L^\infty}\|e_\theta\|_{L^\infty}
            \right)\|\sqrt{\rho}\theta_t\|_{L^2}\|\theta_t\|_{D^1_0}\\
      &\leq&
      Cc_5^2c_0^3\|\sqrt{\rho}\theta_t\|_{L^2}^2+\frac{\kappa^0}{40}\|\theta_t\|_{D^1_0}^2,
    \end{eqnarray*}
        \begin{eqnarray*}
          I_6&\leq& C\left(
          \|\partial_t e_\theta\|_{L^3}\|\rho\|_{L^\infty}\|v\|_{L^\infty}
          +\|e_\theta\|_{L^\infty}\|\rho_t\|_{L^3}\|v\|_{L^\infty}
          \right)\|\nabla\theta\|_{L^2}\|\theta_t\|_{D^1_0}\\
          &&+C\|e_\theta\|_{L^\infty}\|\rho\|_{L^\infty}^\frac{1}{2}\|\sqrt{\rho}\theta_t\|_{L^2}
          \|v_t\|_{D^1_0}\|\nabla \theta\|_{L^3}+\|e_\theta\|_{L^\infty}
          \|\rho\|_{L^\infty}^\frac{1}{2}\|v\|_{L^\infty}\|\sqrt{\rho}\theta_t\|_{L^2}
          \|\nabla \theta_t\|_{L^2}\\
                &\leq&
                Cc_5^2(c_3^2+c_3\|\eta_t\|_{D^1_0})c_0^4\|\nabla
                \theta\|_{L^2}^2+Cc_0^3c_5^2\|\sqrt{\rho}\theta_t\|_{L^2}^2
                +\frac{\kappa^0}{40}\|\nabla \theta_t\|_{L^2}^2\\
                &&+
                \|v_t\|_{D^1_0}^2+Cc_0c_5^2\|\nabla \theta\|_{L^2}
                \|\theta\|_{D^2}\|\sqrt{\rho}\theta_t\|_{L^2}^2,
        \end{eqnarray*}
    $$
    I_7\leq Cc_5\|\nabla v\|_{L^3}\|\nabla v_t\|_{L^2}\|\theta_t\|_{L^6}
            \leq Cc_0^2c_5^2\|\nabla v_t\|^2_{L^2}+\frac{\kappa^0}{20}
            \|\nabla \theta_t\|^2_{L^2},
    $$
         \begin{eqnarray*}
        I_8&\leq& C(\|\mu_t\|_{L^3}+\|\lambda_t\|_{L^3})\|(\nabla
        v)^2\|_{L^2}\|\theta_t\|_{L^6}\\
        &\leq& Cc_5(c_3+c_3^\frac{1}{2}\|\eta_t\|^\frac{1}{2}_{D^1_0})c_0^2\|
        \nabla \theta_t\|_{L^2}\\
        &\leq& Cc_5^2c_0^4(c_3^2+c_3\|\eta_t\|_{D^1_0})+
        \frac{\kappa^0}{20}\|\nabla
    \theta_t\|^2_{L^2},
        \end{eqnarray*}
    $$
    I_9\leq C\|\rho_t\|_{L^3}\|h\|_{L^2}\|\theta_t\|_{L^6}
    \leq Cc_0^2\|h\|_{L^2}\|\nabla \theta_t\|_{L^2}
    \leq Cc_0^4\|h\|^2_{L^2}+\frac{\kappa^0}{20}\|\nabla
    \theta_t\|^2_{L^2},
    $$
               \begin{eqnarray*}
              I_{10}&\leq& C\|h_t\|_{H^{-1}}\|\rho \theta_t\|_{H^1}\\
              &\leq& C\|h_t\|_{H^{-1}}(\|\rho\|^{\frac{1}{2}}_{L^\infty}\|
              \rho^{\frac{1}{2}}\theta_t\|_{L^2}+\|\rho\|_{L^\infty}
              \|\nabla \theta_t\|_{L^2}+\|\nabla \rho\|_{L^3}\|\nabla \theta_t\|_{L^2})
              \\
                    &\leq&
                    Cc_0^2\|h_t\|^2_{H^{-1}}+\|\rho^{\frac{1}{2}}\theta_t\|^2_{L^2}
                    +\frac{\kappa^0}{20}\|\nabla
                    \theta_t\|^2_{L^2},
            \end{eqnarray*}
and
    \begin{eqnarray*}
      I_{11}&\leq&\|\partial_te_\rho\|_{L^3}\|\rho\|_{L^\infty}^2\|\nabla
      v\|_{L^2}\|\theta_t\|_{D^1_0}+\|e_\rho\|_{L^\infty}\|\rho\|_{L^\infty}
      \|\rho_t\|_{L^3}\|\nabla v\|_{L^2}\|\theta_t\|_{D^1_0}\\
      &&    +\|e_\rho\|_{L^\infty}\|\rho\|_{L^\infty}^\frac{3}{2}\|\nabla
      v_t\|_{L^2}\|\sqrt{\rho}\theta_t\|_{L^2}\\
        &\leq&
        Cc_5^2(c_3^2+c_3\|\eta_t\|_{D^1_0})c_0^6+\frac{\kappa^0}{40}\|\theta_t\|_{D^1_0}^2
        +\|\nabla
        v_t\|_{L^2}^2+Cc_0^3c_5^2\|\sqrt{\rho}\theta_t\|_{L^2}^2.
    \end{eqnarray*}
Substituting these estimates into (\ref{2-b-E2.28}), we have
    \begin{eqnarray*}
      &&\frac{d}{dt}\|\sqrt{e_\theta\rho}\theta_t\|^2_{L^2}
      +\kappa^0\|\nabla \theta_t\|^2_{L^2}\\
        &\leq& Cc_2^6c_5^2(1+\|\nabla
        \theta\|^4_{L^2}+\|\theta\|_{D^2}^4)(1+\|\rho^{\frac{1}{2}}\theta_t\|_{L^2}^2)
        +Cc_5^2c_3^2\|\theta\|_{D^2}^2+Cc_3^{-1}\|\eta_t\|_{D^1_0}^2\\
            &&      +Cc_3^3c_5^4\|\theta\|_{D^2}^4+Cc_0^2c_5^2\|v_t\|_{D^1_0}^2
            +Cc_0^4c_5^2(c_3^2+c_3\|\eta_t\|_{D^1_0})(\|\sqrt{\rho}\theta_t\|_{L^2}^2+
            \|\nabla \theta\|_{L^2}^2)\\
        &&+Cc_0^6c_5^2(c_3^2+c_3\|\eta_t\|_{D^1_0})+Cc_0^4\|h\|_{L^2}^2+Cc_0^2\|h_t\|_{H^{-1}}^2\\
        &\leq&Cc_2^8c_3^3c_5^4(1+\|\nabla
        \theta\|^4_{L^2}+\|\theta\|_{D^2}^4)(1+\|\rho^{\frac{1}{2}}\theta_t\|_{L^2}^4)
        +Cc_3^{-1}\|\eta_t\|_{D^1_0}^2+Cc_0^2c_5^2\|v_t\|_{D^1_0}^2\\
        &&+Cc_0^{12}c_3^3c_5^4
        +Cc_0^4\|h\|_{L^2}^2+Cc_0^2\|h_t\|_{H^{-1}}^2,
    \end{eqnarray*}
and then integrating this over $(\tau,t)$, using the estimates
(\ref{2-b-E2.21-2}), (\ref{2-b-E2.24}) and (\ref{2-b-E2.24.2}), we
 have
    \begin{eqnarray}
      &&\|\sqrt{e_\theta\rho}\theta_t(t)\|^2_{L^2}
      +\int^t_\tau\|\nabla \theta_t\|^2_{L^2}ds\nonumber\\
        &\leq& C\|(\sqrt{e_\theta\rho}\theta_t)(\tau)\|^2_{L^2}
        +Cc_0^2c_1c_5^2+Cc_0^{12}c_3^3c_5^4t\nonumber\\
        &&+Cc_2^8c_3^3c_5^4\int^t_\tau(1+\|\nabla
        \theta\|^4_{L^2}+\|\theta\|_{D^2}^4)(1+\|\rho^{\frac{1}{2}}\theta_t\|_{L^2}^4)ds,
            \label{2-b-E2.29}
    \end{eqnarray}
for $0<\tau<t\leq \min(T,T_1)$. To estimate
$\limsup_{\tau\rightarrow0}\|\sqrt{e_\theta\rho}\theta_t(\tau)\|^2_{L^2}$,
we observe from the equation (\ref{2-b-E2.2}) that
    \begin{eqnarray*}
    \int e_\theta\rho \theta_t^2dx&\leq& C\int \left(
    \rho |h|^2+e_\theta\rho |v|^2|\nabla \theta|^2+\rho^{-1}p^2|\nabla
    v|^2\right.\\
            &&\left.
    +\rho^{-1}(\mathrm{div}(\kappa \nabla \theta)+Q(\sigma,\eta,v))^2
    +e_\rho^2\rho^3|\nabla v|^2
    \right)dx
    \end{eqnarray*}
and using the estimates (\ref{2-b-E2.20-1}) and
(\ref{2-b-E2.26-0}), we obtain
    \begin{eqnarray*}
    &&\limsup_{\tau\rightarrow0}\|\sqrt{e_\theta\rho}\theta_t(\tau)\|^2_{L^2}\\
    &\leq&
    C(\|\rho_0\|_{L^\infty}\|h(0)\|_{L^2}^2+\|e_\theta(\sigma(0),\eta(0))\|_{L^\infty}
    \|\rho_0\|_{L^\infty}\|\nabla u_0\|^2_{H^1}\|\nabla \theta_0\|_{L^2}^2
    +\|g_1\|^2_{L^2}\\&&
    +\|p_0^2(\rho_0)(1+\theta_0^2)^2\|_{L^\infty}
    \|\nabla v(0)\|_{L^2}^2+\|e_\rho(\sigma(0),\eta(0))\|_{L^\infty}^2
    \|\rho_0\|_{L^\infty}^3\|\nabla v(0)\|_{L^2}^2)\leq Cc_0^5c_5^2.
    \end{eqnarray*}
Hence, letting $\tau\rightarrow0$ in (\ref{2-b-E2.29}), we deduce
that
    \begin{eqnarray}
      &&\|\rho^{\frac{1}{2}}\theta_t(t)\|^2_{L^2}
      +\int^t_0\|\nabla \theta_t\|^2_{L^2}ds\nonumber\\
        &\leq&Cc_2^8c_3^3c_5^4\int^t_0(1+\|\nabla
        \theta\|^4_{L^2}+\|\theta\|_{D^2}^4)(1+\|\rho^{\frac{1}{2}}\theta_t\|_{L^2}^4)ds
        +Cc_0^4c_1c_5^2,
            \label{2-b-E2.30}
    \end{eqnarray}
for $0\leq t\leq \min(T,T_2)$, where
$T_2=\min(T_1,c_0^{-8}c_3^{-3}c_5^{-2})$. On the other hand, since
    $$
    \frac{d}{dt}\|\nabla \theta\|^2_{L^2}=2\int\nabla \theta\cdot\nabla
    \theta_tdx\leq 2\|\nabla \theta\|_{L^2}\|\nabla \theta_t\|_{L^2},
    $$
it follows that
    $$
    \|\nabla \theta(t)\|_{L^2}^2\leq C\|\nabla \theta_0\|^2_{L^2}
    +C\int^t_0\|\nabla \theta_t\|^2_{L^2}ds,
    $$
for $0\leq t\leq \min(T,T_2)$. Combining this and
(\ref{2-b-E2.30}), we obtain that
    \begin{eqnarray}
      &&\|\nabla \theta(t)\|_{L^2}^2+\|\rho^{\frac{1}{2}}\theta_t(t)\|^2_{L^2}
      +\int^t_0\|\nabla \theta_t\|^2_{L^2}ds\nonumber\\
        &\leq&Cc_2^8c_3^3c_5^4\int^t_0(1+\|\nabla
        \theta\|^4_{L^2}+\|\theta\|_{D^2}^4)(1+\|\rho^{\frac{1}{2}}\theta_t\|_{L^2}^4)ds
        +Cc_0^4c_1c_5^2,
            \label{2-b-E2.31}
    \end{eqnarray}
for $0\leq t\leq \min(T,T_2)$.

To estimate $\|\theta\|_{D^2}$, we use the following elliptic
regularity results: if $(f,g)\in D^1_0\cap D^{1,q}$ is a solution
to the elliptic system
    $$
    -\mathrm{div}(\kappa\nabla f)=F\textrm{ and }-\mathrm{div}(\mu(\nabla g+\nabla g^\top))
      -\nabla(\lambda \mathrm{div}g)=G,
    \textrm{ in }\Omega
    $$
with $(F,G)\in L^2\cap L^q$. From
(\ref{2-b-E2.26-0})-(\ref{2-b-E2.26-1}), we have
      $$
    \sup_{0\leq t\leq \min(T,T_2)}
    \|(\kappa,\mu,\lambda)(\sigma(t,\cdot),\eta(t,\cdot))\|_{C^{\beta}\cap D^{1,q}}\leq c_5,
    \ \beta=1-\frac{3}{q}.
      $$
From Lemma \ref{2-b-L5.5} and \ref{2-b-L5.7}, we have $(f,g)\in
D^{2}$ satisfying
    \begin{equation}
    \|f\|_{D^{2}(\Omega)}\leq c_6\left(\|F\|_{L^{2}(\Omega)}
      +\|\nabla f\|_{L^2(\Omega)}+\|\nabla \kappa
      \|_{L^q(\Omega)}^\frac{q}{q-3}\| f\|_{D^1_0(\Omega)}
      \right),
      \label{2-b-E2.33-1}
    \end{equation}
        \begin{equation}
            \|g\|_{D^{2}(\Omega)}\leq c_6\left(\|G\|_{L^{2}(\Omega)}
      +\|\nabla g\|_{L^2(\Omega)}+\||\nabla \mu|+|\nabla \lambda|
      \|_{L^q(\Omega)}^\frac{q}{q-3}\| g\|_{D^1_0(\Omega)}
      \right),
      \label{2-b-E2.33-2}
    \end{equation}
and  from Lemma \ref{2-b-L5.6} and \ref{2-b-L5.8}, we have
$(f,g)\in D^{2,q}$ satisfying
    \begin{equation}
    \|f\|_{D^{2,q}(\Omega)}\leq c_6\left(\|F\|_{L^{q}(\Omega)}
      +\|\nabla f\|_{L^q(\Omega)}+\|\nabla \kappa
      \|_{L^q(\Omega)}^{\frac{5q-6}{2q-6}}\| f\|_{D^1_0(\Omega)}
      \right),
      \label{2-b-E2.33-1-1}
    \end{equation}
        \begin{equation}
            \|g\|_{D^{2,q}(\Omega)}\leq c_6\left(\|G\|_{L^{q}(\Omega)}
      +\|\nabla g\|_{L^q(\Omega)}+\||\nabla \mu|+|\nabla \lambda|
      \|_{L^q(\Omega)}^{\frac{5q-6}{2q-6}}\| g\|_{D^1_0(\Omega)}
      \right),
      \label{2-b-E2.33-2-1}
    \end{equation}
where the constant $c_6=c_6(c_0)>1$ depends only on  $c_0$.
     It should be noted that the estimates
(\ref{2-b-E2.33-1})-(\ref{2-b-E2.33-2-1}) hold for both bounded
and unbounded domains. Applying the elliptic regularity result
(\ref{2-b-E2.33-1}) to the equation $-\mathrm{div}(\kappa\nabla
e)=F$ in $\Omega$, where $F=\rho(h-e_\theta\theta_t-e_\theta
v\cdot\nabla \theta)-p\mathrm{div}v+Q(\sigma,\eta,v)+e_\rho
\rho^2\mathrm{div}v$, using H\"{o}lder and Sobolev's inequalities,
we have
    \begin{eqnarray}
      \|\theta\|_{D^2}&\leq&Cc_6(\|\rho h\|_{L^2}+
      \|e_\theta\|_{L^\infty}\|\rho \theta_t\|_{L^2}+
      \|e_\theta\|_{L^\infty}\|\rho v\cdot\nabla \theta\|_{L^2}+\|p\mathrm{div}v\|_{L^2}
      \nonumber\\
        &&  +\|Q(\sigma,\eta,v)\|_{L^2}+\|\nabla \theta\|_{L^2}+\|\nabla \kappa
      \|_{L^q(\Omega)}^\frac{q}{q-3}\| \theta\|_{D^1(\Omega)}+
      \|e_\rho\|_{L^\infty}\|\rho\|_{L^\infty}^2\|\nabla v\|_{L^2})\nonumber\\
              &\leq& Cc_0^3c_5c_6(1+\|\rho^{\frac{1}{2}}\theta_t\|_{L^2}+\|\nabla \theta\|_{L^2})
              +Cc_6 c_5^{\frac{q}{q-3}}
            \|\nabla \theta\|_{L^2}\nonumber\\
      &&+Cc_6\|p_0\sqrt{\rho}\|_{L^\infty}\|(1+\theta^2)\nabla v\|_{L^2}
            \nonumber\\
    &\leq& Cc_2^3c_5^{\frac{q}{q-3}}c_6
    (1+\|\rho^{\frac{1}{2}}\theta_t\|_{L^2}+\|\nabla \theta\|_{L^2}^2).
    \label{2-b-E2.41}
    \end{eqnarray}
Combining this and (\ref{2-b-E2.31}), we deduce that
    \begin{eqnarray}
      &&\|\nabla \theta(t)\|_{L^2}^2+\|\rho^{\frac{1}{2}}\theta_t(t)\|^2_{L^2}
      +\int^t_0\|\nabla \theta_t\|^2_{L^2}ds\nonumber\\
        &\leq&Cc_2^{20}c_3^3c_5^{\frac{6q}{q-3}}c_6^4\int^t_0[(1+\|\nabla
        \theta\|^8_{L^2}+\|\rho^{\frac{1}{2}}\theta_t\|_{L^2}^4)
        (1+\|\rho^{\frac{1}{2}}\theta_t\|_{L^2}^4)ds
        +Cc_0^4c_1c_5^2,
            \label{2-b-E2.31-1}
    \end{eqnarray}
for $0\leq t\leq \min(T,T_2)$. Now we define a function $H(t)$ by
    $$
    H(t)=1+\|\rho^{\frac{1}{2}}\theta_t\|^2_{L^2}+\|\nabla
    \theta\|_{L^2}^2.
    $$
Then it follows from (\ref{2-b-E2.31-1}) that
    $$
    H(t)\leq Cc_0^4c_1c_5^2+Cc_2^{20}c_3^3c_5^{\frac{6q}{q-3}}c_6^4\int^t_0
    H^6(s)ds\textrm{ for }0\leq t\leq \min(T,T_2).
    $$
Solving this integral inequality, we easily derive
    $$
    H(t)\leq Cc_0^4c_1c_5^2
    (1-Cc_2^{45}c_3^3c_5^{\frac{11q}{q-3}}c_6^4t)^{-\frac{1}{5}}
    \textrm{ for all small }t\geq0.
    $$
Therefore, taking $T_3=\min(T_2,
(2Cc_2^{45}c_3^3c_5^{\frac{11q}{q-3}}c_6^4)^{-1})$, we conclude
that
    \begin{equation}
      \|\theta(t)\|^2_{D_0^1}+\|\rho^{\frac{1}{2}}\theta_t(t)\|^2_{L^2}
      +\int^t_0\|\theta_t(s)\|^2_{D^1_0}ds\leq A_1c_0^4c_1c_5^2,
      \label{2-b-E2.32}
    \end{equation}
and
    \begin{equation}
      \|\theta(t)\|_{D^2}\leq A_1 c_2^8c_5^{\frac{2q}{q-3}}c_6,
      \label{2-b-E2.34-1}
    \end{equation}
for $0\leq t\leq \min(T,T_3)$. Here, the constant $A_1>1$ is
independent of any one of $c_i$, $i=0,\ldots,6$. From
(\ref{2-b-E2.32})-(\ref{2-b-E2.34-1}), using the
Gagliardo-Nirenberg inequality, we have
    \begin{equation}
      \|\theta(t)-\theta_0\|_{D^1_0} \leq C
      \int^t_0\|\theta_t\|_{D^1_0}ds\leq
      C\left(\int^t_0\|\theta_t\|_{D^1_0}^2ds
      \right)^{\frac{1}{2}}t^{\frac{1}{2}}
             \leq  Cc_0^2c_1^{\frac{1}{2}}c_5t^{\frac{1}{2}}\leq
            Cc_0,\label{2-b-E2.26-101}
    \end{equation}
    \begin{eqnarray}
    \|\theta(t)-\theta_0\|_{D^{1,q}}&\leq&
    C\|\theta(t)-\theta_0\|_{D^1_0}^{\frac{6-q}{2q}}
    \|\theta(t)-\theta_0\|_{D^2}^{\frac{3q-6}{2q}}\nonumber\\
            &\leq&Cc_2^8c_5^{\frac{2q}{q-3}}c_6t^{\frac{6-q}{4q}}
            \leq Cc_0,\label{2-b-E2.26-100}
    \end{eqnarray}
for $0\leq t\leq \min(T,T_4)$, where
$T_4=\min(T_3,c_2^{-\frac{32q}{6-q}}c_5^{-\frac{8q^2}{(q-3)(6-q)}}c_6^{-\frac{4q}{6-q}})$.
Therefor, from (\ref{2-b-E1.7})-(\ref{2-b-E1.8}),
(\ref{2-b-E2.32})  and
(\ref{2-b-E2.26-101})-(\ref{2-b-E2.26-100}), using H\"{o}lder and
Sobolev's inequalities, we obtain
    \begin{equation}
      \sup_{0\leq t\leq \min(T,T_4)}
      \left(\|\theta(t,\cdot)\|_{D^1_0\cap
      D^{1,q}}+\|\theta(t,\cdot)\|_{L^\infty}
      \right)\leq Cc_0,
    \end{equation}
        \begin{equation}
          \sup_{0\leq t\leq \min(T,T_4)}\|p(\rho(t,\cdot),\theta(t,\cdot))\|_{L^\infty}
         \leq Cc_0^3c_5,
        \end{equation}
    \begin{eqnarray}
        &&\sup_{0\leq t\leq \min(T,T_4)}\|\nabla p(\rho(t,\cdot),\theta(t,\cdot))\|_{L^2\cap L^q}
        \nonumber\\
         &\leq&\sup_{0\leq t\leq \min(T,T_4)}\left(
         \|p_\rho\nabla \rho\|_{L^2\cap L^q}+\|p_\theta\nabla \theta\|_{L^2\cap L^q}
         \right)\nonumber\\
                &\leq&\sup_{0\leq t\leq \min(T,T_4)}\left(
         \|p_2(1+\theta^2)\nabla \rho\|_{L^2\cap L^q}+\|p_1\sqrt{\rho}(1+\theta)\nabla \theta\|_{L^2\cap L^q}
         \right)\nonumber\\
         &\leq& Cc_0^3c_5,\label{2-b-E2.35}
    \end{eqnarray}
and
    \begin{eqnarray}
        &&\sup_{0\leq t\leq \min(T,T_4)}\|\partial_t
         p(\rho(t,\cdot),\theta(t,\cdot))\|_{L^2}
         \leq\sup_{0\leq t\leq \min(T,T_4)}\left(
         \|p_\rho \rho_t\|_{L^2}+\|p_\theta\theta_t\|_{L^2}
         \right)\nonumber\\
                &\leq&\sup_{0\leq t\leq \min(T,T_4)}\left(
         \|p_2(1+\theta^2)\rho_t\|_{L^2}+\|p_1\sqrt{\rho}(1+\theta) \theta_t\|_{L^2}
         \right)\leq Cc_1^4c_5^2.\label{2-b-E2.35-1}
    \end{eqnarray}

Hence by virtue of (\ref{2-b-E2.33-1-1}), from
(\ref{2-b-E2.21-2}), (\ref{2-b-E2.25}), (\ref{2-b-E2.24.3}) and
(\ref{2-b-E2.26-0})-(\ref{2-b-E2.26-1}), using H\"{o}lder's
inequality, we deduce that
    \begin{eqnarray}
      &&\int^t_0\|\nabla \theta(s)\|^2_{D^{1,q}}ds\nonumber\\
            &\leq& Cc_6^2\int^t_0(\|\rho h\|_{L^q}^2+
            \|e_\theta\|^2_{L^\infty}\left(\|\rho \theta_t\|^2_{L^q}+
            \|\rho v\cdot\nabla \theta\|^2_{L^q}\right)
            +\|p\mathrm{div}v\|^2_{L^q}\nonumber\\
                &&+\|Q(\sigma,\eta, v)\|^2_{L^q}
                +\|e_\rho\|_{L^\infty}^2\|\rho\|_{L^\infty}^4\|\nabla
                v\|_{L^q}^2
                +\|\nabla \theta\|^2_{L^q}+
                \|\nabla \kappa
                \|_{L^q(\Omega)}^{\frac{5q-6}{q-3}}\| \theta\|_{D^1_0(\Omega)}^2)ds\nonumber\\
      &\leq&Cc_6^2\int^t_0[c_0^2\|h\|^2_{L^q}+c_0^2c_5^2(\|\rho^{\frac{1}{2}}\theta_t\|^2_{L^2}
      +\|\theta_t\|_{D^1_0}^2)
      +c_0^8c_5^2+c_5^2\| \nabla v\|^2_{W^{1,q}}+c_0^2c_5^{\frac{5q-6}{q-3}}]ds\nonumber\\
            &\leq& Cc_0^6c_1c_5^4c_6^2,
            \label{2-b-E2.34-2}
    \end{eqnarray}
for $0\leq t\leq \min(T,T_4)$.

\textbf{III. Estimates for the velocity.}

To derive estimates for the velocity $u$, we differentiate the
equation (\ref{2-b-E2.3}) with respect to $t$, multiply it by
$u_t$ and integrate it over $\Omega$,  using the equation
(\ref{2-b-E2.1}), then derive
    \begin{eqnarray}
      &&\frac{1}{2}\frac{d}{dt}\int\rho|u_t|^2dx+
      \int(\mu[\nabla u_t+\nabla u_t^\top]:\nabla u_t
      +\lambda(\mathrm{div}u_t)^2)dx\nonumber\\
            &=&\int u_t\cdot[-\nabla p_t+(\rho
            f)_t-\rho_tv\cdot\nabla u-\rho
            v_t\cdot\nabla u+\mathrm{div}(\rho v)u_t\nonumber\\
                &&+\mathrm{div}(\mu_t(\nabla u+\nabla u^\top)+\nabla(
                \lambda_t\mathrm{div}u))    ]dx.\label{2-b-E2.55}
    \end{eqnarray}
Using the technique in \cite{Feireisl04}, we can estimate the
second term in the left hand side of (\ref{2-b-E2.55}) as follows:
    \begin{eqnarray}
      &&\int(\mu[\nabla u_t+\nabla u_t^\top]:\nabla u_t
      +\lambda(\mathrm{div}u_t)^2)dx\nonumber\\
            &\geq& \int(\mu[\nabla u_t+\nabla u_t^\top]:\nabla u_t
              -\frac{2}{3}\mu(\mathrm{div}u_t)^2)dx\nonumber\\
    &=& \int\frac{\mu}{2}|\nabla u_t+\nabla u_t^\top-\frac{2}{3}
    \mathrm{div}u_tI|^2dx\nonumber\\
            &\geq& \frac{\mu^0}{2}\int|\nabla u_t+\nabla u_t^\top-\frac{2}{3}
            \mathrm{div}u_tI|^2dx\nonumber\\
    &=& \int(\mu^0[\nabla u_t+\nabla u_t^\top]:\nabla u_t
    -\frac{2}{3}\mu^0(\mathrm{div}u_t)^2)dx\nonumber\\
            &=& \int(\mu^0|\nabla u_t|^2
            +\frac{1}{3}\mu^0(\mathrm{div}u_t)^2)dx\geq
            \int\mu^0|\nabla u_t|^2dx.\label{2-b-E2.55-1}
    \end{eqnarray}
Using H\"{o}lder, Young and Sobolev's inequalities, we easily
deduce that
    \begin{eqnarray*}
      &&\frac{d}{dt}\|\rho^{\frac{1}{2}}u_t\|^2_{L^2}+
      \mu^0\|\nabla u_t\|^2_{L^2}\\
            &\leq &C(\|p_t\|^2_{L^2}+\|\rho\|_{L^\infty}\|
            \rho^{\frac{1}{2}}u_t\|^2_{L^2}+(1+\|\rho\|^2_{L^\infty}
            +\|\nabla \rho\|^2_{L^3})\|f_t\|^2_{H^{-1}}+
            \|\rho_t\|^2_{L^3}\|f\|^2_{L^2}\\
      &&+\|\rho_t\|^2_{L^3}\|v\|^2_{L^\infty}\|\nabla u\|^2_{L^2}+\|\rho\|_{L^\infty}\|v\|_{L^\infty}^2
      \|\rho^{\frac{1}{2}}u_t\|^2_{L^2})+
        Cc_1\|\rho\|_{L^\infty}\|\nabla u\|_{L^2}\|\nabla
      u\|_{D^1}\\
            &&+c_1^{-1}\|\nabla
            v_t\|^2_{L^2}\|\rho^{\frac{1}{2}}u_t\|_{L^2}^2
            +C(\|\mu_t\|^2_{L^3}+\|\lambda_t\|^2_{L^3})\|u\|^2_{D^2}.
    \end{eqnarray*}
Then it follows from the estimates (\ref{2-b-E2.21-2}),
(\ref{2-b-E2.25}), (\ref{2-b-E2.24.3})-(\ref{2-b-E2.25.1}),
(\ref{2-b-E2.26-2}) and (\ref{2-b-E2.35-1}) that
    \begin{eqnarray*}
      &&\frac{d}{dt}\|\rho^{\frac{1}{2}}u_t\|^2_{L^2}+\|\nabla
      u_t\|^2_{L^2}\\
    &\leq& Cc_1^{8}c_5^4+Cc_0^2\|f_t\|^2_{H^{-1}}+(Cc_0^3+c_1^{-1}\|\nabla v_t\|^2_{L^2})
    \|\rho^{\frac{1}{2}}u_t\|^2_{L^2}\\
        &&+C(c_0^6+c_0^2c_1^2)\|
        \nabla u\|^2_{L^2}+Cc_5^2(c_3^2+c_3
        \|\eta_t(t,\cdot)\|_{D^1_0})\|u\|^2_{D^2}.
    \end{eqnarray*}
for $0\leq t\leq \min(T,T_4)$. Hence, using the facts that
    $$
    \rho(\tau)^{-\frac{1}{2}}\left[L(\sigma(\tau),\eta(\tau),u(\tau))+\nabla
    p(\rho(\tau),\theta(\tau))\right]
    \rightarrow g_2
    \textrm{ in }L^2
    \textrm{ as }\tau\rightarrow0
    $$
and
    $$
    \|u(t)\|^2_{D^1_0}\leq C\|u_0\|^2_{D^1_0}+C\int^t_0\|\nabla
    u_t\|^2_{L^2}ds
    \ \textrm{ for }\ 0\leq t\leq \min(T,T_4),
    $$
we easily obtain
    \begin{eqnarray}
      &&\|u(t)\|^2_{D^1_0}+\|\rho^\frac{1}{2}u_t\|^2_{L^2}+\int^t_0\|
      \nabla u_t\|^2_{L^2}ds\nonumber\\
            &\leq& Cc_0^5+C\int^t_0(c_0^4c_1^2+
            c_1^{-1}\|\nabla v_t\|^2_{L^2})(\|u\|^2_{D^1_0}+\|\rho^{\frac{1}{2}}
            u_t\|^2_{L^2})ds\nonumber\\
            &&+Cc_5^2\int^t_0(c_3^2+c_3
        \|\eta_t(t,\cdot)\|_{D^1_0})\|u\|^2_{D^2}ds,
            \label{2-b-E2.37-1}
    \end{eqnarray}
for $0\leq t\leq \min(T,T_4)$.

To estimate $\|u\|_{D^2}$, we apply the elliptic regularity result
(\ref{2-b-E2.33-2}) to the equation (\ref{2-b-E2.3}). Then, it
follows from (\ref{2-b-E2.21-2}), (\ref{2-b-E2.25}),
(\ref{2-b-E2.24.3})-(\ref{2-b-E2.26-1}) and (\ref{2-b-E2.35}) that
    \begin{eqnarray}
      \|u\|_{D^2}&\leq& Cc_6(\|\rho f\|_{L^2}+\|\rho u_t\|_{L^2}+\|
      \rho v\cdot\nabla u\|_{L^2} +\|\nabla p\|_{L^2}\nonumber\\
            &&+\|\nabla u\|_{L^2}
            +\||\nabla\mu|+|\nabla\lambda|\|^{\frac{q}{q-3}}_{L^q}
            \|\nabla u\|_{L^2})\nonumber\\
                &\leq& Cc_6[\|\rho\|_{L^\infty}\|f\|_{L^2}
                +\|\rho\|_{L^\infty}^{\frac{1}{2}}\|\rho^{\frac{1}{2}}u_t\|_{L^2}
              +\|\rho\|_{L^\infty}\|v\|_{L^\infty}\|\nabla u\|_{L^2}\nonumber\\
      &&+c_0^3c_5+
              \|u\|_{D^1_0}+c_5^\frac{q}{q-3}\|\nabla u\|_{L^2}]\nonumber\\
      &\leq& Cc_6[c_0+c_0^{\frac{1}{2}}\|\rho^{\frac{1}{2}}u_t\|_{L^2}
      +c_0^3c_5+(c_0^2+c_5^{\frac{q}{q-3}})
      \|u\|_{D^1_0}]\nonumber\\
      &\leq& Cc_0^3c_5c_6+
        Cc_0^2c_5^{\frac{q}{q-3}}c_6
        (\|u\|_{D^1_0}+\|\rho^{\frac{1}{2}}u_t\|_{L^2}),
        \label{2-b-E2.48}
    \end{eqnarray}
for $0\leq t\leq \min(T,T_4)$. Therefore, substituting this into
(\ref{2-b-E2.37-1}), we obtain
    \begin{eqnarray*}
      &&\|u(t)\|^2_{D_0^1}+\|\rho^\frac{1}{2}u_t\|^2_{L^2}+\int^t_0\|
      \nabla u_t\|^2_{L^2}ds\nonumber\\
            &\leq& Cc_0^5 +C\int^t_0(c_0^4c_5^{\frac{3q}{q-3}}c_6^2( c_1^2c_3^2
            + c_3  \|\eta_t(t,\cdot)\|_{D^1_0})+
            c_1^{-1}\|\nabla v_t\|^2_{L^2})(\|u\|^2_{D^1_0}+\|\rho^{\frac{1}{2}}
            u_t\|^2_{L^2})ds.
    \end{eqnarray*}
and using Gronwall's inequality, we obtain that
    \begin{equation}
       \|u(t)\|^2_{D^1_0}+\|\rho^\frac{1}{2}u_t(t)\|^2_{L^2}
      +\int^t_0\|u_t(s)\|^2_{D^1_0}ds
         \leq  Cc_0^5\exp\{C+Cc_0^4
        c_3^\frac{3}{2}c_5^{\frac{3q}{q-3}}c_6^2t^\frac{1}{2}\}
        \leq Cc_0^5,
      \label{2-b-E2.38}
    \end{equation}
and
    \begin{equation}
    \|u(t)\|_{D^2}\leq
    Cc_0^5c_5^{\frac{q}{q-3}}c_6,
    \label{2-b-E2.38-1}
    \end{equation}
for $0\leq t\leq \min(T,T_4)$.

Moreover, it follows from (\ref{2-b-E2.33-2-1}) that if $0\leq
t\leq \min(T,T_4)$, then
    \begin{eqnarray*}
       \|u\|_{D^{2,q}}
      &\leq&
      Cc_6(c_0\|f\|_{L^q}+c_0(\|\rho^{\frac{1}{2}}u_t\|_{L^2}+\|u_t\|_{D^1_0})
      +c_0^7c_5^{\frac{q}{q-3}}c_6+c_5^\frac{5q-6}{2q-6}\|\nabla
      u\|_{L^2})\\
              &\leq&Cc_6(c_0\|f\|_{L^q}+c_0^{\frac{7}{2}}+c_0\|u_t\|_{D^1_0}
              +c_0^7c_5^{\frac{q}{q-3}}c_6
              +c_5^{\frac{5q-6}{2q-6}}c_0^{\frac{5}{2}}),
    \end{eqnarray*}
and thus
    \begin{equation}
      \int^t_0\|u(s)\|^2_{D^{2,q}}ds\leq Cc_0^7c_6^2,
      \label{2-b-E2.39}
    \end{equation}
for $0\leq t\leq\min(T,T_4)$.

Combining (\ref{2-b-E2.38})-(\ref{2-b-E2.39}), we finally conclude
that
    \begin{equation}
      \|u(t)\|_{D^1_0}+\|\rho^\frac{1}{2}u_t(t)\|_{L^2}
      +\int^t_0\left(\|u_t(s)\|_{D^1_0}^2+\|u(s)\|^2_{D^{2,q}}
      \right)ds\leq A_2c_0^7c_6^2,\label{2-b-E2.40}
    \end{equation}
and
    \begin{equation}
      \|u(t,\cdot)\|_{D^2}\leq
      A_2c_0^5c_5^{\frac{q}{q-3}}c_6,\label{2-b-E2.40-1}
    \end{equation}
     for $0\leq t\leq \min(T,T_4)$. Here, the constant $A_2>1$  is independent of any one of $c_i$,
$i=0,\dots,6$.

\textbf{IV. Conclusion.}

Let us define $c_1$, $c_2$ and $c_3$ by
    \begin{equation}
      c_1=A_2c_0^7c_6^2,
      \   \ c_2=A_2c_0^5c_5^{\frac{q}{q-3}}c_6
      \textrm{ and }
      c_3=A_1 c_2^8c_5^{\frac{2q}{q-3}}c_6.
    \end{equation}
Then choosing any $T_{***}$ such that $0<T_{***}\leq \min(T,T_4)$,
we conclude from  (\ref{2-b-E2.25.1}),
(\ref{2-b-E2.32})-(\ref{2-b-E2.34-1}), (\ref{2-b-E2.34-2}) and
(\ref{2-b-E2.40})-(\ref{2-b-E2.40-1}) that
    \begin{equation}
      \sup_{0\leq t\leq T_{***}}\left(\|\rho(t,\cdot)-\rho^{\infty}\|_{L^6\cap D^1\cap
      D^{1,q}}^2
      +\|\rho_t(t,\cdot)\|_{L^2\cap L^q}
      \right)
      \leq c_4,\label{2-b-E2.42-1}
    \end{equation}
        \begin{equation}
        \sup_{t\in[0,T_{***}]}\left(\|(\sqrt{\rho}\theta_t)(t,\cdot)\|_{L^2}+\|\theta(t,\cdot)\|_{D^2}
         \right)+\int^{T_{***}}_0\|\theta_t(t,\cdot)\|_{D^1_0}^2dt\leq
          c_3,\label{2-b-E2.42-1.1}
        \end{equation}
    \begin{equation}
      \sup_{0\leq t\leq T_{***}}\left(\|\theta(t,\cdot)\|_{D_0^1}+
        \|\sqrt{\rho}u_t(t,\cdot)\|_{L^2}\right)
        +\int^{T_{***}}_0\|\theta(t,\cdot)\|_{D^{2,q}}^2
        dt\leq Cc_7,\label{2-b-E2.42-2}
    \end{equation}
    \begin{equation}
      \sup_{0\leq t\leq T_{***}}
      \|u(t,\cdot)\|_{D^1_0}
      +\int^{T_{***}}_0\left(\|u_t(s,\cdot)\|_{D^1_0}^2+\|u(s,\cdot)\|^2_{D^{2,q}}
      \right)ds\leq c_1,
    \end{equation}
    \begin{equation}
      \sup_{0\leq t\leq T_{***}}
     \|u(t,\cdot)\|_{D^2}
      ds\leq c_2,\label{2-b-E2.42-3}
    \end{equation}
where $c_7= c_0^6c_1c_5^4c_6^2$. Using  similar arguments in
(\ref{2-b-E2.23})-(\ref{2-b-E2.23.2}), we have
    \begin{equation}
    \|(\theta(t,\cdot)-\theta_0(\cdot),u(t,\cdot)-u_0(\cdot))\|_{D^1_0}
    \leq Cc_3^\frac{1}{2}t^{\frac{1}{2}},\label{2-b-E2.42-3.0}
    \end{equation}
and
        \begin{equation}
    \|\rho(t,\cdot)-\rho_0(\cdot)\|_{L^6}
    \leq Cc_0^{2}t^{\frac{1}{3}},\label{2-b-E2.42-3.1}
    \end{equation}
 for $0\leq t\leq T_{***}$.     Hence, it deserves to emphasize
that $c_1, \dots, c_7$ and $T_{***}$ depend only on $c_0$, but not
on the lower bound $\delta$ of the initial density $\rho_0$.

\section{Existence result for the linearized problem}\label{2-b-Sec2-1}\setcounter{equation}{0}
Now we can prove the key lemma to prove Theorem \ref{2-b-T1}.

\begin{lem}\label{2-b-L2.4}
  Assume in addition to (\ref{2-b-E2.21-1})-(\ref{2-b-E2.21-5})
  that the initial data $(\rho_0,\theta_0,u_0)$ satisfy the
  compatibility condition
        \begin{equation}
          -\mathrm{div} (\kappa_0\nabla \theta_0)-Q(\rho_0,\theta_0,u_0)
          =\rho_0^\frac{1}{2}g_1
          \textrm{ and }L(\rho_0,\theta_0,u_0)+\nabla p(\rho_0,\theta_0)=\rho_0^\frac{1}{2}
          g_2,\label{2-b-E2.43}
        \end{equation}
  for some $(g_1,g_2)\in L^2$.
  Assume further that
    $$
    (\sigma,\eta,v)(0)=(\rho_0,\theta_0,u_0),
    \ \sigma^\infty=\rho^\infty,
    $$
    $$
    \sup_{t\in[0,T_{***}]}\|v(t,\cdot)\|_{D^1_0}
    +\int^{T_{***}}_0\left(\|v_t(t,\cdot)\|_{D^1_0}^2+\|v(t,\cdot)\|^2_{D^{2,q}}
      \right)dt\leq c_1
    $$
        $$
    \sup_{t\in[0,T_{***}]}\|v(t,\cdot)\|_{D^2}\leq c_2,
        $$
        $$
        \sup_{t\in[0,T_{***}]}\left(
        \|(\sqrt{\sigma}\eta_t)(t,\cdot)\|_{L^2}+\|\eta(t,\cdot)\|_{D^2}
         \right)+\int^{T_{***}}_0\|\eta_t(t,\cdot)\|_{D^1_0}^2dt\leq
          c_3,
        $$
        $$
    \sup_{t\in[0,T_{***}]}\left(\|\sigma(t,\cdot)-\rho^\infty\|_{L^6\cap D^1\cap
    D^{1,q}}+\|\sigma_t(t,\cdot)\|_{L^2\cap L^q}\right)\leq c_4,
        $$
  for the positive constants $c_1$, $c_2$, $c_3$, $c_4$ and $T_{***}$,
  depend only on $c_0$, chosen as
  before, where
     $$
    c_0=2+\rho^\infty+\|\rho_0-\rho^\infty\|_{L^6\cap D^1\cap D^{1,q}}
    +\|(\theta_0,u_0)\|_{D^1_0\cap D^2}+\|(g_1,g_2)\|^2_{L^2}.
   $$
Then there exists a unique strong solution $(\rho,\theta,u)$ to
the linearized problem (\ref{2-b-E2.1})-(\ref{2-b-E2.6}) in
$[0,T_{***}]$ satisfying the estimates
(\ref{2-b-E2.42-1})-(\ref{2-b-E2.42-3.1}) as well as the
regularity
    $$
    \rho-\rho^\infty\in C([0,T_{***}];L^6\cap D^1\cap D^{1,q}),
    \ \rho_t\in C([0,T_{***}];L^2\cap L^q),
    $$
        $$
        (\theta,u)\in C([0,T_{***}];D_0^1\cap D^2)\cap L^2(0,T_{***};D^{2,q}),
        $$
    $$
    (\theta_t,u_t)\in L^2(0,T_{***};D^1_0)
    \textrm{ and }(\rho^\frac{1}{2}\theta_t,\rho^{\frac{1}{2}}u_t)
    \in L^\infty([0,T_{***}];L^2),
    $$
        $$
    \|\rho(t,\cdot)-\rho_0(\cdot)\|_{L^6\cap D^1\cap D^{1,q}}
    +\|(\theta(t,\cdot)-\theta_0(\cdot),u(t,\cdot)-u_0(\cdot))\|_{D^1_0\cap D^{2}}
    \rightarrow0,\ \textrm{ as }\ t\rightarrow0.
        $$
\end{lem}
\begin{proof}
  We define $\rho_0^\delta=\rho_0+\delta$ for each
  $\delta\in(0,1)$. Then from the compatibility condition
  (\ref{2-b-E2.43}), we derive
    $$
    \left\{\begin{array}{l}
    -\mathrm{div} (\kappa(\rho_0,\theta_0)\nabla
    \theta_0)-Q(\rho_0,\theta_0,
    u_0)
          =(\rho^\delta_0)^\frac{1}{2}g_1^\delta,\\
          L(\rho_0,\theta_0,u_0)+
          \nabla p(\rho^\delta_0, \theta_0) =(\rho^\delta_0)^\frac{1}{2}
          g_2^\delta,
    \end{array}\right.
    $$
  where
         $$
    g_1^\delta=\left(\frac{\rho_0}{\rho_0^\delta}
    \right)^{\frac{1}{2}}g_1 \ \textrm{ and }
    g_2^\delta=\left(\frac{\rho_0}{\rho_0^\delta}
    \right)^{\frac{1}{2}}g_2+\frac{\nabla(p(\rho_0^\delta,\theta_0)-
    p(\rho_0,\theta_0))}
    {(\rho_0^\delta)^\frac{1}{2}}.
        $$
  Since $p_\rho,p_\theta\in
  C^\alpha([0,\infty),L^\infty_{loc}([0,\infty)))$,
$\alpha\in(\frac{1}{2},1
  )$, we obtain that for all small $\delta>0$,
    $$
    c_0\geq1+(\rho^\infty+\delta)+\|\rho^\delta_0-(\rho^\infty+\delta)
    \|_{L^6\cap D^1\cap D^{1,q}}
    +\|(\theta_0,u_0)\|_{D^1_0\cap D^2}+\|(g_1^\delta,g_2^\delta)\|^2_{L^2}.
    $$
Hence, from  previous results for positive initial densities, we
deduce that corresponding to the initial data
$(\rho_0^\delta,\theta_0,u_0)$ with small $\delta>0$, there exists
a unique strong solution $(\rho^\delta,\theta^\delta,u^\delta)$ of
the linearized equations (\ref{2-b-E2.1})-(\ref{2-b-E2.6})
satisfying the local estimates
(\ref{2-b-E2.42-1})-(\ref{2-b-E2.42-3}). From this uniform
estimates on $\delta$, we conclude that there is a subsequence of
the solutions $(\rho^\delta,\theta^\delta,u^\delta)$ which
converges to a limit $(\rho,\theta,u)$ in an obvious weak or
weak$-*$ sense. It is  easy to show that $(\rho,\theta,u)$ is a
weak solution to the linearized equations
(\ref{2-b-E2.1})-(\ref{2-b-E2.6}) in $[0,T_{***}]$. Finally,
thanks to the lower semi-continuity of various norms, we find that
$(\rho,\theta,u)$ also satisfies the estimates
(\ref{2-b-E2.42-1})-(\ref{2-b-E2.42-3.1}). This proves the
existence of a strong solution $(\rho,\theta,u)$ with the
regularity
    \begin{equation}
    \rho-\rho^\infty\in L^\infty([0,T_{***}];L^6\cap D^1\cap D^{1,q}),
    \ \rho_t\in L^\infty([0,T_{***}];L^2\cap L^q),
    \label{2-b-E2.44-1}
    \end{equation}
    \begin{equation}
        (\theta,u)\in L^\infty([0,T_{***}];D_0^1\cap D^2)\cap L^2(0,T_{***};D^{2,q}),
            \label{2-b-E2.44-2}
    \end{equation}
    \begin{equation}
    (\theta_t,u_t)\in L^2(0,T_{***};D^1_0)
    \textrm{ and }(\rho^\frac{1}{2}\theta_t,\rho^{\frac{1}{2}}u_t)
    \in L^\infty(0,T_{***};L^2).
        \label{2-b-E2.44-3}
    \end{equation}

Now we prove the uniqueness of the solution to this regularity
class. Let $(\rho_1,\theta_1,u_1)$ and $(\rho_2,\theta_2,u_2)$ be
two solutions to the problem (\ref{2-b-E2.1})-(\ref{2-b-E2.6})
satisfying the regularity (\ref{2-b-E2.44-1})-(\ref{2-b-E2.44-3}),
and we denote
    $$
    \overline{\rho}=\rho_1-\rho_2,
    \ \overline{\theta}=\theta_1-\theta_2,
    \ \overline{u}=u_1-u_2.
    $$
First, since $\overline{\rho}\in L^\infty(0,T_{***};L^\infty)$ is a
solution to the linear transport equation
$\overline{\rho}_t+\mathrm{div}(\overline{\rho}v)=0$, it follows
from the uniqueness result by DiPerna-Lions in \cite{DiPerna89} that
$\overline{\rho}=0$, that is, $\rho_1=\rho_2:=\rho$ in
$[0,T_{***}]\times\Omega$. Next, to show that $\overline{\theta}=0$
in $[0,T_{***}]\times\Omega$, we multiply the both sides of
    \begin{equation}
      e_\theta(\rho \overline{\theta}_t+\rho v\cdot\nabla \overline{\theta})-
      \mathrm{div}(\kappa\nabla \overline{\theta})+\mathrm{div}v
      (p(\rho,\theta_1)-p(\rho,\theta_2)
      )=0
      \label{2-b-E2.45}
    \end{equation}
by $\overline{\theta}$ and integrate over $(0,t)\times\Omega$.
Then, since $\rho_t+\mathrm{div}(\rho v)=0$ and
$\overline{\theta}(0)=0$ in $\Omega$, we deduce at least formally
that
    \begin{eqnarray}
      &&\frac{1}{2}\int e_\theta\rho|\overline{\theta}|^2 dx
      +\int^t_0\int\kappa|\nabla
      \overline{\theta}|^2dxds\nonumber\\
      &=&\frac{1}{2}\int^t_0\int\left[
      \partial_t e_\theta \rho \overline{\theta}^2+\nabla e_\theta
      \rho v\overline{\theta}^2-\overline{\theta}\mathrm{div}v
      (p(\rho,\theta_1)-p(\rho,\theta_2)    )
      \right]dxds\nonumber\\
            &\leq& \frac{\kappa^0}{2}\int^t_0\int|\nabla
      \overline{\theta}|^2dxds+C\int^t_0\left[\|\partial_te_\theta\|_{L^3}^2+\|\nabla
      e_\theta\|_{L^3}^2
      \|v\|_{L^\infty}^2\right.\nonumber\\
      &&\left.+\|\nabla v\|_{L^3}^2\|p_1\|_{L^\infty}^2(1+\|(\theta_1,\theta_2)
      \|_{L^\infty}^2)\right]\|\sqrt{\rho}\overline{\theta}\|_{L^2}^2ds
      \label{2-b-E2.46}
    \end{eqnarray}
and applying Gronwall's inequality, we also conclude that
$\overline{\theta}=0$ in $(0,t)\times\Omega$. But this argument is
somewhat formal since it is not obvious that
$\rho^{\frac{1}{2}}\overline{\theta}\in L^\infty(0,T_{***};L^2)$
for the case of unbounded domains. Hence we have to justify this
formal argument by deriving  (\ref{2-b-E2.46}) rigorously. For
this purpose, we assume that $\Omega$ is an unbounded domain, let
$\Omega_r=\Omega\cap B_r$, $r>1$, and define
$\overline{\theta}_r\in L^\infty(0,T_{***};H^1_0(\Omega_r))$ by
    $$
    \overline{\theta}_r(t,x)=\overline{\theta}(t,x)\phi_r(x)
    \textrm{ for }(t,x)\in[0,T_{***}]\times\Omega,
    $$
where $\phi_r(x)=\phi(\frac{x}{r})$ is a cut-off function
satisfying $\phi\in C^\infty_0(B_1)$ and $\phi\equiv1$ in
$B_{\frac{1}{2}}$. Then from (\ref{2-b-E2.45}), we derive
    $$
    e_\theta\left(\rho (\overline{\theta}_r)_t+\rho v\cdot\nabla \overline{\theta}_r
    \right)-
      \mathrm{div}(\kappa\nabla \overline{\theta})
      \phi_r+\mathrm{div}v(p(\rho,\theta_1)-p(\rho,\theta_2)
      )\phi_r=e_\theta\rho
      \overline{\theta}v\cdot\nabla\phi_r.
    $$
Hence multiplying this by $\overline{\theta}_r$ and integrating
over $[0,t]\times\Omega$, we obtain that
    \begin{eqnarray}
      &&\frac{1}{2}\int e_\theta\rho|\overline{\theta}_r|^2(t)dx
      +\int^t_0\int\kappa|\nabla
      \overline{\theta}_r|^2dxds\nonumber\\
      &=&\frac{1}{2}\int^t_0\int\left[
      \partial_t e_\theta \rho \overline{\theta}_r^2+\nabla e_\theta
      \rho v\overline{\theta}_r^2-\overline{\theta}_r
      \phi_r\mathrm{div}v(p(\rho,\theta_1)-p(\rho,\theta_2)    )
      \right]dxds
      +I_r(t),\label{2-b-E2.47}
    \end{eqnarray}
where the remainder term $I_r(t)$ satisfies
    \begin{eqnarray*}
      |I_r(t)|&=&\left|\int^t_0\int
      \kappa\overline{\theta}^2
      (\nabla
      \phi_r)^2+e_\theta\overline{\theta}\rho\overline{\theta}_rv\cdot\nabla
      \phi_r
      \right|\\
            &\leq& C\int^t_0\int
            \rho|e_\theta||\overline{\theta}||\overline{\theta}_r||v|
            |\nabla\phi_r|+\kappa\overline{\theta}^2|\nabla\phi_r|^2dxds.
    \end{eqnarray*}
Since $e_\theta,\kappa,\rho\in L^\infty(0,T_{***};L^\infty)$ and
$(\overline{\theta},v)\in L^\infty(0,T_{***};D_0^1)$, it follows
that
    \begin{eqnarray*}
      &&|I_r(t)|\\
      &\leq& \tilde{C}\int^t_0\int\rho|\overline{\theta}_r|^2dxds
      +\frac{\tilde{C}}{r^2}\int^t_0\int_{\Omega_r}
      (\rho|\overline{\theta}|^2|v|^2+\overline{\theta}^2)dxds\\
            &\leq&\tilde{C}\int^t_0\int\rho|\overline{\theta}_r|^2dxds
            +\frac{\tilde{C}}{r^2}\int^t_0\left[
              \|\rho\|_{L^\infty}\|\overline{\theta}^2\|_{L^3}
              \|v^2\|_{L^3}|\Omega_r|^{\frac{1}{3}}
              +\|\overline{\theta}\|_{L^6}^2|\Omega_r|^{\frac{2}{3}}\right]ds\\
    &\leq&\tilde{C}\int^t_0\int\rho|\overline{\theta}_r|^2dxds
              +\tilde{C}
    \end{eqnarray*}
where the positive constant $\tilde{C}$ is independent of $r$.
Therefore, substituting this estimate into (\ref{2-b-E2.47}),
using the similar argument in (\ref{2-b-E2.46}) and applying
Gronwall's inequality, we conclude that
    $$
    \sup_{0\leq t\leq T_{***}}\|\rho^\frac{1}{2}\overline{\theta}(t)\|^2_{L^2(\Omega_{\frac{r}{2
    }})}+\int^{T_{***}}_0\|\nabla\overline{\theta}(t)\|^2_{L^2(\Omega_{\frac{r}{2}})}dt
    \leq \tilde{C}.
    $$
Thus, we have $\rho^{\frac{1}{2}}\overline{\theta}\in
L^\infty(0,T_{***};L^2(\Omega))$ and $\nabla \overline{\theta}\in
L^2(0,T_{***};L^2(\Omega))$. Now, we can estimate $I_r(t)$ again,
    \begin{eqnarray*}
      |I_r(t)|&\leq&\tilde{C}\int^{T_{***}}_0\left(
      \|\nabla \overline{\theta}\|_{L^2}+\|\rho\|_{L^\infty}^\frac{1}{2}
      \|\rho^\frac{1}{2}\overline{\theta}\|_{L^2}\|v\|_{L^\infty}
      \right)\|\nabla
      \overline{\theta}\|_{L^2(\Omega\backslash B_{\frac{r}{2}})}ds\\
            &\leq& \tilde{C}\left(
            \int^{T_{***}}_0\|\nabla \overline{\theta}\|^2_{L^2(\Omega
            \backslash B_{\frac{r}{2}})}
            ds  \right)^\frac{1}{2}\rightarrow0,
            \ \textrm{ as }r\rightarrow\infty.
    \end{eqnarray*}
Hence, letting $r\rightarrow\infty$ in (\ref{2-b-E2.47}), we could
obtain the estimate (\ref{2-b-E2.46}) and $\overline{\theta}=0$ in
$(0,T_{***})\times\Omega$. Similarly, we have $\overline{u}=0$. This
completes the proof of the uniqueness.

Finally, we prove the time-continuity of the solution
$(\rho,\theta,u)$ with the regularities
(\ref{2-b-E2.44-1})-(\ref{2-b-E2.44-3}). From Lemma
\ref{2-b-L2.1}, we have $\rho\in C([0,T_{***}];L^6\cap D^1\cap
D^{1,q})$. From (\ref{2-b-E2.1})  and $v\in
C([0,T_{***}];D^1_0\cap D^{2})$, we can obtain $\rho_t\in
C([0,T_{***}];L^2\cap L^q)$. To show the time-continuity of
$(\theta,u)$, we first observe that $(\theta,u)\in
C([0,T_{***}];D^1_0\cap D^{1,q'})\cap C([0,T_{***}];D^2-weak)$,
$\forall$ $q'\in(3,6)$. From the equations
(\ref{2-b-E2.3})-(\ref{2-b-E2.2}), we also observe that
$((e_\theta\rho \theta_t)_t$, $(\rho u_t)_t)\in
L^2(0,T_{***};H^{-1})$. Since $(e_\theta\rho \theta_t,\rho u_t)\in
L^2(0,T_{***};H^1_0)$, it follows immediately that $(e_\theta\rho
\theta_t,\rho u_t)\in C(0,T_{***};L^2)$. Hence we conclude that
for each $t\in[0,T_{***}]$, $(\theta,u)=(\theta(t),u(t))\in
D_0^1\cap D^2$ is a solution to the elliptic system
$-\mathrm{div}(\kappa\nabla e)=F$ and $L(\sigma,\eta, u)=G$ in
$\Omega$ for some $(F,G)\in C([0,T_{***}];L^2)$. Since
$(\kappa,\mu,\lambda)(\sigma,\eta)\in C([0,T_{***}];C^\beta\cap
D^1\cap D^{1,q})$, using the elliptic regularity estimates
(\ref{2-b-E2.33-1})-(\ref{2-b-E2.33-2-1}), we easily show that
$(\theta,u)\in C([0,T_{***}];D^2)$. From this time-continuity of
the solution $(\rho,\theta,u)$ and the estimates
(\ref{2-b-E2.42-3.0})-(\ref{2-b-E2.42-3.1}), we obtain
$$
    \|\rho(t,\cdot)-\rho_0(\cdot)\|_{L^6\cap D^1\cap D^{1,q}}
    +\|(\theta(t,\cdot)-\theta_0(\cdot),u(t,\cdot)-u_0(\cdot))\|_{D^1_0\cap D^{2}}
    \rightarrow0,\ \textrm{ as }\ t\rightarrow0.
        $$
We have completed the proof of Lemma \ref{2-b-L2.4}.
\end{proof}
\section{Existence result for polytropic fluids}\label{2-b-Sec3}\setcounter{equation}{0}
In this section, we consider the
following initial boundary value problem for a viscous polytropic
fluid:
    \begin{equation}
      \rho_t+\mathrm{div}(\rho u)=0,
      \label{2-b-E3.1}
    \end{equation}
    \begin{equation}
      \rho u_t+\rho u \cdot\nabla u+L(\rho,\theta,u)+\nabla
      p=\rho f, \textrm{ in }(0,\infty)\times\Omega,
      \label{2-b-E3.3}
    \end{equation}
        \begin{equation}
          e_\theta(\rho \theta_t+\rho u\cdot\nabla \theta)-\mathrm{div}(\kappa\nabla \theta)
          +p\mathrm{div}u=Q(\rho,\theta ,u )+\rho h+e_\rho\rho^2\mathrm{div}u,
          \label{2-b-E3.2}
        \end{equation}
            \begin{equation}
          (\rho, \theta, u)|_{t=0}=(\rho_0,\theta_0,u_0)
          \textrm{ in }\Omega,
          \label{2-b-E3.4}
        \end{equation}
    \begin{equation}
      (\theta,u)=(0,0)
      \textrm{ on }(0,\infty)\times\partial \Omega,
      \label{2-b-E3.5}
    \end{equation}
        \begin{equation}
          (\rho,\theta,u)(t,x)\rightarrow(\rho^\infty,0,0)
          \textrm{ as }|x|\rightarrow \infty,
          \label{2-b-E3.6}
        \end{equation}
where $Q(\rho,\theta, u)=\frac{\mu(\rho,\theta)}{2}|\nabla
u+\nabla u^\top|^2+\lambda(\rho,\theta)(\mathrm{div}u)^2$,
$L(\rho,\theta,u)=-\mathrm{div}(\mu(\rho,\theta)(\nabla u+\nabla
u^\top))-\nabla(\lambda(\rho,\theta)\mathrm{div}u)$ and
$(e_\theta,e_\rho,p,\mu,\lambda,\kappa)=(e_\theta,e_\rho,p,\mu,\lambda,\kappa)
(\rho,\theta)$.

\textbf{Proof of Theorem \ref{2-b-T1}}

Our proof will be based on the usual iteration argument and
results in  Lemma \ref{2-b-L2.4}. Let us denote
    $$
    c_0=2+\rho^\infty+\|\rho_0-\rho^\infty\|_{L^6\cap D^1\cap D^{1,q}}
    +\|(\theta_0,u_0)\|_{D^1_0\cap D^2}+\|(g_1,g_2)\|_{L^2}^2,
    $$
and we choose the positive constants $c_1$,  $c_2$, $c_3$, $c_4$
and $T_{4}$ as in Section \ref{2-b-Sec2}, depend only on $c_0$.
Next, let $u^0\in C([0,\infty);D^1_0\cap D^2)\cap
L^2(0,\infty;D^3)$ be the solution to the linear parabolic problem
    $$
    w_t-\Delta w=0
    \textrm{ in }(0,\infty)\times\Omega
    \textrm{ and }w(0)=u_0
    \textrm{ in }\Omega.
    $$
Then, chose a small time $T_5\in (0,T_{4}]$, such that
    $$
    \sup_{0\leq t\leq T_5}
    \|u^0(t)\|_{D^1_0}
    +\int^{T_5}_0\left(\|u^0_t(t)\|_{D^1_0}^2
    +\|u^0(t)\|^2_{D^{2,q}}
    \right)dt\leq c_1,
    $$
    $$
    \sup_{0\leq t\leq T_5}
    \|u^0\|_{D^2}\leq c_2.
    $$
Hence, choosing $T=T_5$ in Section \ref{2-b-Sec2}, it follows from
Lemma \ref{2-b-L2.4} that there exists a unique strong solution
$(\rho^1,e^1,u^1)$ to the linearized problem
(\ref{2-b-E2.1})-(\ref{2-b-E2.6}) with $(\sigma,\eta,v)$ replaced
by $(\rho_0,\theta_0,u^0)$, which satisfies the regularity
estimates (\ref{2-b-E2.42-1})-(\ref{2-b-E2.42-3.1}) with $T_{***}$
replaced by $T_5$. Similarly, we construct approximate solutions
$(\rho^k,\theta^k,u^k)$ as follows: assuming that
$(\rho^{k-1},\theta^{k-1},u^{k-1})$ was defined for $k\geq2$, let
$(\rho^k,\theta^k,u^k)$ be the unique solution to the linearized
problem (\ref{2-b-E2.1})-(\ref{2-b-E2.6}) with $(\sigma,\eta,v)$
replaced by $(\rho^{k-1},\theta^{k-1},u^{k-1})$. Then it also
follows from Lemma \ref{2-b-L2.4} that there exists a constant
$\tilde{C}>1$ such that
    \begin{eqnarray}
   && \sup_{0\leq t\leq T_5}\left(\|\rho^k(t)-\rho^\infty\|_{L^6\cap D^1\cap D^{1,q}}
    +\|\rho^k_t(t)\|_{L^2\cap L^q}
    \right)\nonumber\\
        &&+\sup_{0\leq t\leq T_5}\|(\theta^k,u^k)(t)\|_{D^1_0\cap D^2}
        +\sup_{0\leq t\leq T_5}\|(\sqrt{\rho^k}\theta^k_t,\sqrt{\rho^k}u^k_t)\|_{L^2}\nonumber\\
    &&+\int^{T_5}_0\left(\|(\theta^k_t,u^k_t)(t)\|^2_{D^1_0}+\|
    (\theta^k,u^k)(t)\|^2_{D^{2,q}}
    \right)dt\leq \tilde{C}\label{2-b-E3.9}
    \end{eqnarray}
    \begin{equation}
    \|(\theta^k(t,\cdot)-\theta_0(\cdot),u^k(t,\cdot)-u_0(\cdot))\|_{D^1_0}
    \leq \tilde{C}t^{\frac{1}{2}},
    \end{equation}
and
        \begin{equation}
    \|\rho^k(t,\cdot)-\rho_0(\cdot)\|_{L^6}
    \leq \tilde{C}t^{\frac{1}{3}},
    \end{equation}
for all $k\geq1$. Throughout the proof, we denote by $\tilde{C}$ a
generic constant depending only on $c_0$, but independent of $k$.

From now on, we show that the  sequence $(\rho^k,\theta^k,u^k)$
converges to a solution to the original nonlinear problem
(\ref{2-b-E3.1})-(\ref{2-b-E3.6}) in a strong sense. Let us define
    $$
    \overline{\rho}^{k+1}=\rho^{k+1}-\rho^k,
    \ \overline{\theta}^{k+1}={\theta}^{k+1}-\theta^k,
    \ \overline{u}^{k+1}=u^{k+1}-u^k.
    $$
Then, from (\ref{2-b-E3.1})-(\ref{2-b-E3.2}), we derive equations
for   functions
$(\overline{\rho}^{k+1},\overline{\theta}^{k+1},\overline{u}^{k+1})$,
    \begin{equation}
    \overline{\rho}^{k+1}_t+\mathrm{div}(\overline{\rho}^{k+1}u^k)
    +\mathrm{div}(\rho^k\overline{u}^k)=0,
    \label{2-b-E3.10}
    \end{equation}
    \begin{eqnarray}
      &&\rho^{k+1}\overline{u}^{k+1}_t+\rho^{k+1}u^k\cdot\nabla\overline{u}^{k+1}
      +L(\rho^{k},\theta^{k},\overline{u}^{k+1})\nonumber\\
            &=&\overline{\rho}^{k+1}(f-u^k_t-u^{k-1}\cdot\nabla u^k)
            -\rho^{k+1}\overline{u}^k\cdot\nabla u^k
            +\nabla(p^k-p^{k+1})\nonumber\\
      &&+\mathrm{div}((\mu^{k}-\mu^{k-1})(\nabla u^k+(\nabla
      u^k)^\top))\nonumber\\
            &&+\nabla((\lambda^{k}-
            \lambda^{k-1})\mathrm{div}u^k),
      \label{2-b-E3.12}
    \end{eqnarray}
and
    \begin{eqnarray}
          &&e_\theta^k\left(\rho^{k+1}\overline{\theta}^{k+1}_t+\rho^{k+1}u^k\cdot\nabla
          \overline{\theta}^{k+1}\right)-\mathrm{div}(\kappa^k\nabla
          \overline{\theta}^{k+1})\nonumber\\
                &=&\mathrm{div}((\kappa^{k}-\kappa^{k-1})\nabla \theta^k)
                +Q(\rho^{k},\theta^k,u^k)-Q(\rho^{k-1},
                \theta^{k-1},u^{k-1})\nonumber\\
          &&+\overline{\rho}^{k+1}h+p^k\mathrm{div}u^{k-1}-
          p^{k+1}\mathrm{div}
          u^k+e_\theta^{k-1}\rho^k(\theta^k_t+u^{k-1}\cdot\nabla \theta^k)\nonumber\\
          &&
          -e_\theta^{k}\rho^{k+1}(\theta^{k}_t+u^{k}\cdot\nabla
          \theta^k) +e_\rho^k(\rho^{k+1})^2\mathrm{div}u^k-
                e_\rho^{k-1}(\rho^{k})^2\mathrm{div}u^{k-1},
                \label{2-b-E3.11}
        \end{eqnarray}
where
$(e_\theta^k,e_\rho^k,p^k,\kappa^k,\mu^k,\lambda^k)=(e_\theta,e_\rho
,p,\kappa,\mu,\lambda)(\rho^k,\theta^k)$.

 First, we consider the case that $\rho^\infty>0$. Multiplying
(\ref{2-b-E3.10}) by $\overline{\rho}^{k+1}$ and integrating over
$\Omega$, we obtain (at least formally) that
    \begin{eqnarray*}
      &&\frac{d}{dt}\int|\overline{\rho}^{k+1}|^2dx\\
            &\leq&C\int\left(|\nabla u^k||\overline{\rho}^{k+1}|^2
            +|\nabla \rho^k||\overline{u}^k||\overline{\rho}^{k+1}|
            +\rho^k|\nabla \overline{u}^k||\overline{\rho}^{k+1}|
            \right)dx\\
      &\leq&C\left(\|\nabla
      u^k\|_{W^{1,q}}\|\overline{\rho}^{k+1}\|^2_{L^2}
      +(\|\nabla \rho^k\|_{L^3}+\|\rho^k\|_{L^\infty})\|\nabla
      \overline{u}^k\|_{L^2}\|\overline{\rho}^{k+1}\|_{L^2}
      \right).
    \end{eqnarray*}
Hence, by virtue of Young's inequality, we have
    \begin{equation}
      \frac{d}{dt}\|\overline{\rho}^{k+1}\|_{L^2}^2
      \leq A^k_\epsilon(t)\|\overline{\rho}^{k+1}\|^2_{L^2}
      +\epsilon\|\nabla \overline{u}^k\|^2_{L^2},
      \label{2-b-E3.13}
    \end{equation}
where $A^k_\epsilon(t)=C\|\nabla
u^k(t)\|_{W^{1,q}}+\epsilon^{-1}C(\|\nabla
\rho^k(t)\|^2_{L^3}+\|\rho^k(t)\|^2_{L^\infty})$. From the
estimate (\ref{2-b-E3.9}), we have $A^k_\epsilon(t)\in L^1(0,T_5)$
and $\int^t_0A^k_\epsilon(s)ds\leq \tilde{C}+\tilde{C}_\epsilon t$
for all $k\geq1$ and $t\in[0,T_5]$. Here  we denote by
$\tilde{C}_\epsilon$ a generic positive constant depending only on
$\epsilon^{-1}$ and the parameters of $\tilde{C}$, where
$\epsilon\in(0,1)$ is a small number.

Next, multiplying (\ref{2-b-E3.12}) by $\overline{u}^{k+1}$ and
integrating over $\Omega$, recalling that
    \begin{equation*}
      (\rho^{k+1})_t+\mathrm{div}(\rho^{k+1}u^k)=0
      \textrm{ in }\Omega,
    \end{equation*}
and using the technique in the proof of (\ref{2-b-E2.55-1}), we
have (at least formally) that
    \begin{eqnarray*}
      &&\frac{1}{2}\frac{d}{dt}\int\rho^{k+1}|\overline{u}^{k+1}|^2dx
      +\mu^0\int|\nabla\overline{u}^{k+1}|^2dx\\
            &\leq &C\int\left[|\overline{\rho}^{k+1}|(|f|
            +|u^k_t|+|u^{k-1}\cdot\nabla u^k|)|\overline{u}^{k+1}|
            +|\rho^{k+1}||\overline{u}^{k}||\nabla
            u^k||\overline{u}^{k+1}|
            \right.\\
      &&+|p^{k+1}-p^k|
      |\nabla \overline{u}^{k+1}|
            \left.+(|\mu^k-\mu^{k-1}|
            +|\lambda^k-\lambda^{k-1}|
              )|\nabla u^k||\nabla \overline{u}^{k+1}|\right]dx.
    \end{eqnarray*}
Then, using H\"{o}lder, Sobolev and Young's inequalities, it also
follows from (\ref{2-b-E3.9}) that
    \begin{eqnarray}
      &&\frac{d}{dt}\|\sqrt{\rho^{k+1}}\overline{u}^{k+1}\|_{L^2}^2
      +\mu^0\|\nabla\overline{u}^{k+1}\|_{L^2}^2\nonumber\\
            &\leq&
            B^k_\epsilon(t)\left(\|\overline{\rho}^{k+1}\|_{L^2}^2
                        +\|\overline{\rho}^{k}\|_{L^2}^2
            +\|\sqrt{\rho^{k+1}}\overline{u}^{k+1}\|_{L^2}^2
                        +\|\sqrt{\rho^{k+1}}\overline{\theta}^{k+1}\|_{L^2}^2
            \right)+\tilde{C_\epsilon}\|\sqrt{\rho^{k}}\overline{\theta}^{k}\|_{L^2}^2
            \nonumber\\
      &&
        +\frac{\epsilon^2}{8}\|\nabla\overline{\theta}^k\|_{L^2}^2
      +\epsilon\|\nabla\overline{u}^k\|_{L^2}^2+C\int|\overline{\rho}^{k+1}|
      |u^k_t||\overline{u}^{k+1}|dx,
      \label{2-b-E3.16}
    \end{eqnarray}
where $B^k_\epsilon(t)\in L^1(0,T_5)$ such that
$\int^t_0B^k_\epsilon(s)ds\leq \tilde{C}+\tilde{C}_\epsilon t$ for
$0\leq t\leq T_5$ and $k\geq1$.

 Finally, multiplying (\ref{2-b-E3.11}) by $\overline{\theta}^{k+1}$,
integrating over $\Omega$ and recalling that
    \begin{equation*}
      (\rho^{k+1})_t+\mathrm{div}(\rho^{k+1}u^k)=0
      \textrm{ in }\Omega,
    \end{equation*}
we obtain (at least formally) that
    \begin{eqnarray*}
      &&\frac{1}{2}\frac{d}{dt}\int e_\theta^k\rho^{k+1}|\overline{\theta}^{k+1}|^2dx+\kappa^0
      \int|\nabla \overline{\theta}^{k+1}|^2dx\\
            &\leq &C\int\left\{|\kappa^{k}-\kappa^{k-1}||\nabla \theta^k||\nabla \overline{\theta}^{k+1}|
            \right.+(|\mu^{k}|+|\lambda^{k}|)\left(|\nabla u^k|+|\nabla u^{k-1}|
            \right)|\nabla \overline{u}^k||\overline{\theta}^{k+1}|\\
                &&+|\overline{\theta}^{k+1}||\nabla u^{k-1}|^2\left(
                |\mu^{k}-\mu^{k-1}|
                +|\lambda^{k}-\lambda^{k-1}|
                \right)+|\overline{\rho}^{k+1}||h||\overline{\theta}^{k+1}|\\
                        &&+
                        |\overline{\theta}^{k+1}|\left|p^k\mathrm{div}u^{k-1}-
                      p^{k+1}\mathrm{div}
                      u^k\right.
          +e_\theta^{k-1}\rho^k(\theta^k_t+u^{k-1}\cdot\nabla \theta^k)
          -e_\theta^{k}\rho^{k+1}(\theta^{k}_t+u^{k}\cdot\nabla
          \theta^k)\nonumber\\
                &&\left.\left.+e_\rho^k(\rho^{k+1})^2\mathrm{div}u^k-
                e_\rho^{k-1}(\rho^{k})^2\mathrm{div}u^{k-1}
                +\frac{1}{2}(\partial_t e_\theta^k+
                u^k\cdot \nabla e_\theta^k)\rho^{k+1}\overline{\theta}^{k+1}\right|\right\}dx.
    \end{eqnarray*}
From (\ref{2-b-E3.9}), using H\"{o}lder, Sobolev and Young's
inequalities, we have
     \begin{eqnarray}
      &&\frac{d}{dt}\|\sqrt{e_{\theta}^k\rho^{k+1}}\overline{\theta}^{k+1}\|_{L^2}^2
      +\kappa^0\|\nabla \overline{\theta}^{k+1}\|_{L^2}^2\nonumber\\
            &\leq&
            D^k_{\epsilon}(t)\left(\|\overline{\rho}^{k+1}\|_{L^2}^2+
            \|\overline{\rho}^{k}\|_{L^2}^2+
            \|
            \sqrt{\rho^{k+1}}\overline{\theta}^{k+1}\|_{L^2}^2
            \right)+\tilde{C}_\epsilon\|
            \sqrt{\rho^{k}}\overline{\theta}^{k}\|_{L^2}^2\nonumber\\
                    &&+\tilde{C}\|\sqrt{\rho^k}\overline{u}^k\|_{L^2}^2+\tilde{C}\|\nabla
                    \overline{u}^k\|_{L^2}^2
                    +\frac{\kappa^0\epsilon}{8}\|\nabla \overline{\theta}^{k}\|_{L^2}^2
                    +\tilde{C}\int|\overline{\rho}^{k+1}||\theta^k_t||\overline{\theta}^{k+1}|dx
            \label{2-b-E3.15}
            \end{eqnarray}
where $D^k_{\epsilon}(t)\in L^1(0,T_5)$ such that $\int^t_0
D^k_{\epsilon}(s)ds\leq \tilde{C}+\tilde{C}_{\epsilon}t$ and
$k\geq1$.

Therefore, combining (\ref{2-b-E3.13})-(\ref{2-b-E3.15}), defining
    $$
    \psi^{k+1}(t)=\|\overline{\rho}^{k+1}\|_{L^2}^2
    +\frac{\epsilon}{\kappa^0}\|\sqrt{e_\theta^k\rho^{k+1}}\overline{\theta}^{k+1}\|_{L^2}^2
    +\|\sqrt{\rho^{k+1}}\overline{u}^{k+1}\|_{L^2}^2,
    $$
we deduce that
    \begin{eqnarray}
      &&\frac{d}{dt}(\psi^{k+1}+\|\overline{\rho}^k\|_{L^2}^2)+\epsilon\|\nabla
      \overline{\theta}^{k+1}\|_{L^2}^2+\mu^0\|
      \nabla \overline{u}^{k+1}\|_{L^2}^2\nonumber\\
            &\leq& E^k_\epsilon(t)\left(\psi^{k+1}+\|\overline{\rho}^k\|_{L^2}^2
            \right)+\tilde{C}_\epsilon\psi^{k}
        +\frac{\epsilon^2}{4}\|\nabla
      \overline{\theta}^{k}\|_{L^2}^2
            +3\epsilon \tilde{C}
            (\|\nabla \overline{u}^{k}\|_{L^2}^2+
            \|\nabla \overline{u}^{k-1}\|_{L^2}^2)\nonumber\\
                &&+\tilde{C}
            \int|\overline{\rho}^{k+1}|\left(\epsilon
            |\theta^k_t||\overline{\theta}^{k+1}|+|u^k_t||\overline{u}^{k+1}|
            \right)dx,
            \label{2-b-E3.17}
    \end{eqnarray}
where $E^k_\epsilon(t)\in L^1(0,T_5)$ such that $\int^t_0
E^k_\epsilon(s)ds\leq \tilde{C}+\tilde{C}t$ for $0\leq t\leq T_5$
and $k\geq 2$.

To estimate the last integral term in (\ref{2-b-E3.17}), we assume
 that $\Omega$ is an unbounded domain in
$\mathbb{R}^3$, and Claim that there exists a large radius $R_0>0$
and a small time $T_6\in (0,T_5]$ independent of k, such that
    \begin{equation}
    \frac{1}{4}\rho^\infty
      \leq \rho^{k+1}(t,x)\leq 4\rho^\infty,
      \ \forall\ (t,x)\in[0,T_6]\times(\Omega\backslash
      B_{R_0}).
      \label{2-b-E3.19}
    \end{equation}
From (\ref{2-b-E3.19}), we can  estimate the integral term in
(\ref{2-b-E3.17}) as follows: for $0\leq t\leq T_6$,
    \begin{eqnarray*}
      &&\tilde{C}\int_{\Omega\cap B_{R_0}}|\overline{\rho}^{k+1}|\left(\epsilon
      |\theta^k_t||\overline{\theta}^{k+1}|+|u^k_t||\overline{u}^{k+1}|
      \right)dx\\
            &\leq& \tilde{C}
            \left(\|\nabla \theta^{k}_t\|_{L^2}^2+\|\nabla
            u^k_t\|_{L^2}^2
            \right)\|\overline{\rho}^{k+1}\|_{L^2}^2
            +\frac{1}{4}\left(\epsilon\|\nabla
            \overline{\theta}^{k+1}\|_{L^2}^2+\mu^0\|\nabla
            \overline{u}^{k+1}\|_{L^2}^2
            \right)
    \end{eqnarray*}
and
    \begin{eqnarray*}
      &&\tilde{C}\int_{\Omega\backslash B_{R_0}}|\overline{\rho}^{k+1}|\left(\epsilon
      |\theta^k_t||\overline{\theta}^{k+1}|+|u^k_t||\overline{u}^{k+1}|
      \right)dx\\
            &\leq&\tilde{C}\|\overline{\rho}^{k+1}\|_{L^2}
            \left(\epsilon\|\theta^k_t\|_{L^6}\|\overline{\theta}^{k+1}\|_{L^3}
            +\|u^k_t\|_{L^6}\|\overline{u}^{k+1}\|_{L^3}
            \right)\\
      &\leq& \frac{\tilde{C}\|\overline{\rho}^{k+1}\|_{L^2}}
      {(\rho^{\infty})^{\frac{1}{4}}}\left(\epsilon\|\nabla \theta^k_t\|_{L^2}
      \|\sqrt{\rho^{k+1}}\overline{\theta}^{k+1}\|_{L^2}^\frac{1}{2}
      \|\nabla \overline{\theta}^{k+1}\|_{L^2}^\frac{1}{2}\right.\\
            &&\left.
      +\|\nabla u^k_t\|_{L^2}\|\sqrt{\rho^{k+1}}\overline{u}^{k+1}\|_{L^2}^\frac{1}{2}
      \|\nabla\overline{u}^{k+1}\|_{L^2}^\frac{1}{2}
      \right)\\
            &\leq&\tilde{C}(\rho^\infty)
            \left(\|\nabla \theta^k_t\|_{L^2}^2+
            \|\nabla u^k_t\|_{L^2}^2+1\right)\psi^{k+1}
            +\frac{1}{4}\left(\epsilon\|\nabla
            \overline{\theta}^{k+1}\|_{L^2}^2+\mu^0\|\nabla
            \overline{u}^{k+1}\|_{L^2}^2
            \right)
    \end{eqnarray*}
where the constant $\tilde{C}(\rho^\infty)$ depend also on
$\rho^\infty$. Therefore, substituting these estimates into
(\ref{2-b-E3.17}), we obtain
    \begin{eqnarray}
      &&\frac{d}{dt}(\psi^{k+1}+\|\overline{\rho}^k\|_{L^2}^2)+\frac{\epsilon}{2}\|\nabla\overline{\theta}^{k+1}\|_{L^2}^2
      +\frac{\mu^0}{2}\|\nabla
      \overline{u}^{k+1}\|_{L^2}^2\nonumber\\
      &\leq& F^k_\epsilon(t)\left(\psi^{k+1}+\|\overline{\rho}^k\|_{L^2}^2\right)+
      \tilde{C}_\epsilon\psi^{k}+
      \frac{\epsilon^2}{4}\|\nabla
      \overline{\theta}^{k}\|_{L^2}^2+3\epsilon\tilde{C}(\|\nabla
      \overline{u}^k\|_{L^2}^2+\|\nabla
      \overline{u}^{k-1}\|_{L^2}^2)
      \label{2-b-E3.21}
    \end{eqnarray}
where $F^k_\epsilon\in L^1(0,T_6)$, such that $\int^t_0
F^k_\epsilon(s)ds\leq \tilde{C}(\rho^\infty)+\tilde{C}_\epsilon t$
for $k\geq2$ and $0\leq t\leq T_6$. It is easy to see that this
estimate holds also for a bounded domain $\Omega$ since we can
choose a sufficiently large $R_0$ satisfying $\Omega\subset
B_{R_0}$.

Now, we can prove (\ref{2-b-E3.19}) as follows. Since
$\rho_0-\rho^\infty\in L^6\cap D^1\cap D^{1,q}\hookrightarrow
C_0$, where the space $C_0$ consists of continuous functions on
$\Omega$ vanishing at infinity, we can choose a large radius
$R_0>1$ independent of $k$, so that
    \begin{equation}
      \frac{1}{2}\rho^\infty
      \leq \rho_0(x)\leq 2\rho^\infty,
      \ \forall\ x\in\Omega\backslash B_{\frac{R_0}{2}}.
      \label{2-b-E3.18}
    \end{equation}
From Lemma \ref{2-b-L2.1}, we have
    \begin{equation}
      \rho^{k+1}(t,x)=\rho_0(U^{k+1}(0,t,x))\exp\left\{
      -\int^t_0\mathrm{div}u^k(s,U^{k+1}(s,t,x))ds\right\},
      \label{2-b-E3.20}
    \end{equation}
where $U^{k+1}(t,s,x)$ is the solution to the initial value
problem
    $$
      \left\{\begin{array}{rlll}
      \frac{\partial}{\partial
      t}U^{k+1}(t,s,x)&=&u^k(t,U^{k+1}(t,s,x)),&t\in[0,T_5],\\
            U^{k+1}(s,s,x)&=&x,&s\in[0,T_5],\ x\in\Omega.
      \end{array}
      \right.
    $$
Moreover, in view of (\ref{2-b-E3.9}), we deduce that
    $$
    \int^t_0|\mathrm{div}u^k(s,U^{k+1}(s,t,x))|ds
    \leq\int^t_0\|\nabla u^k(s)\|_{L^\infty}ds\leq \tilde{C}t^{\frac{1}{2}}
    $$
and
    \begin{eqnarray*}
      &&|U^{k+1}(0,t,x)-x|=|U^{k+1}(0,t,x)-U^{k+1}(t,t,x)|\\
            &\leq& \int^t_0|u^{k}(\tau,U^{k+1}(\tau,t,x))|
            \tau\leq \tilde{C}t,
    \end{eqnarray*}
for all $(t,x)\in [0,T_5]\times\Omega$. Letting
$T_6=\min(T_5,\tilde{C}^{-2}(ln2)^2,\frac{1}{2}\tilde{C}^{-1}R_0)$
is a small positive time dependent only on $T_5$ and the
parameters $\tilde{C}$, we obtain
    $$
    \int^t_0|\mathrm{div}u^k(s,U^{k+1}(s,t,x))|ds
    \leq \ln 2,\ \textrm{ and }
     \  |U^{k+1}(0,t,x)-x|\leq \frac{R_0}{2},
    $$
 for all $(t,x)\in [0,T_6]\times\Omega$. In particular, it follows that
if $0\leq t\leq T_6$ and $x\in\Omega\backslash B_{R_0}$, then
$U^{k+1}(0,t,x)\in \Omega\backslash B_{\frac{R_0}{2}}$. Hence the
desired result (\ref{2-b-E3.19}) follows immediately from
(\ref{2-b-E3.18})-(\ref{2-b-E3.20}).

Now, recalling that $\psi^{k+1}(0)=0$ and using Gronwall's
inequality, we deduce from (\ref{2-b-E3.21}) that
    \begin{eqnarray*}
        &&    \psi^{k+1}+\|\overline{\rho}^k\|_{L^2}+
      \int^{t}_0\frac{\epsilon}{2}\|\nabla\overline{\theta}^{k+1}\|_{L^2}^2ds
      +\int^{t}_0\frac{\mu^0}{2}\|\nabla
      \overline{u}^{k+1}\|_{L^2}^2ds\\
      &\leq&\tilde{C}_\epsilon\int^t_0\psi^{k}(s)ds\\
            &&+\left(\frac{\epsilon^2}{4}\int^t_0\|\nabla
      \overline{\theta}^{k}\|_{L^2}^2ds+3\epsilon\tilde{C}\int^t_0(\|\nabla
      \overline{u}^k\|_{L^2}^2+\|\nabla
      \overline{u}^{k-1}\|_{L^2}^2)ds
    \right)\exp\left(\tilde{C}(\rho^\infty)+\tilde{C}_\epsilon
    t\right)\\
        &\leq&\tilde{C}_\epsilon\int^t_0\psi^{k}(s)ds\\
                &&+\tilde{C}\epsilon
    \left(\frac{\epsilon}{2}\int^t_0\|\nabla
      \overline{\theta}^{k}\|_{L^2}^2ds+\frac{\mu^0}{2}\int^t_0(\|\nabla
      \overline{u}^k\|_{L^2}^2+\|\nabla
      \overline{u}^{k-1}\|_{L^2}^2)ds
    \right)\exp\left(\tilde{C}(\rho^\infty)+\tilde{C}_\epsilon
    t\right),
    \end{eqnarray*}
for $k\geq2$ and $0\leq t\leq T_6$. Hence choosing small constants
$\epsilon>0$ satisfying
$\tilde{C}\epsilon\exp(\tilde{C}(\rho^\infty))=\frac{1}{8}$, and
$T_7\in(0,T_6]$ satisfying $\exp(\tilde{C}_\epsilon  T_7)=2$, we
easily deduce that
        \begin{eqnarray*}
        &&    \psi^{k+1}+\|\overline{\rho}^k\|_{L^2}+
      \int^{t}_0\frac{\epsilon}{2}\|\nabla\overline{\theta}^{k+1}\|_{L^2}^2ds
      +\int^{t}_0\frac{\mu^0}{2}\|\nabla
      \overline{u}^{k+1}\|_{L^2}^2ds\\
        &\leq&\tilde{C}\int^t_0\psi^{k}(s)ds
                +\frac{1}{4}
    \left(\frac{\epsilon}{2}\int^t_0\|\nabla
      \overline{\theta}^{k}\|_{L^2}^2ds+\frac{\mu^0}{2}\int^t_0(\|\nabla
      \overline{u}^k\|_{L^2}^2+\|\nabla
      \overline{u}^{k-1}\|_{L^2}^2)ds
    \right).
    \end{eqnarray*}
for $k\geq2$ and $0\leq t\leq T_7$. Using  similar arguments as in
the proof of (\ref{2-b-E3.13})-(\ref{2-b-E3.15}), we can obtain
$\sup_{0\leq t\leq T_{7}}\psi^2(t)\leq \tilde{C}$ easily. Fixing a
large $K>1$, summing this over $2\leq k\leq K$ and using
Gronwall's inequality, we obtain
    \begin{eqnarray*}
    &&\sum_{k=2}^K\sup_{0\leq t\leq T_{7}}\psi^{k+1}(t)+
    \sum_{k=2}^K\int^{T_{7}}_0
    \left(\epsilon\|\nabla \overline{\theta}^{k+1}\|_{L^2}^2
    +\mu^0\|\nabla \overline{u}^{k+1}\|_{L^2}^2
    \right)dt\\
    &\leq& \tilde{C}\left(1+\int^{T_{7}}_0\psi^2(s)ds\right)\leq \tilde{C}.
    \end{eqnarray*}
Thus, we have
    $$
    \sum_{k=1}^\infty\sup_{0\leq t\leq T_*}\psi^{k+1}(t)+
    \sum_{k=1}^\infty\int^{T_*}_0
    \left(\epsilon\|\nabla \overline{\theta}^{k+1}\|_{L^2}^2
    +\mu^0\|\nabla \overline{u}^{k+1}\|_{L^2}^2
    \right)dt\leq \tilde{C}<\infty,
    $$
where $T_*=T_7$.

Therefore, we conclude that the full sequence
$(\rho^k,\theta^k,u^k)$ converges to a limit $(\rho,\theta,u)$ in
the following strong sense:
    \begin{equation}
      \left\{
      \begin{array}{lll}
      \rho^k-\rho^1\rightarrow \rho-\rho^1&in&L^\infty(0,T_*;L^2)\\
      ({\theta}^{k},u^k)\rightarrow(\theta,u)&in&L^2(0,T_*;D^1_0).
      \end{array}
      \right.\label{2-b-E3.22}
    \end{equation}
It is easy to show that the limit $(\rho,\theta,u)$ is a weak
solution to the original nonlinear problem
(\ref{2-b-E3.1})-(\ref{2-b-E3.6}). Furthermore, it follows from
(\ref{2-b-E3.9}) that $(\rho,\theta,u)$ satisfies the following
regularity estimate:
    \begin{eqnarray*}
      &&\sup_{0\leq t\leq T_*}
      \|(\sqrt{\rho}\theta_t,\sqrt{\rho}u_t)(t)
      \|_{L^2}+\int^{T_*}_0\left(\|(\theta_t,u_t)(t)\|_{D^1_0}^2+
      \|(\theta,u)\|_{D^{2,q}}^2
      \right)dt\\
            &&+\sup_{0\leq t\leq T_*}\left(
            \|\rho(t)-\rho^\infty\|_{L^6\cap D^1\cap D^{1,q}}+\|\rho_t(t)\|_{L^2\cap L^q}
            +\|(\theta,u)(t)\|_{D^1_0\cap D^2}
            \right)\leq \tilde{C}.
    \end{eqnarray*}
 Then adapting same
arguments in the Lemma \ref{2-b-L2.4}, we can easily prove the
time-continuity of the solution $(\rho,\theta,u)$ and
        $$
    \|\rho(t,\cdot)-\rho_0(\cdot)\|_{L^6\cap D^1\cap D^{1,q}}
    +\|(\theta(t,\cdot)-\theta_0(\cdot),u(t,\cdot)-u_0(\cdot))\|_{D^1_0\cap D^{2}}
    \rightarrow0,\ \textrm{ as }\ t\rightarrow0.
        $$ This
proves the existence of a strong solution. Now we prove the
uniqueness of the strong solutions. Let $(\rho_1,\theta_1,u_1)$
and $(\rho_2,\theta_2,u_2)$ be two strong solutions to the problem
(\ref{2-b-E3.1})-(\ref{2-b-E3.6}) with the regularities
(\ref{2-b-E3.8-1})-(\ref{2-b-E3.8-3}) and we denote by
$(\overline{\rho},\overline{\theta},\overline{u})$ their
difference. Then following the argument used to derive
(\ref{2-b-E3.21}), we easily deduce that
    \begin{eqnarray*}
      &&\frac{d}{dt}\left(\|\overline{\rho}\|^2_{L^2}
      +\|\sqrt{e_\theta(\rho_1,\theta_1)\rho_1}\overline{\theta}\|_{L^2}^2
      +\|\sqrt{\rho_1}\overline{u}\|_{L^2}^2
      \right)+\frac{\kappa^0}{2}\|\overline{\theta}\|_{D^1_0}
      +\tilde{C}\|\overline{u}\|_{D^1_0}\\
            &\leq&G(t)\left(\|\overline{\rho}\|^2_{L^2}
            +\|\sqrt{e_\theta(\rho_1,\theta_1)\rho_1}\overline{\theta}\|_{L^2}^2
            +\|\sqrt{\rho_1}\overline{u}\|_{L^2}^2
            \right)
    \end{eqnarray*}
where $G(t)\in L^1(0,T_*)$. Therefore, using Gronwall's
inequality, we conclude that
$\overline{\rho}=\overline{\theta}=\overline{u}=0$ in
$(0,T_{***})\times\Omega$. This completes the proof of the
existence and uniqueness part for the case that $\rho^\infty>0$.

Now we consider the case that $\rho^\infty=0$. To prove the
convergence in this case, we need to modify slightly  previous
arguments. First, multiplying (\ref{2-b-E3.10}) by
$\mathrm{sgn}(\overline{\rho}^{k+1})|\overline{\rho}^{k+1}|^\frac{1}{2}$
and integrating over $\Omega$, we obtain (at least formally) that
    \begin{eqnarray*}
    \frac{d}{dt}\int|\overline{\rho}^{k+1}|^\frac{3}{2}dx&\leq&
    C\int|\nabla u^k||\overline{\rho}^{k+1}|^\frac{3}{2}
    +(|\nabla \rho^k||\overline{u}^k|+\rho^k|\nabla \overline{u}^k|
    )|\overline{\rho}^{k+1}|^\frac{1}{2}dx\\
            &\leq&C\|\nabla u^k\|_{W^{1,q}}
            \|\overline{\rho}^{k+1}\|_{L^\frac{3}{2}}^\frac{3}{2}
            +C\|\rho^k\|_{L^6\cap D^1}\|\nabla \overline{u}^k\|_{L^2}
            \|\overline{\rho}^{k+1}\|_{L^\frac{3}{2}}^\frac{1}{2}.
    \end{eqnarray*}
Hence, multiplying this by
$\|\overline{\rho}^{k+1}\|_{L^\frac{3}{2}}^\frac{1}{2}$, we obtain
    \begin{equation}
      \frac{d}{dt}\|\overline{\rho}^{k+1}\|_{L^\frac{3}{2}}^2
      \leq
      A^k_\epsilon(t)\|\overline{\rho}^{k+1}\|_{L^\frac{3}{2}}^2+\epsilon\|
      \nabla \overline{u}^k\|_{L^2}^2,
      \label{2-b-E3.23}
    \end{equation}
where $A^k_\epsilon(t)=C\|\nabla
u^k(t)\|_{W^{1,q}}+\epsilon^{-1}C\| \rho^k(t)\|_{D^1_0}^2$. The
uniform bound (\ref{2-b-E3.9}) implies that
$\int^t_0A^k_\epsilon(s)ds\leq \tilde{C}+\tilde{C}_\epsilon t$ for
all $k\geq1$ and $t\in [0,T_5]$. In a similar manner, we can also
show that
    \begin{equation}
    \frac{d}{dt}\|\overline{\rho}^{k+1}\|_{L^2}^2
      \leq
      B^k_\epsilon(t)\|\overline{\rho}^{k+1}\|_{L^2}^2+\epsilon\|
      \nabla \overline{u}^k\|_{L^2}^2,
      \label{2-b-E3.24}
    \end{equation}
where $B^k_\epsilon(t)\in L^1(0,T_5)$ such that
$\int^t_0B^k_\epsilon(s)ds\leq \tilde{C}+\tilde{C}_\epsilon t$ for
all $k\geq1$ and $t\in [0,T_5]$.

Next, multiplying (\ref{2-b-E3.11}) by $\overline{\theta}^{k+1}$
and integrating over $\Omega$,
%we also obtain formally that
%    \begin{eqnarray*}
%      &&\frac{1}{2}\frac{d}{dt}\int\rho^{k+1}|\overline{\theta}^{k+1}|^2dx
%      +\kappa^0\int|\nabla \overline{\theta}^{k+1}|^2dx\\
%            &\leq& C\left(\|\mu(\rho^{k},\theta^k)\|_{L^\infty}+
%            \|\lambda(\rho^{k},\theta^k)\|_{L^\infty}
%            \right)
%            \left(\|\nabla u^k\|_{L^3}+\|\nabla u^{k-1}\|_{L^3}
%            \right)\|\nabla \overline{u}^k\|_{L^2}\|\nabla
%            \overline{\theta}^{k+1}\|_{L^2}\\
%       &&+C\int\left[\left(|\mu(\rho^k,\theta^k)-\mu(\rho^{k-1},\theta^{k-1})|+
%       |\lambda(\rho^k,\theta^k)-\lambda(\rho^{k-1},\theta^{k-1})|
%       \right)|     \overline{\theta}^{k+1}|
%            |\nabla u^{k-1}|^2\right]dx\\
%            &&+C\int\left(|\kappa(\rho^k,\theta^k)-\kappa(\rho^{k-1},\theta^{k-1})|
%            |\nabla \overline{\theta}^{k+1}||\nabla
%            \theta^k|\right)dx\\
%       &&+C\|\overline{\rho}^{k+1}\|_{L^\frac{3}{2}\cap L^2}
%       \left(\|\theta^k_t\|_{D^1_0}+\|h\|_{L^3}+\|\nabla u^{k-1}\|_{H^1}
%       \|\nabla \theta^k\|_{H^1}
%       \right)\|\nabla \overline{\theta}^{k+1}\|_{L^2}\\
%            &&+C\|\nabla u^k\|_{L^\infty}
%            \|\sqrt{\rho^{k+1}}\overline{\theta}^{k+1}\|_{L^2}^2
%            +C\|\rho^{k+1}\|_{L^\infty}^\frac{1}{2}\|\nabla
%            \overline{u}^k\|_{L^2}
%            \|\nabla
%            \theta^k\|_{H^1}\|\sqrt{\rho^{k+1}}\overline{\theta}^{k+1}\|_{L^2}
%    \end{eqnarray*}
using  similar arguments as in the proof of (\ref{2-b-E3.15}), we
have
    \begin{eqnarray}
      &&\frac{d}{dt}\|\sqrt{e_\theta^k\rho^{k+1}}\overline{\theta}^{k+1}\|_{L^2}^2
      +\kappa^0\|\nabla \overline{\theta}^{k+1}\|_{L^2}^2\nonumber\\
            &\leq&D^k_\epsilon(t)
            \left(\|\overline{\rho}^{k+1}\|_{L^\frac{3}{2}\cap
            L^2}^2+\|\overline{\rho}^{k}\|_{
            L^2}^2+\|\sqrt{\rho^{k+1}}\overline{\theta}^{k+1}\|_{L^2}^2
            \right)+\tilde{C}\|\sqrt{\rho^k}\overline{u}^k\|_{L^2}^2\nonumber\\
            &&+\tilde{C}_\epsilon\|\sqrt{\rho^{k}}\overline{\theta}^{k}\|_{L^2}^2+\frac{\kappa^0\epsilon}{8}\|\nabla \overline{\theta}^{k}\|_{L^2}^2
            +\tilde{C}\|\nabla \overline{u}^k\|_{L^2}^2,
            \ \forall k\geq1,
            \label{2-b-E3.25}
    \end{eqnarray}
where $D^k_\epsilon(t)\in L^1(0,T_5)$ such that
$\int^{T_5}_0D^k_\epsilon(s)ds\leq \tilde{C}+\tilde{C}_\epsilon
t$.

Finally, from (\ref{2-b-E3.12}), using  similar arguments as in
the proof of (\ref{2-b-E3.16}), we easily deduce that
    \begin{eqnarray}
      &&\frac{d}{dt}\|\sqrt{\rho^{k+1}}\overline{u}^{k+1}\|_{L^2}^2
      +\mu^0\|\nabla \overline{u}^{k+1}\|_{L^2}^2\nonumber\\
            &\leq&E^k_\epsilon(t)
            \left(\|\overline{\rho}^{k+1}\|_{L^\frac{3}{2}\cap
            L^2}^2+\|\overline{\rho}^{k}\|_{
            L^2}^2+\|\sqrt{\rho^{k+1}}\overline{u}^{k+1}\|_{L^2}^2
            +\|\sqrt{\rho^{k+1}}\overline{\theta}^{k+1}\|_{L^2}^2
            \right)\nonumber\\
                    &&+\tilde{C}_\epsilon\|\sqrt{\rho^{k}}\overline{\theta}^{k}\|_{L^2}^2
                    +\frac{\epsilon^2}{8}\|\nabla \overline{\theta}^{k}\|_{L^2}^2
            +\epsilon\|\nabla \overline{u}^k\|_{L^2}^2,
            \label{2-b-E3.26}
    \end{eqnarray}
where $E^k_\epsilon(t)\in L^1(0,T_5)$ such that
$\int^{T_5}_0E^k_\epsilon(s)ds\leq \tilde{C}+\tilde{C}_\epsilon t$
for $k\geq1$.

From (\ref{2-b-E3.23})-(\ref{2-b-E3.26}), using same arguments
 in the proof of (\ref{2-b-E3.22}), we can show that there exists a small time
  $T_*\in(0,T_5]$,
  such that the
sequence $(\rho^k,\theta^k,u^k)$ also converges to a limit
$(\rho,\theta,u)$ in the sense of
    $$
      \left\{
      \begin{array}{lll}
      \rho^k-\rho^1\rightarrow \rho-\rho^1&in&L^\infty(0,T_*;L^2\cap L^\frac{3}{2})\\
      ({\theta}^{k},u^k)\rightarrow(\theta,u)&in&L^2(0,T_*;D^1_0).
      \end{array}
      \right.
    $$
 Then adapting
previous arguments for the case $\rho^\infty>0$, we can easily
prove that $(\rho,\theta,u)$ is a unique solution to the problem
(\ref{2-b-E3.1})-(\ref{2-b-E3.6}) with the regularities
(\ref{2-b-E3.8-1})-(\ref{2-b-E3.8-4}).

It remains to prove the blow-up criterion (\ref{2-b-E3.13-1}). To
prove this, suppose that $T^*<\infty$, and let us introduce define
by
    \begin{eqnarray*}
    \Phi(t)&=&1+\|\sqrt{\rho}\theta_t(t,\cdot)\|_{L^2}+
    \int^t_0\|\nabla \theta_t(s,\cdot)\|_{L^2}^2ds\\
    &&+\|\rho(t,\cdot)-\rho^\infty(\cdot)
    \|_{L^6\cap D^1\cap D^{1,q}}+\|\theta(t,\cdot)\|_{D^1_0\cap D^{1,6}}+\|u(t,\cdot)\|_{
      D^{1}_0},
    \end{eqnarray*}
        \begin{eqnarray*}
        \Psi(t)&=&\|\rho_t(t,\cdot)\|_{L^2\cap L^q}+\|\rho^{\frac{1}{2}}u_t(t,\cdot)\|_{L^2}
      +\|(\theta,u)(t,\cdot)\|_{ D^2}\\
            &&+\int^{t}_0\|u_t(\tau,\cdot)\|_{D_0^1}^2+\|(\theta,u)(\tau,\cdot)\|_{D^{2,q}}^2
               d\tau,
        \end{eqnarray*}
and    $
    J(t)=\Phi(t)+\Psi(t),
    $
for $0<t<T^*$. Let $\tau$ be a fixed time in $(0,T^*)$. Then
$(\rho,\theta,u)$ is a strong solution to the equations
(\ref{2-b-E3.1})-(\ref{2-b-E3.6}) in $[\tau, T^*)\times\Omega$,
which satisfy (\ref{2-b-E3.8-1})-(\ref{2-b-E3.8-3}).

If (\ref{2-b-E3.13-1}) is not hold. We can find a constant $K$
such that
    \begin{equation}
      J(\tau)+\sup_{0\leq t\leq T^*}\Phi(t)\leq K.
    \end{equation}
Hence, following  similar arguments as in Section \ref{2-b-Sec2}
with $T=T^*$, we can prove the analogues of
 (\ref{2-b-E2.37-1}), (\ref{2-b-E2.48}) and
(\ref{2-b-E2.39}): for each $t\in(\tau, T^*)$,
    \begin{eqnarray}
      &&\|(\sqrt{\rho}u_t)(t,\cdot)\|_{L^2}^2+\int^t_0\|u_t(s,\cdot)\|_{D^1_0}^2ds\nonumber\\
      &\leq& C_K(1+T^*)+\int^t_0\left[C_K(1+\|\nabla \theta_t\|_{L^2}^2)(
      \|u\|_{D^2}^2+\|\sqrt{\rho}u_t\|_{L^2}^2)+\|u\|_{D^{2,q}}^2\right]ds,
      \label{2-b-E3.33}
    \end{eqnarray}
        \begin{equation}
          \|u(t,\cdot)\|_{D^2}\leq
          C_K(1+\|\sqrt{\rho}u_t\|_{L^2}),
          \label{2-b-E3.33-1}
        \end{equation}
and
        \begin{equation}
          \|u(t,\cdot)\|_{D^{2,q}}^2\leq
          C_K(\epsilon)(1+\|f\|_{L^q}^2+\|\sqrt{\rho}u_t\|_{L^2})
          +\epsilon\|\nabla u_t\|_{L^2}^2.\label{2-b-E3.38}
        \end{equation}
 Combining (\ref{2-b-E3.33})-(\ref{2-b-E3.38}), we have
    \begin{eqnarray*}
    &&\|(\sqrt{\rho}u_t)(t,\cdot)\|_{L^2}^2+\int^t_0\|u_t(s,\cdot)\|_{D^1_0}^2ds\\
    &\leq& C_K(1+T^*)+C_K\int^t_0(1+\|\nabla \theta_t\|_{L^2}^2)
     \|\sqrt{\rho}u_t\|_{L^2}^2ds,
      \ \forall\ t\in(\tau, T^*),
    \end{eqnarray*}
and using the Gronwall's inequality, we obtain that for each
$t\in(\tau,T^*)$,
    $$
    \|(\sqrt{\rho}u_t)(t,\cdot)\|_{L^2}^2+\int^t_0\|u_t(s,\cdot)\|_{D^1_0}^2ds
    \leq C_K(1+T^*).
    $$
From this estimate and (\ref{2-b-E3.33-1})-(\ref{2-b-E3.38}), we
have
    $$
    \|u(t,\cdot)\|_{D^1}+\int^t_0\|u(s,\cdot)\|_{D^{2,q}}^2ds\leq
    C_K(1+T^*),\ \forall\ t\in(\tau, T^*).
    $$
And using the estimates (\ref{2-b-E2.25.1}),
(\ref{2-b-E2.33-1})-(\ref{2-b-E2.33-2-1}), we obtain
    $$
    \|\rho_t(t,\cdot)\|_{L^2\cap L^q}\leq C_K,
    $$
    $$
    \|\theta(t,\cdot)\|_{ D^2}+\int^{t}_0\|\theta(s,\cdot)\|_{D^{2,q}}^2
    ds\leq C_K(1+T^*),
    $$
and
    \begin{equation}
      J(t)\leq C_K(1+T^*),
    \end{equation}
for all $t\in(\tau,T^*)$. It is a contradiction, because the
maximality of $T^*$ implies that $J(t)\rightarrow\infty$ as
$t\rightarrow T^*$. So (\ref{2-b-E3.13-1}) must be hold. This
completes the proof of Theorem \ref{2-b-T1}.

\begin{rem}
By  similar arguments in Lemma \ref{2-b-L2.4}, using the cut-off
function, we could prove that $\bar{\rho}^k$, $ \bar{\rho}\in
L^\infty(0, T_*;L^2)$ or $\bar{\rho}^k$, $ \bar{\rho}\in
L^\infty(0, T_*;L^{\frac{3}{2}}\cap L^2)$, and
$\sqrt{\rho^k}\bar{\theta}^k$, $ \sqrt{\rho^k}\bar{u}^k$,
$\sqrt{\rho_1}\bar{\theta}$, $\sqrt{\rho_1}\bar{u}\in L^\infty(0,
T_*;L^2)$   rigorously.
\end{rem}

\section{Some regularity results on elliptic system}\label{2-b-Sec4}
\setcounter{equation}{0} In this section, we derive some
regularity estimates for the so-called Lam\'{e} system:
    \begin{equation}
     P(\mu,\lambda,u)= -\mathrm{div}(\mu(\nabla u+\nabla u^\top))
      -\nabla(\lambda \mathrm{div}u)=F,
     \label{2-b-E5.1}
    \end{equation}
and the elliptic system:
    \begin{equation}
      -\mathrm{div}(\kappa\nabla e)=G,
      \label{2-b-E5.2}
    \end{equation}
 where $\Omega$ is a bounded or unbounded
domain in $\mathbb{R}^3$, $(\kappa,\mu,\lambda)(x)$ satisfy
    \begin{equation}
    \mu(x)\geq\mu^0>0,\ \ 3\lambda(x)+2\mu(x)\geq  0,
    \ \ \kappa(x)\geq \kappa^0>0,\ \forall x\in\Omega.
    \end{equation}
  The result  in Theorem \ref{2-b-T1} rely crucially on these
estimates.

\begin{lem}\label{2-b-L5.1}
  Assume that $\Omega$ is a bounded domain in $\mathbb{R}^3$ with
  smooth boundary, and let $u\in D^1_0(\Omega)\cap D^{1,q}(\Omega)$ be a weak
  solution to the system (\ref{2-b-E5.1}), where $1<q<\infty$. If
  $F\in L^q(\Omega)$, $(\mu,\lambda)\in C^\beta(\Omega)$ with  $\beta\in(0,1)$, and $(|\nabla \mu|+|\nabla \lambda|)
      |\nabla u|\in L^q(\Omega)$, then $u\in D^{2,q}(\Omega)$ and
    \begin{equation}
      \|u\|_{D^{2,q}(\Omega)}\leq C\left(\|F\|_{L^{q}(\Omega)}
      +\|\nabla u\|_{L^q(\Omega)}+\|(|\nabla \mu|+|\nabla \lambda|)
      |\nabla u|\|_{L^q(\Omega)}
      \right),
    \end{equation}
where $C=C(q,\mu^0,\Omega,\|(\mu,\lambda)\|_{C^\beta(\Omega)})$ is
a positive constant.
\end{lem}
\begin{proof}
  First, we consider the domain $\Omega_R=\Omega\cap
  B_R(x_0)\subset\subset\Omega$. We choose a function $\chi_0\in C^\infty_0(\Omega_R)$
  such that $\chi_0|_{\Omega_\frac{R}{2}}\equiv1$ and $0\leq \chi_0\leq 1$.
  Letting $v_0=\chi_0 u$, from (\ref{2-b-E5.1}), we have
    \begin{equation}
  -\mu(x_0)\Delta
  v_0-(\lambda(x_0)+\mu(x_0))\nabla\mathrm{div}v_0
  =F_1+F_2,
    \end{equation}
where
    \begin{eqnarray}
      F_1&=&\chi_0\left[F+\nabla \mu\cdot(\nabla u+\nabla
      u^\top)+\nabla\lambda\mathrm{div}u
      \right]-\mu(2(\nabla\chi_0\cdot \nabla) u+\Delta \chi_0u)\nonumber\\
            &&
            -(\lambda+\mu)(\nabla\chi_0\mathrm{div}u+\nabla
            u\cdot\nabla\chi_0
            +u\cdot\nabla\nabla\chi_0)
    \end{eqnarray}
and
    \begin{equation}
      F_2=(\mu(x)-\mu(x_0))\Delta v_0+(\lambda(x)+\mu(x)-\lambda(x_0)-\mu(x_0
  ))\nabla \mathrm{div}v_0.
    \end{equation}
From a well-known elliptic theory due to S. Agmon, A. Douglis and
L. Nirenberg in \cite{Agmon64}, we have
    \begin{equation}
    \|v_0\|_{D^{2,q}(\Omega)}\leq C_1\left(
    \|F_1\|_{L^q(\Omega)}+\|F_2\|_{L^q(\Omega)}
    \right)
    \end{equation}
where
$C_1=C_1(q,\mu^0,\Omega,\|(\mu,\lambda)\|_{C^\beta(\Omega)})$.
Since $\mu,\lambda\in C^\beta$,  We have
    $$
    \|F_2\|_{L^q(\Omega)}\leq C_2\sup_{x\in B_R}|x-x_0|^\beta\|v_0\|_{D^{2,q}
    (\Omega)},
    $$
and
    \begin{equation}
    \|v_0\|_{D^{2,q}(\Omega)}\leq C_1\left(
    \|F_1\|_{L^q(\Omega)}+C_2\sup_{x\in B_R}|x-x_0|^\beta\|v_0\|_{D^{2,q}(\Omega)},
    \right)
    \end{equation}
where $C_2=C_2(\|(\mu,\lambda)\|_{C^\beta(\Omega)})$. We can
choose $R_0>0$ satisfying $C_1C_2R_0^\beta<\frac{1}{2}$, and
obtain
    \begin{eqnarray}
    \|u\|_{D^{2,q}(\Omega_\frac{R}{2})}&\leq&
    \|v_0\|_{D^{2,q}(\Omega)}\leq 2C_1 \|F_1\|_{L^q}\nonumber\\
            &\leq& C_3\left(\|F\|_{L^q(\Omega)}
            +\|\nabla u\|_{L^q(\Omega)}+\|(|\nabla \mu|+|\nabla \lambda|)
      |\nabla u|\|_{L^q(\Omega)}
            \right)    \label{2-b-E4.8}
    \end{eqnarray}
where
$C_3=C_3(q,\mu^0,\Omega,\|(\mu,\lambda)\|_{C^\beta(\Omega)})$ and
for all $0<R\leq R_0$. Here, we use the Poincar\'{e}'s inequality
    $$
    \|u\|_{L^q(\Omega)}\leq C\|\nabla u\|_{L^q(\Omega)},
    $$
since $u|_{\partial \Omega}=0$. Second, we assume the
$\Omega_1=B_1(0)\cap \mathbb{R}^3_+$ and consider the domain
$\Omega_r=B_r\cap \mathbb{R}^3_+$. We choose a function $\chi_1\in
C^\infty_0(B_r)$
  such that $\chi_1|_{\Omega_\frac{r}{2}}\equiv1$ and $0\leq \chi_1\leq 1$.
  Letting $v_1=\chi_1 u$, using  similar arguments in the proof of
  (\ref{2-b-E4.8}), we have that there exists a constant $r_0>0$,
  and
    \begin{equation}
    \|u\|_{D^{2,q}(\Omega_\frac{r}{2})}\leq C_4\left(\|F\|_{L^q(\Omega)}
            +\|\nabla u\|_{L^q(\Omega)}+\|(|\nabla \mu|+|\nabla \lambda|)
      |\nabla u|\|_{L^q(\Omega)}
            \right)    \label{2-b-E4.9}
    \end{equation}
where
$C_4=C_4(q,\mu^0,\Omega_1,\|(\mu,\lambda)\|_{C^\beta(\Omega)})$
and for all $0<r\leq r_0$.

For the general bounded domain $\Omega$, since $\partial \Omega
\in C^\infty$, from (\ref{2-b-E4.9}), we easily obtain that for
each $x_0\in
\partial\Omega$, there exists $r_{x_0}>0$ such that
    \begin{equation}
    \|u\|_{D^{2,q}(B_\frac{r}{2}(x_0)\cap \Omega)}\leq C_5\left(\|F\|_{L^q(\Omega)}
            +\|\nabla u\|_{L^q(\Omega)}+\|(|\nabla \mu|+|\nabla \lambda|)
      |\nabla u|\|_{L^q(\Omega)}
            \right)    \label{2-b-E4.10}
    \end{equation}
where
$C_5=C_5(q,\mu^0,\Omega,\|(\mu,\lambda)\|_{C^\beta(\Omega)})$ and
for all $0<r\leq r_{x_0}$. Since $\Omega$ is a bounded domain,
using the Finite Covering Theorem, we have
    $$
    \|u\|_{D^{2,q}(\Omega)}\leq C_6\left(\|F\|_{L^q(\Omega)}
            +\|\nabla u\|_{L^q(\Omega)}+\|(|\nabla \mu|+|\nabla \lambda|)
      |\nabla u|\|_{L^q(\Omega)}
            \right)
    $$
where
$C_6=C_6(q,\mu^0,\Omega,\|(\mu,\lambda)\|_{C^\beta(\Omega)})$.
\end{proof}

\begin{lem}\label{2-b-L5.2}
  Let $\Omega_R=B_R$ or $B_R\cap \mathbb{R}^3_+$ with $R\geq1$, and
   $u\in D^1_0(\Omega_R)\cap D^{1,q}(\Omega_R)$ be a weak
  solution to the system (\ref{2-b-E5.1}), where $1<q<\infty$. If
  $F\in L^q(\Omega_R)$, $(\mu,\lambda)\in C^\beta
  (\Omega_R)$ with  $\beta\in(0,1)$, and $(|\nabla \mu|+|\nabla
\lambda|)
      |\nabla u|\in L^q(\Omega_R)$, then $u\in D^{2,q}(\Omega_R)$ and
    \begin{equation}
      \|u\|_{D^{2,q}(\Omega_R)}\leq C\left(\|F\|_{L^{q}(\Omega_R)}
      +\|\nabla u\|_{L^q(\Omega_R)}+\|(|\nabla \mu|+|\nabla \lambda|)
      |\nabla u|\|_{L^q(\Omega_R)}
      \right),
    \end{equation}
where
$C=C(q,\mu^0,\Omega_1,\|(\mu,\lambda)\|_{C^\beta(\Omega_R)})$ is a
positive constant, independent of $R$.
\end{lem}
\begin{proof}
  If we define $p=\frac{3}{1-\beta}>3$, $v\in W^{1,q}_0(\Omega_1)$ by $v(x)=u(Rx)$,
  $\mu_R(x)=R^{\frac{3}{p}-1}\mu(Rx)$ and  $\lambda_R(x)=R^{\frac{3}{p}-1}\lambda(Rx)$ for $x\in
  \Omega_1$, then we have:
    $$
    P(\mu_R,\lambda_R,v)=R^{\frac{3}{p}+1}F(Rx)=F_R
    \textrm{ for }x\in \Omega_1.
    $$
  Hence, it follows from Lemma \ref{2-b-L5.1} and the fact
$\|(\mu_R,\lambda_R)\|_{C^\beta(\Omega_1)}\leq
\|(\mu,\lambda)\|_{C^\beta(\Omega_R)}$ that
    $$
    \|v\|_{D^{2,q}(\Omega_1)}\leq C_7\left(\|F_R\|_{L^{q}(\Omega_1)}
      +\|\nabla v\|_{L^q(\Omega_1)}+\|(|\nabla \mu_R|+|\nabla \lambda_R|)
      |\nabla v|\|_{L^q(\Omega_1)}
      \right),
    $$
where the positive constant
$C_7=C_7(q,\mu^0,\Omega_1,\|(\mu,\lambda)\|_{C^\beta(\Omega_R)})$
is independent of $R$. Converting this back into the unscaled
variables, we obtain
    \begin{eqnarray*}
    \|u\|_{D^{2,q}(\Omega_R)}&\leq& C_7\left(R^{\frac{3}{p}-1}\|F\|_{L^{q}(\Omega_R)}
      +R^{-1}\|\nabla u\|_{L^q(\Omega_R)}+R^{\frac{3}{p}-1}\|(|\nabla \mu|+|\nabla \lambda|)
      |\nabla u|\|_{L^q(\Omega_R)}
      \right)\\
            &\leq&C_7\left(\|F\|_{L^{q}(\Omega_R)}
      +\|\nabla u\|_{L^q(\Omega_R)}+\|(|\nabla \mu|+|\nabla \lambda|)
      |\nabla u|\|_{L^q(\Omega_R)}
      \right).
    \end{eqnarray*}.
\end{proof}

\begin{lem}\label{2-b-L5.3}
  Let $\Omega_R=B_R\cap \Omega$, where $\Omega$ is an exterior
  domain in $\mathbb{R}^3$ with smooth boundary. Choose a fixed number $R_1\geq1$ such
  that $\Omega^c\subset B_{R_1}$ and assume $R>2R_1$. Assume that
   $u\in D^1_0(\Omega_R)\cap D^{1,q}(\Omega_R)$ be a weak
  solution to the system (\ref{2-b-E5.1}), where $1<q<\infty$. If
  $F\in L^q(\Omega_R)$, $(\mu,\lambda)\in C^\beta
 (\Omega_R)$ with  $\beta\in(0,1)$, and $(|\nabla \mu|+|\nabla \lambda|)
      |\nabla u|\in L^q(\Omega_R)$, then $u\in D^{2,q}(\Omega_R)$ and
    \begin{eqnarray}
      \|u\|_{D^{2,q}(\Omega_R)}&\leq& C\left(\|F\|_{L^{q}(\Omega_R)}
      +\|\nabla u\|_{L^q(\Omega_R)}+\|(|\nabla \mu|+|\nabla \lambda|)
      |\nabla u|\|_{L^q(\Omega_R)}\right.\nonumber\\
            &&\left.+\|(|\nabla \mu|+|\nabla \lambda|)
      |u|\|_{L^q(\Omega_{2R_1})}
      \right),
    \end{eqnarray}
where the positive constant
$C=C(q,\mu^0,R_1,\|(\mu,\lambda)\|_{C^\beta(\Omega_R)})$ is
independent of $R$.
\end{lem}
\begin{proof}
Choosing a cut-off function $\phi\in C^\infty_c(\mathbb{R}^3)$ such
that $\phi\equiv1$ in $B_{R_1}$ and $\phi\equiv0$ in $B_{2R_1}^c$,
we define
    $$
    v=\phi u
    \textrm{ and }
    w=(1-\phi)u:=\psi u.
    $$
First, observing that $v\in W^{1,q}_0(\Omega_{2R_1})$, we deduce
from Lemma \ref{2-b-L5.1} and Poincar\'{e}'s inequality
    $$\|u\|_{L^q(\Omega_{2R_1})}\leq C\|\nabla u\|_{L^q(\Omega_{2R_1})},$$
(since $u|_{\partial \Omega}=0$),    that
    \begin{eqnarray}
      &&\|v\|_{D^{2,q}({\Omega_{2R_1}})}\nonumber\\
      &\leq&C\left( \|P(\mu,\lambda,v)\|_{L^{q}(\Omega_{2R_1})}
      +\|\nabla v\|_{L^q(\Omega_{2R_1})}+\|(|\nabla \mu|+|\nabla \lambda|)
      |\nabla v|\|_{L^q(\Omega_{2R_1})}
      \right)\nonumber\\
            &\leq&C\left( \|F\|_{L^{q}(\Omega_{2R_1})}
            +\|\nabla u\|_{L^q(\Omega_{2R_1})}+\|(|\nabla \mu|+|\nabla \lambda|)
            (|\nabla u|+|u|)\|_{L^q(\Omega_{2R_1})}      \right).
            \label{2-b-E4.13}
    \end{eqnarray}
Next, to estimate $w$, we observe that $w=0$ in $B_{R_1}$ and
$w\in D^1_0\cap D^{1,q}(B_R)$. Then it follows from Lemma
\ref{2-b-L5.2} and Poincar\'{e}'s inequality that
    \begin{eqnarray}
      &&\|w\|_{D^{2,q}({\Omega_{R}})}\nonumber\\
      &\leq&C\left( \|P(\mu,\lambda,w)\|_{L^{q}(\Omega_{R})}
      +\|\nabla w\|_{L^q(\Omega_{R})}+\|(|\nabla \mu|+|\nabla \lambda|)
      |\nabla w|\|_{L^q(\Omega_{R})}
      \right)\nonumber\\
            &\leq&C\left( \|F\|_{L^{q}(\Omega_{R})}
            +\|\nabla u\|_{L^q(\Omega_{R})}+\|(|\nabla \mu|+|\nabla \lambda|)
            |\nabla u|\|_{L^q(\Omega_{R})}\right.\nonumber\\
    &&+\left.   \|(|\nabla \mu|+|\nabla \lambda|)
            | u|\|_{L^q(\Omega_{2R_1})}   \right).
            \label{2-b-E4.14}
    \end{eqnarray}
Combining (\ref{2-b-E4.13}) and (\ref{2-b-E4.14}), we finish the
proof.
\end{proof}

Since the constant $C$ in Lemma \ref{2-b-L5.2}-\ref{2-b-L5.3} is
independent of $R$, using the domain expansion technique as in
 \cite{Lions96} and the fact that the solution to the
system (\ref{2-b-E5.1}) is unique in $D^1_0(\Omega)$,  we easily
prove
\begin{lem}\label{2-b-L5.4}
  Let $\Omega$ be the whole space $\mathbb{R}^3$, the half space $\mathbb{R}^3_+$
  or an exterior domain in $\mathbb{R}^3$ with smooth boundary.
  Assume that   $u\in D^1_0(\Omega)\cap D^{1,q}(\Omega)$ be a weak
  solution to the system (\ref{2-b-E5.1}), where $1<q<\infty$. If
  $F\in L^q(\Omega)$, $(\mu,\lambda)\in C^\beta
  (\Omega)$ with  $\beta\in(0,1)$, and $(|\nabla \mu|+|\nabla \lambda|)
      |\nabla u|\in L^q(\Omega)$, then $u\in D^{2,q}(\Omega)$
      satisfying
    \begin{equation}
      \|u\|_{D^{2,q}(\Omega)}\leq C\left(\|F\|_{L^{q}(\Omega)}
      +\|\nabla u\|_{L^q(\Omega)}+\|(|\nabla \mu|+|\nabla \lambda|)
      |\nabla u|\|_{L^q(\Omega)}
      \right),
    \end{equation}
when $\Omega$ is the whole space $\mathbb{R}^3$ or the half space
$\mathbb{R}^3_+$; and
    \begin{eqnarray}
      \|u\|_{D^{2,q}(\Omega)}&\leq& C\left(\|F\|_{L^{q}(\Omega)}
      +\|\nabla u\|_{L^q(\Omega)}+\|(|\nabla \mu|+|\nabla \lambda|)
      |\nabla u|\|_{L^q(\Omega)}\right.\nonumber\\
            &&\left.+\|(|\nabla \mu|+|\nabla \lambda|)
      |u|\|_{L^q(\Omega_{2R_1})}
      \right),
    \end{eqnarray}
when $\Omega$ is an exterior domain in $\mathbb{R}^3$ with smooth
boundary, and $\Omega^c\subset B_{R_1}$ with  $R_1>1$, where
$C=C(q,\mu^0,\|(\mu,\lambda)\|_{C^\beta(\Omega)})$.
\end{lem}
\begin{lem}\label{2-b-L5.5}
  Let $\Omega$ be the whole space $\mathbb{R}^3$, the half space
  $\mathbb{R}^3_+$,
   an exterior domain in $\mathbb{R}^3$ or
    a bounded domain in $\mathbb{R}^3$ with smooth boundary.
  Assume
  $F\in L^2(\Omega)$, $(\mu,\lambda)\in C^\beta\cap
D^{1,p}(\Omega)$ with  $p>3$, $\beta\in(0,1)$,  and   $u\in
D^1_0(\Omega)\cap D^{2}(\Omega)$ be a weak
  solution to the system (\ref{2-b-E5.1}), then we have
    \begin{equation}
      \|u\|_{D^{2}(\Omega)}\leq C\left(\|F\|_{L^{2}(\Omega)}
      +\|\nabla u\|_{L^2(\Omega)}+\||\nabla \mu|+|\nabla \lambda|
      \|_{L^p(\Omega)}^\frac{p}{p-3}\| u\|_{D^1_0(\Omega)}
      \right),
    \end{equation}
where $C=C(p,\mu^0,\|(\mu,\lambda)\|_{C^\beta(\Omega)})$ is
positive constant.
\end{lem}
\begin{proof}
First, we consider the case that $\Omega$ is an exterior domain in
$\mathbb{R}^3$ with smooth boundary. Using H\"{o}lder,
Galiardo-Nirenberg
  and Young's inequalities,  we obtain
    \begin{eqnarray}
      &&C\|(|\nabla \mu|+|\nabla \lambda|)
      |\nabla u|\|_{L^2(\Omega)}\nonumber\\
            &\leq& C\||\nabla \mu|+|\nabla
              \lambda|\|_{L^p(\Omega)}
              \|\nabla u\|_{L^\frac{2p}{p-2}(\Omega)}\nonumber\\
      &\leq& C\||\nabla \mu|+|\nabla
      \lambda|\|_{L^p(\Omega)}
      \| u\|_{D^1_0(\Omega)}^\frac{p-3}{p}
      \|u\|_{D^2(\Omega)}^\frac{3}{p}\nonumber\\
            &\leq& C\||\nabla \mu|+|\nabla
      \lambda|\|_{L^p(\Omega)}^\frac{p}{p-3}
      \| u\|_{D^1_0(\Omega)}+\frac{1}{4}
      \|u\|_{D^2(\Omega)},
      \label{2-b-E4.18}
    \end{eqnarray}
and
    \begin{eqnarray}
      &&\|(|\nabla \mu|+|\nabla \lambda|)
      |u|\|_{L^2(\Omega_{2R_1})}\nonumber\\
            &\leq& C\||\nabla \mu|+|\nabla
              \lambda|\|_{L^p(\Omega_{2R_1})}
              \| u\|_{L^\frac{2p}{p-2}(\Omega_{2R_1})}\nonumber\\
                        &\leq& C\||\nabla \mu|+|\nabla
                         \lambda|\|_{L^p(\Omega_{2R_1})}
                        \| u\|_{L^6(\Omega_{2R_1})}\nonumber\\
      &\leq& C\||\nabla \mu|+|\nabla
      \lambda|\|_{L^p(\Omega)}
      \| u\|_{D^1_0(\Omega)}.
                  \label{2-b-E4.19}
    \end{eqnarray}
From Lemma \ref{2-b-L5.4} and the estimates
(\ref{2-b-E4.18})-(\ref{2-b-E4.19}), we have
    \begin{equation}
      \|u\|_{D^{2}(\Omega)}\leq C\left(\|F\|_{L^{2}(\Omega)}
      +\|\nabla u\|_{L^2(\Omega)}+\||\nabla \mu|+|\nabla \lambda|
      \|_{L^p(\Omega)}^\frac{p}{p-3}\| u\|_{D^1_0(\Omega)}
      \right).
    \end{equation}
By  similar arguments, we can easily obtain the result in the case
that $\Omega$ is the whole space $\mathbb{R}^3$, the half space
  $\mathbb{R}^3_+$ or
    a bounded domain in $\mathbb{R}^3$ with smooth boundary.
\end{proof}
\begin{lem}\label{2-b-L5.6}
  Let $\Omega$ be the whole space $\mathbb{R}^3$, the half space
  $\mathbb{R}^3_+$,
   an exterior domain in $\mathbb{R}^3$ or
    a bounded domain in $\mathbb{R}^3$ with smooth boundary.
  Assume
  $F\in L^p(\Omega)$, $(\mu,\lambda)\in C^\beta\cap
D^{1,p}(\Omega)$
    and   $u\in D^1_0(\Omega)\cap D^{1,p}(\Omega)\cap D^{2.p}(\Omega)$ be a weak
  solution to the system (\ref{2-b-E5.1}), where $p>3$ and $\beta\in(0,1)$, then we have
    \begin{equation}
      \|u\|_{D^{2,p}(\Omega)}\leq C\left(\|F\|_{L^{p}(\Omega)}
      +\|\nabla u\|_{L^p(\Omega)}+\||\nabla \mu|+|\nabla \lambda|
      \|_{L^p(\Omega)}^{\frac{5p-6}{2p-6}}\| u\|_{D^{1}_0(\Omega)}
      \right),
    \end{equation}
where $C=C(p,\mu^0,\|(\mu,\lambda)\|_{C^\beta(\Omega)})$ is a
positive constant.
\end{lem}
\begin{proof}
First, we consider the case that $\Omega$ is an exterior domain in
$\mathbb{R}^3$ with smooth boundary. Using H\"{o}lder,
Galiardo-Nirenberg
  and Young's inequalities,  we obtain
    \begin{eqnarray}
      &&C\|(|\nabla \mu|+|\nabla \lambda|)
      |\nabla u|\|_{L^p(\Omega)}\nonumber\\
            &\leq& C\||\nabla \mu|+|\nabla
              \lambda|\|_{L^p(\Omega)}
              \|\nabla u\|_{L^\infty(\Omega)}\nonumber\\
      &\leq& C\||\nabla \mu|+|\nabla
      \lambda|\|_{L^p(\Omega)}
      \| u\|_{D^1_0(\Omega)}^\frac{2p-6}{5p-6}
      \|u\|_{D^{2,p}(\Omega)}^\frac{3p}{5p-6}\nonumber\\
            &\leq& C\||\nabla \mu|+|\nabla
      \lambda|\|_{L^p(\Omega)}^{\frac{5p-6}{2p-6}}
      \| u\|_{D^1_0(\Omega)}+\frac{1}{4}
      \|u\|_{D^{2,p}(\Omega)},
      \label{2-b-E4.18-1}
    \end{eqnarray}
and
    \begin{eqnarray}
      &&C\|(|\nabla \mu|+|\nabla \lambda|)
      |u|\|_{L^p(\Omega_{2R_1})}\nonumber\\
            &\leq& C\||\nabla \mu|+|\nabla
              \lambda|\|_{L^p(\Omega_{2R_1})}
              \| u\|_{L^\infty(\Omega_{2R_1})}\nonumber\\
      &\leq& C\||\nabla \mu|+|\nabla
      \lambda|\|_{L^p(\Omega)}
      \| u\|_{D^1_0(\Omega)}^\frac{4p-6}{5p-6}
      \|u\|_{D^{2,p}(\Omega)}^\frac{p}{5p-6}\nonumber\\
            &\leq& C\||\nabla \mu|+|\nabla
      \lambda|\|_{L^p(\Omega)}^{\frac{5p-6}{4p-6}}
      \| u\|_{D^1_0(\Omega)}+\frac{1}{4}
      \|u\|_{D^{2,p}(\Omega)}.
                  \label{2-b-E4.19-1}
    \end{eqnarray}
From Lemma \ref{2-b-L5.4} and the estimates
(\ref{2-b-E4.18-1})-(\ref{2-b-E4.19-1}), we have
    \begin{equation}
      \|u\|_{D^{2}(\Omega)}\leq C\left(\|F\|_{L^{2}(\Omega)}
      +\|\nabla u\|_{L^2(\Omega)}+\||\nabla \mu|+|\nabla \lambda|
      \|_{L^p(\Omega)}^{\frac{5p-6}{2p-6}}\| u\|_{D^1_0(\Omega)}
      \right).
    \end{equation}
By  similar arguments, we can easily obtain the result in the case
that $\Omega$ is the whole space $\mathbb{R}^3$, the half space
  $\mathbb{R}^3_+$ or
    a bounded domain in $\mathbb{R}^3$ with smooth boundary.
\end{proof}
Using  similar arguments in Lemmas \ref{2-b-L5.5}-\ref{2-b-L5.6},
we can easily obtain:
\begin{lem}\label{2-b-L5.7}
  Let $\Omega$ be the whole space $\mathbb{R}^3$, the half space
  $\mathbb{R}^3_+$,
   an exterior domain in $\mathbb{R}^3$ or
    a bounded domain in $\mathbb{R}^3$ with smooth boundary.
  Assume
  $G\in L^2(\Omega)$, $\kappa\in C^\beta\cap
D^{1,p}(\Omega)$ with  $p>3$, $\beta\in(0,1)$,  and   $e\in
D^{1}_0(\Omega)\cap D^{2}(\Omega)$ be a weak
  solution to the system (\ref{2-b-E5.2}), then we have
    \begin{equation}
      \|e\|_{D^{2}(\Omega)}\leq C\left(\|G\|_{L^{2}(\Omega)}
      +\|\nabla e\|_{L^2(\Omega)}+\|\nabla \kappa
      \|_{L^p(\Omega)}^\frac{p}{p-3}\| e\|_{D^1_0(\Omega)}
      \right),
    \end{equation}
where $C=C(p,\kappa^0,\|\kappa\|_{C^\beta(\Omega)})$ is a positive
constant.
\end{lem}
\begin{lem}\label{2-b-L5.8}
  Let $\Omega$ be the whole space $\mathbb{R}^3$, the half space
  $\mathbb{R}^3_+$,
   an exterior domain in $\mathbb{R}^3$ with smooth boundary or
    a bounded domain in $\mathbb{R}^3$ with smooth boundary.
  Assume
  $G\in L^p(\Omega)$, $\kappa\in C^\beta\cap
D^{1,p}(\Omega)$
    and   $e\in D^1_0(\Omega)\cap D^{1,p}(\Omega)\cap D^{2.p}(\Omega)$ be a weak
  solution to the system (\ref{2-b-E5.2}) with  $p>3$ and $\beta\in(0,1)$, then we have
    \begin{equation}
      \|e\|_{D^{2,p}(\Omega)}\leq C\left(\|G\|_{L^{p}(\Omega)}
      +\|\nabla e\|_{L^p(\Omega)}+\|\nabla \kappa
      \|_{L^p(\Omega)}^{\frac{5p-6}{2p-6}}\| e\|_{D^{1}_0(\Omega)}
      \right),
    \end{equation}
where $C=C(p,\kappa^0,\|\kappa\|_{C^\beta(\Omega)})$ is a positive
constant.
\end{lem}
\begin{rem}
  We can remove the assumptions in Lemma
  \ref{2-b-L5.5}-\ref{2-b-L5.8} that $(u,e)\in D^2$ or $D^{2,q}$.
  Indeed, we can use the usual iteration argument
  to prove that $(u,e)\in D^2$ or $D^{2,q}$ rigorously.
\end{rem}

\section{Blow-up}\label{2-c-s1}\setcounter{equation}{0}
In this section, we show the blow-up of the smooth solution to the
compressible Navier-Stokes equations in $\mathbb{R}^n(n\geq1)$
when the initial density has compactly support and the initial
total momentum is nonzero. We denote by $X(\alpha,t)$ the particle
trajectory starting at $\alpha$ when $t=0$, that is
    $$
    \frac{d}{dt}X(\alpha,t)=u(X(\alpha,t),t)
    \textrm{ and }X(\alpha,0)=\alpha.
    $$
We set
        $
    \Omega(0)=\mathrm{supp}\rho_0
        $
and
    $
    \Omega(t)=\{x=X(\alpha,t)|\alpha\in \Omega(0)\}.
    $
From the transport equation (\ref{2-c-E1.1}), one can easily show
that
    $
    \mathrm{supp}\rho(x,t)=\Omega(t)
    $
and hence from the equation of state (\ref{2-c-E1.4-1}) that
    $$
    p(x,t)=\theta(x,t)=e(x,t)=0,
   \ \textrm{ if }x\in \Omega(t)^c.
    $$
Therefore, from the equation (\ref{2-c-E1.2}), we observe that
    $$
    \frac{\mu(0,0)}{2}|\nabla u+\nabla
          u^\top|^2+\lambda(0,0)(\mathrm{div } u)^2=0,
   \ \textrm{ if }x\in \Omega(t)^c.
    $$

The following lemma was shown in \cite{Cho2004,xin}, when
$\mu,\lambda$ are constants. Therefore, the proof is very similar
to that in \cite{Cho2004,xin}, and is just modified to our case
that $(\mu,\lambda)$ satisfies (\ref{2-c-E1.4}), so we omit the
detail.
\begin{lem}\label{2-c-L2.1}
  We assume that $(\mu,\lambda)$ satisfies (\ref{2-c-E1.4}), and
   $(\rho,S,u)\in
C^1([0,T];H^k)$, $k>[\frac{n}{2}]+2$, is a solution to the Cauchy
problem (\ref{2-c-E1.1})-(\ref{2-c-E1.2}).  Then
    $$
    u(x,t)=0
    \textrm{ in } x\in \Omega(t)^c.
    $$
 Moreover, $\Omega(t)=\Omega(0)$ for all $0 \leq t \leq T$.
\end{lem}

\textbf{Proof of Theorem \ref{2-c-T1.1}}

We assume that $\Omega(t)=\mathrm{supp}\rho(\cdot,t)$ is contained
in a ball $B_{r(t)}$. Multiplying $x_i$ to (\ref{2-c-E1.1}) and
integrating it over $\mathbb{R}^n$, we get the identity
    \begin{equation}
      \frac{d}{dt}\int_{\mathbb{R}^n}\rho x_idx
      =\int_{\mathbb{R}^n}\rho u_idx.
      \label{2-c-E3.1}
    \end{equation}
If we integrate  (\ref{2-c-E1.3})  over $\mathbb{R}^n$, then we also
obtain the identity
    \begin{equation}
      \frac{d}{dt}\int_{\mathbb{R}^n}\rho udx
      =0.\label{2-c-E3.2}
    \end{equation}
Integrating (\ref{2-c-E3.1}) and (\ref{2-c-E3.2}) over $[0,t]$,
respectively, we obtain the following identities
    \begin{equation}
      \int_{\mathbb{R}^n}\rho x_idx=
            \int_{\mathbb{R}^n}\rho_0 x_idx
      +\int^t_0      \int_{\mathbb{R}^n}\rho u_idxds,
      \label{2-c-E3.3}
    \end{equation}
            \begin{equation}
              \int_{\mathbb{R}^n}\rho udx=
            \int_{\mathbb{R}^n}\rho_0 u_0dx.
            \label{2-c-E3.4}
            \end{equation}
From (\ref{2-c-E3.3})-(\ref{2-c-E3.4}), we get
    \begin{equation}
      \int_{\mathbb{R}^n}\rho x_idx=
            \int_{\mathbb{R}^n}\rho_0 x_idx
      +t \int_{\mathbb{R}^n}\rho_0 u_{0i}dx=m_{1i}+m_{2i}t.
      \label{2-c-E3.5}
    \end{equation}

On the other hand, since $\Omega(t)\subset B_{r(t)}$, we can
estimate the bound of the L.H.S. of (\ref{2-c-E3.5}) as follows,
    \begin{eqnarray}
      &&\left|\int_{\mathbb{R}^n}\rho x_idx
      \right|=\left|\int_{\Omega(t)}\rho x_idx
      \right|=\left|\int_{|x|\leq r(t)}\rho x_idx
      \right|\nonumber\\
                &\leq& r(t)\int_{\mathbb{R}^n}\rho dx
                  =r(t)\int_{\mathbb{R}^n}\rho_0 dx=r(t)m_0.
    \label{2-c-E3.6}
    \end{eqnarray}
From Lemma \ref{2-c-L2.1}, we have $\Omega(t)=\Omega(0)$ for
$t\in[0,T]$ and hence $r(t)$ can be chosen to be $r_0$. Therefore,
we have
    $$
    m_0r_0\geq m_{1i}+m_{2i}t\geq -m_0r_0,
    \ i=1,\ldots,n.
    $$
This completes the proof of the theorem \ref{2-c-T1.1}.

\end{document}